\documentclass[12pt]{amsart}
\usepackage{amsmath,amsfonts,amssymb,amscd}
\usepackage{latexsym}
\def\dfont{\textit}
\def\~{{\bf --}}

\newcommand{\BR}{{\mathbb R}}
\newcommand{\BQ}{{\mathbb Q}}
\newcommand{\BC}{{\mathbb C}}
\newcommand{\BP}{{\mathbb P}}
\newcommand{\BZ}{{\mathbb Z}}
\newcommand{\BN}{{\mathbb N}}

\newcommand{\cH}{{\mathcal H}}
\newcommand{\cA}{{\mathcal A}}
\newcommand{\cB}{{\mathcal B}}
\newcommand{\ccF}{{\mathfrak F}}
\newcommand{\cD}{{\mathcal D}}
\newcommand{\cL}{{\mathcal L}}
\newcommand{\cF}{{\mathcal F}}
\newcommand{\cP}{{\mathcal P}}
\newcommand{\cX}{{\mathcal X}}
\newcommand{\cY}{{\mathcal Y}}

\newcommand{\Z}{{\mathbb Z}}

\newcommand{\C}{{\mathbb C}}

\def\HH{\mbox{${\mathcal H}$\kern-5.2pt${\mathcal H}$}}

\newtheorem{theorem}{Theorem}[section]
\newtheorem{proposition}[theorem]{Proposition}
\newtheorem{definition}[theorem]{Definition}
\newtheorem{lemma}[theorem]{Lemma}
\newtheorem{corollary}[theorem]{Corollary}

\newtheorem{theorem }{Theorem}[section]
\newtheorem{proposition }[theorem]{Proposition}
\newtheorem{definition }[theorem]{Definition}
\newtheorem{lemma }[theorem]{Lemma}
\newtheorem{corollary }[theorem]{Corollary}
\newtheorem{notation }[theorem]{Notation}
\newtheorem{remark }[theorem]{Remark}
\newtheorem{example }[theorem]{Example}

\newtheorem{ theorem}{Theorem}[section]
\newtheorem{ proposition}[theorem]{Proposition}
\newtheorem{ definition}[theorem]{Definition}
\newtheorem{ lemma}[theorem]{Lemma}
\newtheorem{ corollary}[theorem]{Corollary}
\newtheorem{ notation}[theorem]{Notation}
\newtheorem{ remark}[theorem]{Remark}
\newtheorem{ example}[theorem]{Example}

 \newcommand{\rmk}{{\bf Comment.\ }}

\def\for{\  \hbox{ for } \ }

\def\and{\  \hbox{ and } \ }
\def\and{\  \hbox{ and } \ }

\def\equal{\stackrel{\,\mathbf{def}}{= \kern-3pt =}}

\def\la{\lambda}

\def\ga{\gamma}
\def\ep{\epsilon}

\def\de{\delta}

\def\si{\sigma}

\newcommand{\bS}{{\mathbf S}}
\newcommand{\bH}{{\mathbf H}}
\newcommand{\bF}{{\mathbf F}}
\newcommand{\bE}{{\mathbf E}}

\def\0{\mathbf{0}}

\def\çF{\mathcal{F}}

\def\m{\mathcal{M}}

\def\p{\mathcal{P}}
\def\cP{\mathcal{P}}
\def\a{\mathcal{A}}

\def\e{\mathcal{E}}

\def\b{\mathcal{B}}

\newcommand{\Aut}{\operatorname{Aut}}

\newcommand{\Ind}{\operatorname{Ind}}

\newcommand{\Rad}{\operatorname{Rad}}

\newcommand{\id}{\operatorname{id}}

\newcommand{\sq}{\phantom{1}\hfill$\qed$}
\newcommand{\Rea}{\Re}
\newcommand{\Ima}{\Im}

\newcommand{\modd}{\mbox{\,mod\,}}
\newcommand{\lr}{\langle}
\newcommand{\rr}{\rangle}
\newcommand{\eps}{\varepsilon}
\newcommand{\phk}{\phi^{(k)}}
\newcommand{\psk}{\psi^{(k)}}
\newcommand{\Res}{\mbox{Res}\;}
\newcommand{\sgn}{\mbox{sgn}}

\def\HH{\mathfrak{H}}

\def\HH{\hbox{${\mathcal H}$\kern-5.2pt${\mathcal H}$}}

\font\smm=msbm10 at 12pt 
\def\symbol#1{\hbox{\smm #1}}
\def\lsmash{{\symbol n}}

\def\#{\sharp}

\title [Double Hecke Algebra]
{From Double Hecke Algebra to Fourier Transform}
\author [Ivan Cherednik] {Ivan Cherednik\ $^\dag$} 
\author [Viktor Ostrik] {Viktor Ostrik\ $^\ddag$}
\date{November 2001}

\thanks{$^\dag\ $ Partially supported by NSF grants 
DMS-0200276}
\thanks{$^\ddag\ $ Partially supported by NSF 
grant DMS-0098830}

\address[I. Cherednik]{Department of Mathematics, UNC 
Chapel Hill, NC 27599, USA\\
chered@math.unc.edu}
\address[V. Ostrik]{Department of Mathematics, MIT, 
77 Mass.Ave,
Cambridge, MA 02139, USA\\
ostrik@math.mit.edu}

\begin{document}
\maketitle

The paper is mainly based on the series of lectures on
the one-dimensional double Hecke algebra delivered by
the first author at Harvard University in 2001. 
It also contains the material of other 
talks (MIT, University Paris 6) and new results. 
The most interesting is the classification of 
finite dimensional representations.
  
Concerning the proofs, we followed the principle,  
{\em the more proof the better},
which the first author adopted during his studies in
combinatorics. Quite a few theorems were
proved twice in the lectures (sometimes even three times),
using different tools. This is true as well in the paper.
  
\bigskip
{\bf Methods.} There are deep classical
origins and relations to the theory of special
functions, old and new.
We did not try to reconstruct 
systematically the history of the subject and review the 
connections either in the
lectures or in the paper. There are many
comments but they are fragmentary.

We recommend \cite{M3, O2, C6, C10} to those
who are interested in recent developments.
Concerning Riemann and other great masters of Fourier analysis, 
see e.g., \cite{Ed}. The book \cite{An} is a good introduction to
the theory of $q$-functions. See also \cite{HO1},\cite{He}
about the relations to the spherical functions,
and Helgason's notes about Harish-Chandra, 
the creator of the harmonic analysis
on the symmetric spaces.
 
The focus is on the $q$-Fourier transform and the corresponding
representations of the double affine Hecke algebra.
This transform was introduced about five years ago as a
$q$-generalization of the classical Hankel and
Harish-Chandra transforms with deep relations to the
Macdonald orthogonal polynomials, combinatorics,
$p$-adic representations, Gaussian sums, and Verlinde algebras.

The objective of the lectures was to convince the audience that
the $q$-theory provides unification as well as simplification.
In the paper, we closely follow the notes of the lectures. 
The exposition is based on the recent theory of the so-called
nonsymmetric Macdonald polynomials. See 
\cite{D, O2, M3, KnS, C3, C4}.
Opdam mentions in \cite{O2} that a definition of the nonsymmetric 
polynomials (in the differential setup) was given in Heckman's 
unpublished lectures.  

The case of $A_1$ is considered in the paper, with
a reservation about the last section devoted to relations to
the $p$-adic theory, which are discussed in complete generality.
However the proofs are mainly of a general nature
and can be transfered to the case of arbitrary (reduced)
root systems. The paper is designed to be an
introduction to the  general multidimensional theory.
There are some exceptions, mainly
when the explicit formula for Macdonald's
truncated theta function was used, which has no reasonable
multidimensional counterparts. 
See \cite{C5,C10} for the general theory.

\vfil
{\bf Results.}
Let us mention
the most interesting ones.

a) Detailed proof of the main formula for the $q$-Mellin
transform found by the second author. It helps in
controlling the spaces of analytic functions involved
in the formula and is closely related to the
analytic theory of the shift operator from \cite{C9}.

b) Integrating the Gaussian over $\BR$ with respect to
the $q$-measure, which leads to the Appell functions.
It was performed by P.~Etingof.
The paper is mainly about imaginary integration
and Jackson-type summation. 
Real integration completes the picture. 

c) A new theory of the
$q$-Fourier transform for the scalar product without the
conjugation $q\mapsto q^{-1}.$ This transform
directly generalizes the one from \cite{O2}. In a way, 
Opdam's $w_0$ is replaced by $T_{w_0},$ which is closely
connected with  \cite{O3,O4}.

d) The classification of irreducible finite dimensional
representations of the double  
Hecke algebra, including the Fourier-invariant ones. 
It is closely connected with \cite{C7,C10}.

e) The limit to the $p$-adic Macdonald--Matsumoto theory 
\cite{Mat} of
spherical functions. It is partially based on the lecture
of the first author at MIT (V.~Kac's seminar). 
The transform from c) is used here.
It develops the corresponding results from \cite{C3}.

\vfil
Actually the lectures and the paper are more about the technique
than about the results themselves. The main theorem of
the course is very simple to formulate:
{\em $q$-Fourier transform is self-dual}.  The problem was
to justify, clarify, extend, and apply this fact.

The participants of the lectures helped a lot to achieve
this goal. Our special thanks go to
Pavel Etingof and David Kazhdan who organized the course,
to Alexander Braverman,
Dmitry Nikshych, and Alexander Polyshchuk.

The main part of the paper was prepared when the first
author was appointed by the Clay Mathematics Institute. 

{\small
\tableofcontents
}

\bigskip
\setcounter{equation}{0}
\section{From Euler's integral to Gaussian sum}
Practically all results in the paper are parts
of the following general program: {\em to connect
the Gaussian sums at roots of unity with the Gauss 
integrals in one theory}. The key points of this
program will be explained in the following three sections. 

The technique is actually quite elementary. One needs
to calculate a couple of intermediate $q$-integrals, similar
to those considered by Ramanujan 
(see Askey's paper in \cite{An}). There is nothing here
beyond classical calculus. 
However it is not surpising that
the connection between the Euler-Gauss integral 
and the Gaussian sums was not established
in the $19$-th century.  The $q$-functions and
$q$-integrals reached a proper level of maturity later,
as well as general understanding of the importance of roots
of unity and $p$-adic methods.

\subsection{Euler's integral and Riemann's zeta}
We begin with one of the most famous classical formulas:
\begin{equation}\label{f11}
2\int_0^\infty e^{-x^2}x^{2k}dx =\Gamma (k+\frac{1}{2}),
\end{equation}
where $k$ is a complex number with $\Rea k>-\frac{1}{2}$.
Actually it is the best way to introduce the 
$\Gamma$-function, so it is more of a definition than a formula.
Indeed, it readily results in

{(i)} the functional equation $\Gamma (x+1)=x\Gamma (x),$

{(ii)} the meromorphic continuation of $\Gamma(x)$
to all complex $x,$

{(iii)} the infinite Weierstrass product formula for $\Gamma.$

For $k=0,$ this
formula reduces to the Poisson integral
\begin{equation}\label{f12}
2\int_0^\infty e^{-x^2}dx=\sqrt{\pi}.
\end{equation}
Other remarkable special values are
$$\Gamma (n+1)=n!\hbox{ \ and\ } \Gamma (n+\frac{1}{2})=
\frac{(2n)!}{2^{2n}(n)!}\sqrt{\pi},\ n\in \BZ_+.$$

Formula (\ref{f11}) has tremendous applications in both
mathematics and physics. In the first place, it gave birth to
{analytic number theory}. The
"perturbation" of (\ref{f11})
$$\int_{-\infty}^\infty \frac{|x|^{2k}}{e^{x^2}-1}dx =
\zeta (k+\frac{1}{2})
\Gamma (k+\frac{1}{2}),\Rea k>0, $$
leads to the analytic continuation of Riemann's
zeta function $\zeta(s)$ and the functional equation.

This is due to Euler and Riemann. It is interesting
to mention that the left-hand side
did appear in physics (Landau and Lifshitz)
as a perturbation of (\ref{f11}).

There are three celebrated closely
connected and entirely open problems concerning
the behavior of $\zeta(k+1/2)$ in the critical strip
$\{-1/2<\Rea k <1/2\}:$

{(i)} the Riemann hypothesis: all zeros are $k$-imaginary,

{(ii)} the distribution of the zeros beyond the  $T\log T$--formula,

{(iii)} the growth estimates in this strip as $\Ima k\to \infty.$

Intensive computer experiments confirm (i)
and strongly indicate that the imaginary parts of the zeros
are distributed randomly subject to the approximate
formula $T\log T$ for the number of zeros between $0$ and
$iT.$ The latter means that essentially
there is nothing to expect beyond the Riemann hypothesis.
This problem is relatively recent. 
Note that there are also conjectural formulas for the 
higher correlators (Bogomolny {\em et.al.}), which are also in 
full harmony with the computer simulations (Odlyzko).

For instance, the connection of the zeros of the $\zeta$
on the critical line with the so-called Gramm points is
too good to be true, as well as 
similar sophisticated qualitative "laws". 
Indeed, all of them were rejected by computers.

I think that quite a few specialists
in the 19th century and at the beginning of the 20th century
suspected that the zeros of the zeta function
are far from being regular.
Now we are almost certain that they are
completely irrational, modulo Riemann's hypothesis.

\subsection{Fourier analysis}
A much more "rational" extension of (\ref{f11}) is the {Fourier
analysis}. Instead of perturbing the Gaussian towards the zeta,
we simply multiply it by the Bessel function. The corresponding
integral, the Hankel transform of the Gaussian, can be
calculated without
difficulties. The theory is plane and square. All standard
facts about the usual Fourier transform take place.
The above program becomes: {\it to connect
the Hankel transform and the Fourier transform
on $\BZ_N$ in one theory.} 

A natural step is to
go from
the Hankel transform to its trigonometric
variant, the Harish-Chandra spherical
transform. Unfortunately the latter is worse than
the Hankel transform. Actually the only features which
completely survive in the Harish-Chandra theory 
are as follows:

(i) the analytic description of the Fourier-image
of the space of compactly supported functions, 

(ii) and the Harish-Chandra inversion formula, generalizing
the selfduality of the Hankel transform.

Technically, we replace the measure $|x|^{2k}dx,$
which makes the Bessel functions pairwise orthogonal, by
$|\sinh(x)|^{2k},$ which governs the orthogonality
of the hypergeometric and spherical functions. The first
serious problem is that the integral
$$2\int_0^\infty e^{-x^2}\sinh(x)^{2k}dx$$
becomes transcendental
apart from $k\in \BZ/2.$  When the Gaussian is multiplied
by the generalized spherical or hypergeometric functions
the integral gets even worse.
It can be calculated exactly only when
$k=0$ and in the so-called group case $k=1.$
It is transcendental even for the
classical spherical functions in the orthogonal
and symplectic cases (when $k=1/2,2$).

Unfortunately we need the latter integral and similar
ones, because they are 
nothing but the Harish-Chandra transforms
of the Gaussian, which are important to know at almost all 
levels of the harmonic analysis on the symmetric spaces.

A bypass was suggested in \cite{C2,C5}. We go from
the Bessel function to the next {\em but one} level, the
{\em basic hypergeometric function}, a $q$-generalization
of the classical  hypergeometric function.
The resulting transform is self-dual, and has
all other good properties of the Hankel transform.
The paper contains a complete algebraic theory of this
$q$-transform and elements of the analytic theory.

To explain this development we will begin with the
Hankel transform, including elementary properties
of Bessel functions from scratch. 
These functions were a must for quite a few
generations of mathematicians but not anymore.

To put things in perspective, let us start with 
the following multidimensional counterpart of (\ref{f11}).

\subsection{Selberg--Mehta--Macdonald formula}
Let $\BR^n$ be a Euclidean vector space with
scalar product $(\cdot,\cdot)$
and $R\subset \BR^n$ a reduced irreducible root system.
For any root
$\alpha,$ we set
$$\alpha^\vee =\frac{2\alpha}{(\alpha,\alpha)},\;
\tilde \alpha =\frac
{\sqrt{2}\alpha}{|\alpha|}.$$
Note that $(\tilde \alpha,\tilde \alpha)=2.$
We will need $\rho =\frac{1}{2}\sum_{\alpha >0}\alpha.$
Given $\alpha \in R$ and
$x\in \BR^n,$ let $\tilde x_\alpha =(x,\tilde \alpha).$
By  $dx$ we mean
the standard measure on the
Euclidean space $\BR^n$.

\begin{theorem}\label{t11}
(Mehta--Macdonald integral)
For any complex number $k$ such that
$\Rea k>-\frac{1}{h^\vee}$, where $h^\vee$ is
the dual Coxeter number of $R$, we
have the identity
\begin{equation}\label{f13}
\int_{\BR^n}\prod_{\alpha \in R}|\tilde
x_\alpha|^ke^{-(x,x)/2}dx=
(2\pi)^{n/2}\prod_{\alpha >0}
\frac{\Gamma (k(\rho,\alpha^\vee)+k+1)}{\Gamma(
k(\rho,\alpha^\vee)+1)}.
\end{equation}
\end{theorem}

This formula was conjectured by Mehta for the root
system of type $A_n$. Bombieri readily deduced it from
Selberg's integral.
Then Macdonald formulated it in \cite{M1} for
all root systems and proved for the classical
systems $B, C, D.$
Finally Opdam found a uniform proof using
the technique of shift operators \cite{O1}.

\medskip
{\bf Example.} Let us check the Mehta--Macdonald formula
for the root system of type $A_1$.
In this case $n=1,$ and there is only one positive root
$\alpha : \; (\alpha,\alpha)=2.$
So $\tilde \alpha=\alpha^\vee=\alpha,\;
\rho=\frac{\alpha}{2},\;
(\rho,\alpha^\vee)=1$. Let $u=\tilde x_\alpha$.
Then formula (\ref{f13}) reads as 
$$\int_{-\infty}^\infty|u|^{2k}e^{-u^2/4}\frac{du}{\sqrt{2}}=
(2\pi)^{1/2}\frac{\Gamma(2k+1)}{\Gamma(k+1)}.$$
Changing the variables and using the classical doubling formula
for the $\Gamma$-function
$$2^{2k}\Gamma(k+\frac{1}{2})=
\sqrt{\pi}\frac{\Gamma(2k+1)}{\Gamma(k+1)},$$
we get (\ref{f11}). Thus  (\ref{f13}) is a generalization
of formula (\ref{f11}) to higher rank root systems.

\subsection{Hankel transform}
For any complex number $k\not \in -\frac{1}{2}-
\BZ_+,$ we define the function
$$\phi^{(k)}(t)=\sum_{m=0}^\infty\frac{t^{2m}
\Gamma(k+\frac{1}{2})}{m!\Gamma
(k+m+\frac{1}{2})}.$$
The convergence is for all $t$ and fast.
It is clear that $\phi^{(k)}(0)=1$. Function $\phi^{(k)}(t)$
is related to
the classical Bessel function $J_k(t)$ by the formula
\begin{equation}\label{fbesselj}
\phi^{(k)}(t)=\Gamma(k+\frac{1}{2})t^{-k+1/2}J_{k-1/2}(2it).
\end{equation}

\begin{theorem}\label{t12}
For any complex $k$ with
$\Rea k>-\frac{1}{2},$
\begin{equation}\label{f14}
\int_{-\infty}^\infty e^{-x^2}\phi^{(k)}(\lambda x)
\phi^{(k)}(\mu x)|x|^{2k}
dx=\phi^{(k)}(\lambda \mu)e^{\lambda^2+\mu^2}
\Gamma(k+\frac{1}{2}),
\end{equation}
which reduces to (\ref{f11}) as $\lambda =\mu =0.$
\end{theorem}

\begin{definition}
The real Hankel transform of
an even function $f(x)$ is
given by the formula
$$\bH_{re}^k(f)(\lambda)=
\frac{1}{\Gamma(k+\frac{1}{2})}\int_{-\infty}^\infty
\phi^{(k)}(x\lambda)f(x)|x|^{2k}dx.$$
\end{definition}

The Hankel transforms are even functions of $\lambda$;
the transforms of odd $f(x)$ are zero.
Formula (\ref{f14}) states that
\begin{equation}\label{f15}
\bH_{re}^k(\phi^{(k)}(\mu x)e^{-x^2})=
e^{\mu^2}(\phi^{(k)}(\mu \lambda)
e^{\lambda^2}).
\end{equation}

\begin{definition}
The imaginary Hankel transform
of an even function $f(x)$ is
$$\bH_{im}^k(f)(\lambda)=
\frac{1}{i\Gamma(k+\frac{1}{2})}\int_{i\BR}
\phi^{(k)}(x\lambda)f(x)|x|^{2k}dx.$$
\end{definition}

Let $V_{re}^k$ be the linear span of the functions
$\phi^{(k)}(\lambda x)e^{-x^2}$ treated as functions
of $x\in \BR.$
Similarly,  $V_{im}^k$ is the span of
$\phi^{(k)}(\lambda x)e^{\lambda^2}$
considered as functions of $\lambda\in i\BR.$
It is clear that
$\bH_{re}^k$ maps $V^k_{re}$ to $V^k_{im}$ and
$\bH_{im}^k$ maps $V^k_{im}$
to $V^k_{re}$. As an immediate consequence of (\ref{f14}), 
we get the inversion formula.

\begin{corollary}\label{c11}
The maps $\bH_{re}^k$ and $\bH_{im}^k$ are
inverse to each other.
\end{corollary}

Here we can replace the spaces $V_{re}^k$ and $V_{im}^k$ by their
suitable completions. These completions are large enough
as can be seen from the following corollary.

\begin{corollary}\label{c12}
Let $U_{re}$ (respectively $U_{im}$)
be the linear spans
of functions in the form $p(x)e^{-x^2}$
(respectively $p(\lambda)e^{\lambda^2}$)
where $p$ runs through even polynomials.
Then $\bH_{re}$ maps $U_{re}$ to
$U_{im}$, $\bH_{im}$ maps $U_{im}$ to $U_{re}$,
$\bH_{im}\bH_{re}=id$ and
$\bH_{re}\bH_{im}=id$. Moreover,
$$\bH_{re}:\; U_{re}^{\le n} \mapsto
U_{im}^{\le n}, \
\bH_{im}:\; U_{im}^{\le n}\mapsto U_{re}^{\le n},
$$
if $U_{re}^{\le n}\subset U_{re}$
(respectively $U^{\le n}_{im}\subset U_{im}$) are subspaces
generated by
$p(x)e^{-x^2}$ (respectively $p(\lambda)e^{\lambda^2})$ for
all even polynomials $p$ of degree $\le n.$
\end{corollary}

{\bf Proof.} Differentiating formula (\ref{f15}) $\,2n$ times
in terms of the variable $\mu$ and setting
$\mu=0,$ we get that $\bH_{re}(x^{2n})$ is 
$e^{\lambda^2}$ times a polynomial in $\lambda.$
A better and more constructive proof
will be presented below. \sq\medskip

\subsection{Gaussian sums} Let $N$ be a natural number,
$q$  a primitive $N$-th root of 1. We will also need $q^{1/4},$
which will be picked in primitive $4N$-th roots of $1$. We
will consider the so-called {\em generalized Gaussian sums}
$$\tau =\sum_{j=0}^{2N-1}q^{j^2/4}.$$
Note that $q^{j^2/4}$ depends only on the residue of $j$
modulo $2N,$
so we can assume that the summation index runs through
$\BZ\modd 2N$. A natural choice for $q^{1/4}$ is
$e^{\pi \imath/2N}.$
In this case, we have the celebrated formula of
Gauss:
\begin{equation}\label{f16}
\tau=(1+\imath)\sqrt{N}.
\end{equation}

{\bf Comment.}
The standard definition of the
Gaussian sum is $\tau'=\sum_{j=0}^{n-1}q^{j^2}$. Gauss
proved that for $q=e^{2\pi i/n},$
$$\tau'=\left\{ \begin{array}{cl}\sqrt{n}&\mbox{if}\
n\modd 4=1,\\ i\sqrt{n}&
\mbox{if}\ n\modd 4=3,\\
 (1+i)\sqrt{n}&\mbox{if}\ n\modd 4=0,\\ 0&\mbox{if}\;
n\modd 4=2.\end{array}\right.$$
Formula (\ref{f16}) corresponds
to $n=4N.$  It somewhat resembles Fresnel's integral
$$\int_{-\infty}^\infty e^{ix^2}dx=(1+i)\sqrt{\frac{\pi}{2}},$$
which can be obtained from (\ref{f12}) by changing the variable
$x\mapsto \frac{1+i}
{\sqrt{2}}x$ and shifting the contour of integration. 
In a sense,
$\sqrt{2N}$ substitutes for  $\sqrt{\pi}$ when we switch 
from the real theory to roots of unity.
\sq
\bigskip

There are two levels in the theory of Gaussian sums:

(a) checking that $\tau \bar \tau=2N,$ which is simple
and entirely conceptual,

(b) finding $\arg(\tau),$ which is not.

Let us recall (a).
We consider the space $V$ of
complex valued functions on $\BZ\modd 2N,$
pick a primitive $2N$-th root
$q^{1/2}$ of $1,$ and introduce the Fourier transforms
$\bF_+$ and $\bF_-$ on $V$:
$$\bF_\pm (f)(\lambda)=\sum_{x=0}^{2N-1}q^{\pm x\lambda /2}f(x).$$
Then

(i) $\bF_-\bF_+\ =\ \bF_+\bF_-\ =\ \mbox{multiplication by}\; 2N,$

(ii) $\bF_+(\gamma)=\tau \gamma^{-1}$ and
$\bF_-(\gamma^{-1})=\bar \tau \gamma,$

where $\gamma (x)=q^{x^2/4}$ for a primitive $4N$-th root of unity
$q^{1/4}.$

Indeed,
$$\sum_{\lambda \in \BZ\modd 2N}q^{(x-y)\lambda /2}=2N$$
for $x=y$ and zero otherwise. We get (i). As for (ii), it suffices
to examine $\ \bF_+(\gamma)$
$$=\sum_{x\in \BZ\modd 2N}q^{x\lambda /2}q^{x^2/4}=
\sum_{x\in \BZ\modd 2N}
q^{(x+\lambda)^2/4}q^{-\lambda^2/4}=\tau q^{-\lambda^2/4}.$$

The relation $\tau \bar \tau =2N$ is immediate from (i) and (ii).
By the way, it is easy to show that $\bar \tau=Tr(\bF_+)$,
but it is not helpful for calculating $\arg(\tau).$

We are going to deduce (\ref{f16}) from a $q$-variant of
(\ref{f11}).
There will be two interesting corollaries:

(1) a formula for $\tau$ for any primitive root $q$,

(2) a generalization of (\ref{f16}) to arbitrary integral $k.$

\noindent
The former is closely connected with known formulas, 
the latter is new.

\bigskip
\setcounter{equation}{0}
\section{Imaginary integration}
In this section, we will discuss the imaginary integration,
which is called the {\em compact case} in the context of
symmetric spaces.

Let us first switch from $\BR$ to $i\BR.$
Recall our main formula
$$\int_{-\infty}^\infty e^{-x^2}|x|^{2k}dx=\Gamma (k+\frac{1}{2}),
\ \Rea k>-1/2.$$
Changing the variable $x\mapsto ix,$
\begin{equation}\label{f21}
\frac{1}{i}\int_{i\BR}e^{x^2}|x|^{2k}dx=\Gamma (k+\frac{1}{2}),
\ \Rea k>-1/2.
\end{equation}
These two formulas are of course equivalent
but their $q$-counterparts are not. It is similar to the
Harish-Chandra theory.
The real integration corresponds
to the so-called {\em noncompact case},
which is very different from the compact case.

{\bf Comments.} (i) 
There is an improvement of (\ref{f21}):
\begin{equation}\label{f22}
\int_{-\varepsilon +i\BR}e^{x^2}(-x^2)^kdx=\Gamma (k+\frac{1}{2})
\cos (\pi k),
\end{equation}
which holds for {\em all} complex $k$ provided
that $\varepsilon>0.$
It can be readily deduced from the well-known
classical definition of the $\Gamma$-function
$$\frac{1}{2}\int_C\frac{(-z)^k e^{-z}}{i\sqrt{-z}}dz=
\Gamma(k+\frac{1}{2})\cos (\pi k).
$$
Change the variable $z\mapsto -x^2.$ Here the path of the
integration $C$ begins at $z=-\varepsilon i+\infty$, moves to the
left down the positive real axis to $-\varepsilon i$, then circles
the origin and returns along the positive real axis to 
$\varepsilon i+\infty$.
The branches of the
$\log$ in the expression $(-z)^k=e^{k\log(-z)}$ and 
$\sqrt{-z}$ are standard, with the cutoff at 
-$\BR_+.$

(ii) The usage of the Gaussian $e^{x^2}$ makes the classical
formula in terms of $e^{-z}$ somewhat more elegant, 
but of course is not significant. It is
different in the $q$-theory, which does require the 
formula in terms of the Gaussian. This 
is directly connected with the appearance of the double
Hecke algebra. The conjugation by Gaussian is 
an automorphism of the latter. The substitution $z=-x^2$
does not seem meaningful from the viewpoint of this algebra.

\subsection{Macdonald's measure}
Let $q$ be a real number, $0<q<1$. We set
$q=e^{-1/a}$ for a real number $a>0.$ The $q$-counterpart of
$x^{2k}$ is the following function of
complex variable $x$:
\begin{equation}\label{f23}
\delta_k(x)=\prod_{j=0}^\infty
\frac{(1-q^{2x+j})(1-q^{-2x+j})}{(1-q^{2x+j+k})
(1-q^{-2x+j+k})}.
\end{equation}
If $k$ is a non-negative integer, formula (\ref{f23}) reads as
\begin{equation}\label{f24}
\delta_k(x)=\prod_{j=0}^{k-1}(1-q^{2x+j})(1-q^{-2x+j}).
\end{equation}

\begin{theorem}\label{t21}
Assume that $\Rea k>0$. Then $\delta_k(x)$ is
regular for $x\in i\BR$ and 
\begin{equation}\label{f25}
\frac{1}{i}\int_{i\BR}q^{-x^2}\delta_k(x)dx=
2\sqrt{\pi a}\prod_{j=0}^\infty
\frac{1-q^{k+j}}{1-q^{2k+j}}.
\end{equation}
\end{theorem}

{\bf Comment.} The right and left-hand sides of formula (\ref{f25})
are well defined
for any number $k$ such that $\Rea k\not \in -\BZ_+.$
However formula (\ref{f25}) {\em is
not} true when $\Rea k<0$.
\sq\medskip

Before proving the theorem,
let us get (\ref{f21}) as the limit of (\ref{f25})
when $q\to 1.$

\begin{proposition}\label{p21}
As $q\to 1,$ the leading term of the left-hand side of
formula (\ref{f25}) is
$\left( \frac{a}{4}\right)^{-k}\sqrt{a}
\int_{i\BR}e^{z^2}(-z^2)^kdz,$ and
the leading term of the right-hand side of this
formula is $2a^{-k}\sqrt{\pi a}
\frac{\Gamma(2k)}{\Gamma(k)}$.
Hence (\ref{f25}) implies (\ref{f21}).
\end{proposition}

{\bf Proof.} Introduce the $q$-counterpart of the 
$\Gamma$-function:
$$\Gamma_q(x)=(1-q)^{1-x}\prod_{j=0}^\infty
\frac{1-q^{x+j}}{1-q^{1+j}}.$$
It is well-known (and not difficult to prove) that
$\lim_{q\to 1}\Gamma_q(x)=
\Gamma (x)$. Hence the right-hand side of formula (\ref{f25}) as
$q\to 1$ is $2a^{-k}\sqrt{\pi a}\frac{\Gamma (2k)}{\Gamma (k)}$
asymptotically.

Now consider the left-hand side
of formula (\ref{f25}). Upon the change of variable
$x\mapsto \sqrt{a}z,$
$$\frac{1}{i}\int_{i\BR}q^{-x^2}\delta_k(x)dx=
\frac{1}{i}\int_{i\BR}e^{z^2}
\delta_k(\sqrt{a}z)\sqrt{a}dz.$$

\begin{lemma}\label{l21}
(Stirling--Moak) One has
$$\lim_{a\to \infty}\delta_k(\sqrt{a}z)
\left( \frac{a}{4}\right)^k=(-z^2)^k,$$ 
where the standard branch of
the logarithm is taken for $(-z^2)^k$ and $k$ is
arbitrary complex.
\end{lemma}

It is obvious in the case $k\in \BZ_+$. Indeed, asymptotically,
$$\delta_k(\sqrt{a}z)=
\prod_{j=0}^{k-1}(1-e^{j/a}e^{2z/\sqrt{a}})(1-e^{j/a}
e^{-2z/\sqrt{a}})\approx \left(-\frac{4z^2}{a}\right)^k.$$
We omit the proof for general $k.$ Actually this
formula can be applied only if one can estimate 
the remainder. The most convenient formula
for the latter is due to Moak. We will not discuss
it because it is applied in the presence of the Gaussian 
in the integrand when the convergence is simple.\sq\medskip

It follows from the lemma that the left-hand side of 
formula (\ref{f25}) is asymptotically
$\left( \frac{a}{4}\right)^{-k}\sqrt{a}\int_{i\BR}e^{z^2}
(-z^2)^kdz$. Comparing it with the asymptotical 
behavior of the 
right-hand side and using the doubling formula
$$2^{2k-1}\Gamma(k+1/2)=\sqrt{\pi}\Gamma(2k)/\Gamma(k),$$
we get (\ref{f21}) from (\ref{f25}). \sq\medskip

\subsection{Meromorphic continuations} 
It is instructional to calculate the difference between the
right-hand side and the left-hand side of formula (\ref{f25}) as 
$-1<k<0$. 
It is based on the general technique of meromorphic 
continuation of the integrals depending on parameters "to the
left" and can be applied to many functions in place of the
Gaussian. 
\medskip

{\bf Comment.}
A good example is the integral
$\int_{i\BR}(q^{x^2}-1)^{-1}\delta_k(x)dx,$
i.e.,  a certain $q$-variant of Riemann's zeta,
up to $q$-Gamma factors. 
It has a meromorphic
continuation in $k$ to the left 
(and other interesting properties).
See \cite{C9}.
Surprisingly, it {\em does not} converge to
$\zeta(k+1/2)\Gamma(k+1/2)$ for $\Rea k<0,$ 
as one can expect.
Its limit as $a\to \infty$ is a combination of 
$\Gamma$-factors (the $\zeta$ does not appear at all). 
We can get the desired limiting behavior
by switching from $(q^{x^2}-1)^{-1}$
to $(q^{x^2}+1)^{-1}.$ In this case, the integral
does converge to what can be expected.
Calculating these and similar limits requires
the technique of analytic continuation based
on the shift operator, which will be considered later
in this section.
\sq\medskip

We will involve
the intermediate function
$$f_\eps(k)=\frac{1}{i}\int_{\eps+i\BR}q^{-x^2}\delta_k(x)dx,$$
where $\eps>0$.
Let $f(k)$ be the left-hand side of (\ref{f25}), $g(k)$ the
right-hand side. The difference $f_\eps-f$ is $2\pi i$ 
times the sum of the
residues of $\frac{1}{i}q^{-x^2}\delta_k(x)$ over its
poles   in the strip
$0<\Rea x<\eps.$ The calculation depends on the domain of $k.$
For instance, $f_\eps$ and $f(k)$ coincide for $\Rea k>2\eps$
because the strip contains no $x$-poles for such $k.$

Let $0<\Rea k<2\eps.$ Then
the poles are $\frac{k}{2}+n\pi ia,\; n\in \BZ$.
So
$$f_\eps=f+\Pi_k,\; \Pi_k=\frac{2\pi i}{i}\Res_{x=k/2}\delta_k(x)
\sum_{n=-\infty}^\infty q^{-(k/2+n\pi ia)^2}.$$
Since $\lim^{\,x=k/2+\eps\,}_{\ \,\eps\to 0}\,(1-q^{-2x+k})=-a/2,$
$$\Res_{x=k/2}= -(a/2)\prod_{j=0}^\infty
\frac{(1-q^{k+j})(1-q^{-k+j})}{(1-q^{2k+j})(1-q^{j+1})}.$$
On the other hand,
\begin{equation}
\sum_{n=-\infty}^\infty q^{-(k/2+n\pi ia)^2}=
(2\sqrt{\pi a})^{-1}\sum_{m=-\infty}^\infty q^{m^2}q^{mk},
\label{ffunceq}
\end{equation}
thanks to the functional equation for the theta function.
Finally,

$$\Pi_k=-\frac{\sqrt{\pi a}}{2}\left(
\sum_{m=-\infty}^\infty q^{m^2}q^{mk}\right)
\prod_{j=0}^\infty
\frac{(1-q^{k+j})(1-q^{-k+j})}{(1-q^{2k+j})(1-q^{j+1})}.$$

If $-2\eps<\Rea k<0,$
the poles are $-\frac{k}{2}+n\pi ia,\; n\in \BZ,$ and
$f_\eps=f-\Pi_{-k}.$ The minus of $-\Pi_{-k}$ reflects the
change $a/2\mapsto -a/2$ in the calculation of the residues.

The function $f_\eps(k)$ is well defined for $|\Rea k|<2\eps$ and
equals $g(k)+\Pi_k$ in this strip. Recall that $g(k)$ is 
given by the product formula
meromorphic for all complex $k,$ and coinciding with $f(k)$ as
$\Rea k>0.$
Therefore
\begin{equation}\label{ffkgk}
f(k)=f_\eps+\Pi_{-k}=
g(k)+\Pi_k+\Pi_{-k} \hbox{\ for\ } -1<\Rea k<0.
\end{equation}

\subsection{Using the constant term}
We can reformulate Theorem \ref{t21} in a purely algebraic way.
Indeed, the function $q^x$ is periodic with 
the  period $\omega=2\pi ia$.
So is $\delta_k(x).$ Actually its period is $\omega/2$ but
we will use $\omega$ in the following calculation:
$$\frac{1}{i}\int_{i\BR}q^{-x^2}\delta_k(x)dx=
\frac{1}{i}\int_0^\omega
\sum_{n\in \BZ}q^{-(x+n\omega)^2}\delta_k(x)dx$$
$$=\frac{1}{i}\int_0^\omega
\sum_{n\in \BZ}q^{n^2/4+nx}\frac{\delta_k(x)}{2\sqrt{\pi a}}dx.$$

The second equality results from the functional equation for the
theta function, in a form slightly different from
(\ref{ffunceq}).

Since
$$\int_0^\omega q^{nx}dx=
\left\{ \begin{array}{clr}0&\mbox{if}&n\neq 0,\\
2\pi ia&\mbox{if}&n=0,\end{array}\right.$$
we need to represent $\delta_k(x)$
as a Laurent series and then find the {\em constant term} $CT$ of 
$\widehat{\gamma_-}(x)\widehat{\delta}_k(x),$ where
$$\widehat{\gamma_-}(x)\equal\sum_{n\in \BZ}q^{n^2/4+nx},$$
$\widehat{\delta}_k(x)$ is the expansion of
$\delta_k(x)$ as a Laurent series of $q^x.$
To be more exact, we will take the expansion with the 
coefficients which are power series in terms of the 
variables $q, q^k.$ Then it is determined uniquely.

Generally speaking, a product of two Laurent series is not 
well defined, but here we have perfect
convergence thanks to the $q^{n^2/4}$
from $\widehat{\gamma_-}(x).$ Finally,

$$\frac{1}{i}\int_{i\BR}q^{-x^2}\delta_k(x)dx=\frac{2\pi ia}
{i2\sqrt{\pi a}}CT(\widehat{\gamma_-}(x)\widehat{\delta}_k(x)),$$
and Theorem \ref{t21} is equivalent to

\begin{theorem}\label{t22}
We have an equality:
$$CT\left(\widehat{\gamma_-}(x)\widehat{\delta}_k(x)\right)=
2\prod_{j=0}^\infty \frac{1-q^{k+j}}{1-q^{2k+j}}.$$
\end{theorem}

The expansion of $\delta_k$ can be calculated explicitly.
Let us introduce the function
\begin{equation}\label{formmu}
\mu_k(x)=\prod_{j=0}^\infty
\frac{(1-q^{2x+j})(1-q^{-2x+j+1})}{(1-q^{2x+k+j})
(1-q^{-2x+k+j+1})}.
\end{equation}
It is closely connected with $\delta_k:$
\begin{equation}\label{fmudelta}
\delta_k(x)=\frac{\mu_k(x)+\mu_k(-x)}{1+q^k}.
\end{equation}
Following Macdonald, we may call it 
{\em truncated theta function} since
for $k\in \BZ_+,$
$$\mu_k(x)=\prod_{j=0}^{k-1}(1-q^{2x+j})(1-q^{-2x+j+1}).$$
As $\BZ_+\ni k\to\infty,$ it becomes the classical theta function.

\begin{theorem}\label{t23}
(Ramanujan's $_1\Psi_1$-summation)
$$\mu_k(x)/CT(\mu_k(x))=
1+\frac{q^k-1}{1-q^{k+1}}(q^{2x}+q^{1-2x})+\cdots$$
$$+\frac{(q^k-1)(q^k-q)(q^k-q^2)\cdots
(q^k-q^m)}{(1-q^{k+1})(1-q^{k+2})
(1-q^{k+3})\cdots (1-q^{k+m})}(q^{2mx}+q^{m-2mx})+\cdots .$$
\end{theorem}

{\bf Proof.} Use
$\mu_k(x+\frac{1}{2})=\frac{q^{2x+k}-1}{q^{2x}-q^k}\mu_k(x).$
 \sq\medskip

\begin{theorem}\label{t24}
(Constant Term Conjecture) 
\begin{equation}\label{f26}
CT(\mu_k(x))=\frac{(1-q^{k+1})^2(1-q^{k+2})^2\cdots }
{(1-q^{2k+1})(1-q^{2k+2})\cdots (1-q)(1-q^2)\cdots }.
\end{equation}
\end{theorem}

{\bf Proof.} First,
\begin{equation}\label{f27}
CT(\mu_{k+1}(x))=\frac{(1-q^{2k+1})(1-q^{2k+2})}
{(1-q^{k+1})^2}CT(\mu_k(x)).
\end{equation}
Indeed, $\mu_{k+1}(x)=(1-q^{2x+k})(1-q^{-2x+k+1})\mu_k(x)$ and
Theorem \ref{t23} results in
$$CT(\mu_{k+1}(x))=
CT((1-q^{2x+k}-q^{-2x+k+1}+q^{2k+1})\mu_k(x))$$
$$=(1+q^{2k+1})CT(\mu_k(x))-2q^{k+1}\frac{q^k-1}
{1-q^{k+1}}CT(\mu_k(x))$$
$$=\frac{CT(\mu_k(x))}{1-q^{k+1}}
((1+q^{2k+1}-q^{k+1}-q^{3k+2}-2q^{2k+1}+
2q^{k+1})$$
$$=CT(\mu_k(x))\frac{(1+q^{k+1})(1-q^{2k+1})}{1-q^{k+1}}.$$

Let us denote the right-hand side of (\ref{f26}) by $c_k.$
Then $CT(\mu_k)/c_k$ is a periodic function in terms of
$k$ with the period $1.$ Hence the expansion of $c_k$
in terms of $q,q^k$ is invariant under the substitution
$q^k\mapsto qq^k.$  So it does not depend on
$k$ and is a series in terms of $q.$ It must
be $1$ because $\mu_0=1.$  \sq\medskip

{\bf Comment.} The above calculation is well-known
(see \cite{An, AI}). The celebrated constant term conjecture 
due to Macdonald is for arbitrary root systems.
It was proved first for $A_n,$ then for $BC_n$ by
Kadell, and then for $G_2,F_4$ using computers.
The case of $E_6$
appeared to be beyond the capacity of modern computers.
The proof from \cite{C1} is based on the shift operators,
generalizing those introduced by Opdam in the differential setup.
Let us formulate the statement for an arbitrary reduced
irreducible root system $R=\{\alpha\}\subset \BR^n.$

In terms of the one-dimensional $\mu$ above:
$$\mu^R_k(x)=\prod_{\alpha >0}\mu_k(x_{\alpha}/2).$$
The Macdonald conjecture reads as
$$CT(\mu_k^R(x))=\prod_{j=1}^\infty \prod_{\alpha >0}\frac{(1-
q^{k(\rho,\alpha^{\vee})+j})^2}
{(1-q^{k(\rho,\alpha^\vee)+j+k})
(1-q^{k(\rho,\alpha^\vee)+j-k})}.$$

There are no good formulas for other coefficients of
$\mu^R_k.$ Generally speaking, they are not $q$-products.
It was the main problem with calculating $CT(\mu^R).$
If Theorem \ref{t23} existed, it would be straightforward.

\subsection{Shift operator}
We are going to prepare the main tool
for proving Theorem \ref{t22}.
The {\em shift operator} is
$$S(f(x))\equal \frac{f(x-\frac{1}{2})-f(x+\frac{1}{2})}
{q^x-q^{-x}}.$$
It is a $q$-variant of the differentiation
$x^{-1}d/dx.$ It plays an important role in the theory of the
so-called basic hypergeometric function.
See e.g., [AI].

Let $g_k=\prod_{j=0}^\infty\frac{1-q^{k+j}}
{1-q^{2k+j}}$. Given a function $f(x)$ defined on the imaginary
line, its {\em q-Mellin transform} is introduced by the formula
$$\Psi_k(f)=\frac{1}{ig_k}\int_{i\BR}f(x)\delta_k(x)dx.$$
The function $\Psi_k$ is an analytic function of the variable
$k$ in the half-plane $\Rea k>0,$ provided we have the integrability.

Let us examine the behavior of the
$q$-Mellin transform under the  $q$-variant of the integrating by parts.

\begin{theorem}\label{t25}
Assume that the function $f(x)$ is
analytic in an open
neighborhood of the strip $|\Rea x|\le 1$.
Provided that the
integrals below are well defined,
\begin{equation}\label{f28}
\Psi_k(f)=(1-q^{k+1})\Psi_{k+1}(f)+q^{k+3/2}\Psi_{k+2}(S^2(f)).
\end{equation}
\end{theorem}

{\bf Proof} is direct. Later we will give a better
one, with modest calculations. However the theorem
will be applied to various classes of functions, so 
an explicit proof is helpful to control 
the analytic matters.

First, we check that $\ S^2(f)\ $
$$=\frac{f(x-1)}{(q^x-q^{-x})(q^{x-1/2}-q^{-x+1/2})}
+\frac{f(x+1)}{(q^x-q^{-x})(q^{x+1/2}-q^{-x-1/2})}$$
$$-\frac{f(x)}{q^x-q^{-x}}\left(\frac{1}
{q^{x-1/2}-q^{-x+1/2}}+
\frac{1}{q^{x+1/2}-q^{-x-1/2}}\right).$$

Second, we will
change the variables and move the contour of integration.
The resulting formulas are as follows:

\begin{align}
&\int_{i\BR}S^2(f)\delta_k(x)dx \notag\\
&=\int_{i\BR}\frac{f(x)\delta_k(x+1)dx}
{(q^{x+1}-q^{-x-1})(q^{x+1/2}-
q^{-x-1/2})}-\\
&-\int_{i\BR}\frac{f(x)\delta_k(x)dx}
{(q^x-q^{-x})(q^{x+1/2}-q^{-x-1/2})}+\notag\\
&+\int_{i\BR}\frac{f(x)\delta_k(x-1)dx}
{(q^{x-1}-q^{-x+1})(q^{x-1/2}-
q^{-x+1/2})}-\notag\\
&-\int_{i\BR}\frac{f(x)\delta_k(x)dx}
{(q^x-q^{-x})(q^{x-1/2}-
q^{-x+1/2})}.\notag
\end{align}

Let us denote the first two terms in the 
right-hand side by $A$ and
the second two terms (lines) by $B.$
So the integral is $A+B.$

Third,
$$\delta_k(x+1)=
\frac{(1-q^{2x+k})(1-q^{2x+1+k})(1-q^{-2x-2})(1-q^{-2x-1})}
{(1-q^{-2x-2+k})(1-q^{-2x-1+k})
(1-q^{2x})(1-q^{2x+1})}\delta_k(x),$$
and
$$\delta_k(x-1)=
\frac{(1-q^{-2x+k})(1-q^{-2x+1+k})(1-q^{2x-2})(1-q^{2x-1})}
{(1-q^{2x-2+k})(1-q^{2x-1+k})
(1-q^{-2x})(1-q^{-2x+1})}\delta_k(x).$$

Fourth,
$$A=-\int_{i\BR}\frac{(q^{x+1/2}+q^{-x-1/2})
(1-q^{2k-1})f(x)\delta_k(x)dx}
{(1-q^{-2x-2+k})(1-q^{-2x-1+k})(q^x-q^{-x})q^{2x+1}},$$
and
$$B=\int_{i\BR}\frac{(q^{x-1/2}+q^{-x+1/2})
(1-q^{2k-1})f(x)\delta_k(x)dx}
{(1-q^{2x-2+k})(1-q^{2x-1+k})(q^x-q^{-x})q^{-2x+1}}.$$

Now,
$$\frac{A+B}{(1-q^{2k-1})}=\int_{i\BR}Gf(x)\delta_k(x)dx$$
where
$$G=$$ $$\frac{1+q+q^k+q^{2x}+q^{-2x}-q^{k-2}-
q^{2k-3}-q^{2k-2}-q^{2x+2k-2}-q^{-2x+2k-2}}
{(1-q^{2x-2+k})(1-q^{2x-1+k})(1-
q^{-2x-2+k})(1-q^{-2x-1+k})q^{3/2}}.$$

Regrouping the terms:
$$\frac{A+B}{(1-q^{2k-1})}=(q^{1-k}-q^{k-1})
\int_{i\BR}\frac{q^{-1/2}f(x)\delta_k(x)dx}
{(1-q^{2x-1+k})(1-q^{-2x-1+k})}$$
$$+\int_{i\BR}\frac{q^{-3/2}(1-q^{2k-3})
(1+q^{k-1})(1+q^{2-k})f(x)\delta_k(x)dx}{(1-
q^{2x-2+k})(1-q^{2x-1+k})(1-q^{-2x-2+k})
(1-q^{-2x-1+k})}.$$

Taking into account that
$$\delta_{k-1}(x)=\frac{\delta_k(x)}
{(1-q^{2x-1+k})(1-q^{-2x-1+k})},\ \delta_{k-2}(x)$$
$$=\frac{\delta_k(x)}
{(1-q^{2x-2+k})(1-q^{2x-1+k})
(1-q^{-2x-2+k})(1-q^{-2x-1+k})},$$
and that
$$\frac{g_k}{g_{k-1}}=(1-q^{2k-1})(1+q^{k-1}),$$
$$\frac{g_k}{g_{k-2}}=
(1-q^{2k-3})(1+q^{k-2})(1-q^{2k-1})(1+q^{k-1}),$$
we come to the formula
$$q^{k-1/2}\Psi_k(S^2(f))=
\Psi_{k-2}(f)-(1-q^{k-1})\Psi_{k-1}(f),$$
which is equivalent to the statement of the theorem. \sq\medskip


\subsection{Applications}
The simplest example is $f(x)=1.$ We get the formula
$$CT(\delta_k(x))=
\frac{(1-q^{k})(1-q^{2k+1})}{(1-q^{k+1})(1-q^{2k})}
(CT(\delta_{k+1}(x)),$$
which is equivalent to the formula (\ref{f27}) above.

Now let $f(x)=q^{-x^2}$.
It is easy to check that $S(q^{-x^2})=
q^{-1/4}q^{-x^2}$. Therefore
$$\Psi_k=(1-q^{k+1})\Psi_{k+1}+q^{k+1}\Psi_{k+2}
$$
for $\Psi_k=\Psi_k(q^{-x^2})$ or, equivalently,
$$\Psi_k-\Psi_{k+1}=-q^{k+1}(\Psi_{k+1}-\Psi_{k+2}).$$

Introducing the function
$$\phi_k=(-1)^kq^{(k+1)k/2}(\Psi_k-
\Psi_{k+1}),$$ 
we see that it is periodic with the period $1$
in the variable $k.$ In particular, 
it can be extended to the whole complex plane, where it
is analytic.
The function $\phi_k$ is
also quasi-periodic with the period $\omega$ because
$\Psi_k$ is $\omega$-periodic by construction:
$$
\phi_{k+\omega}=q^{k\omega+\omega^2/2}\phi_k\ =\
e^{-2\pi i k}e^{-\pi i \omega}\phi_k.
$$
Recall that $\omega=2\pi ia, q=e^{-1/a}.$ 
These are the defining properties of the classical
$$
\vartheta_3(k)=\sum_{m=-\infty}^{\infty}e^{2\pi i m k}
e^{\pi i m^2\omega}. 
$$
For instance, $\phi_k$ has only one
zero in the fundamental parallelogram if it is not zero
identically. The zero is known to be $(\omega+1)/2.$

However it is straightforward to check that the $\phi_k$
vanish at $k=0.$ Indeed,
$$\Psi_0=\frac{2}{i}\int_{i\BR}q^{-x^2}dx=2\sqrt{\pi a},$$
$$\Psi_1=\frac{1}{i(1-q)}
\int_{i\BR}q^{-x^2}(1-q^{2x})(1-q^{-2x})dx
$$
$$=\frac{1}{i(1-q)}
\left(2\int_{i\BR}q^{-x^2}dx-\int_{i\BR}q^{-x^2+2x}dx-
\int_{i\BR}q^{-x^2-2x}dx\right)$$
$$=\frac{2}{i}\int_{i\BR}q^{-x^2}dx=\Psi_0.$$

Therefore $\phi=0$ and
$\Psi_k$ is a constant, which is the value of 
$\Psi=\Psi_k(q^{-x^2})$ at zero, that is $2\sqrt{\pi a}.$ \sq\medskip

{\bf Comment.}
The above argument can be simplified
using a reformulation of Theorem \ref{t22} 
in terms of the Laurent series.
First, we check a variant of Theorem \ref{t25}
for $\omega$-periodic even functions $f(x)$ and the
integration over the imaginary period. Then
we replace integrating over the period by taking
the constant term and switch entirely to the 
Laurent expansions.
Third, we apply the shift--formula to
$f(x)=\widehat{\gamma_-}(x).$ In this setting,
$\phi_k$ is not considered as an analytic
function, but becomes a formal series in terms of the variables
$q$ and $q^k.$ For such series, the equality 
$\phi_{k+1}=-q^{k+1}\phi_k$
immediately implies that $\phi_k=0$.

\bigskip
\setcounter{equation}{0}
\section{Jackson and Gaussian sums}
We are going to change the imaginary integration in the main
formula of the previous section to the Jackson summation.

Recall that the
imaginary $q$-Mellin transform is
\begin{align}
&\Psi_k(f)=\frac{1}{\imath g_k}\int_{\imath \BR}
f(x)\delta_k(x)dx,\ g_k=\prod_{j=0}^\infty
\frac{1-q^{k+j}}{1-q^{2k+j}},\notag\\
&\delta_k(x)=\prod_{j=0}^\infty
\frac{(1-q^{2x+j})(1-q^{-2x+j})}
{(1-q^{2x+k+j})(1-q^{-2x+k+j})}.\notag
\end{align}

Using this transform,
Theorem \ref{t21} reads as 
$$\Psi_k(q^{-x^2})=2\sqrt{\pi a},\hbox{\ where\ }
q=e^{-1/a},\, \Rea k>0.
$$

The next theorem will be very close to Theorem \ref{t21}.
There was one place in its proof where we used special 
features of the imaginary integration more than absolutely
necessary. Let us first somewhat improve it.

We introduced 
$\phi_k=\Psi_k(q^{-x^2})-\Psi_{k+1}(q^{-x^2}),$
checked that $\phi_k=-q^{k+1}\phi_{k+1},$ and found that
$\phi_0=0$ and therefore $\phi_k=0$ for all $k\in \BZ_+.$
This part remains unchanged. Now it 
is easy to see that the function
$\Psi_k(q^{-x^2})$ is analytic in terms of 
the $K=q^k$ considered as a new variable
in a neighborhood of $K=0.$ So is $\phi_k.$ 
However the latter has infinitely many zeros in a 
neighborhood of $K=0,$ which results in $\phi_k=0$. 
The rest of the proof is the same as in the previous 
section. 

\subsection{Sharp integration} For an analytic function
$f(x)$ in a neighborhood of the positive real axis,
we define its {\em sharp $q$-Mellin transform}
\begin{equation}\label{fpsish}
\Psi_k^\#(f)=\frac{1}{\imath g_k}\int_Cf(x)\delta_k(x)dx,
\end{equation}
where the path of the integration $C$ begins at
$z=-\eps \imath+\infty$, moves to the left down the positive
real axis to $-\eps \imath$, then moves up to $\eps \imath$
and returns along the positive real axis to
$\eps \imath+\infty.$ 

The behavior of $C$ near $0$ is not important
in the $q$-theory. Recall that the classical paths of 
this type used for $\Gamma$ and $\zeta$ 
{\em must} go around zero. 
Our $C$ is like a pencil aimed at zero from $+\infty$;
that is why we call it sharp.

\begin{theorem}\label{t31}
Provided that $\Rea k>0$ and $|\Ima k|<2\eps$  (equivalently,
$k/2$ sits inside $C),$ 
\begin{equation}\label{f31}
\Psi_k^\#(q^{x^2})=(-a\pi )\prod_{j=1}^\infty
\frac{(1-q^{j+k})(1-q^{j-k-1})}
{(1-q^j)^2}\sum_{j=-\infty}^\infty q^{(k-j)^2/4}.
\end{equation}
\end{theorem}

{\bf Comment.} 
The sum $\sum_{j=-\infty}^\infty q^{(k-j)^2/4}$
can be expressed as the infinite product
$$q^{k^2/4}\prod_{j=1}^\infty (1-q^{j/2})(1+q^{j/2-1/4+k/2})
(1+q^{j/2-1/4-k/2})$$ by Jacobi's triple product formula.
\sq\medskip

{\bf Proof}  begins with the following lemma.
Recall that $q=e^{-1/a}$ and $\omega =2\pi \imath a.$

\begin{lemma}\label{l31}
Assume that $f(x)$ is analytic for
all $x\in \BC$ and  that
$\lim_{\xi \to \infty}f(x+\xi)e^{c\xi}=0$
for arbitrary $c>0$ and $x\in \BC.$  Then
the function $\Psi_k^\#(f)$ has an analytic continuation
to all $k\in \BC$
and, moreover, vanishes at the points from 
$\BZ_++\BZ \omega.$ 
\end{lemma}

{\bf Proof.} By Cauchy's Theorem, the integral (\ref{f31})
is the following sum of residues:
$$\Psi_k^\#(f)=(-a\pi )\prod_{j=1}^\infty
\frac{(1-q^{j-k-1})}{(1-q^j)}\times$$
$$\times \sum_{j=0}^\infty
\frac{1-q^{j+k}}{1-q^k}\prod_{l=1}^j\frac{1-q^{l+2k-1}}
{1-q^l}f\left(\frac{k+j}{2}\right)q^{-kj}.$$
The right-hand side gives the desired analytic continuation.
\sq\medskip

Using the above formula we can reformulate
Theorem \ref{t31} in a purely algebraic way:
\begin{align}\label{f32}
&\sum_{j=0}^\infty q^{(k-j)^2/4}\frac{1-q^{j+k}}{1-q^k}
\prod_{l=1}^j\frac{1-q^{l+2k-1}}{1-q^l}\\
&=\prod_{j=1}^\infty
\frac{1-q^{j+k}}{1-q^j}\left( \sum_{j=-\infty}^\infty
q^{(k-j)^2/4}\right).\notag
\end{align}

The left-hand side is a {\em Jackson sum\,}, i.e., the summation
of the values of a given function (here $q^{(k-j)^2/4}$)
over $\BZ$ with some weights. There is a proof of the theorem 
directly in terms of the Jackson summation, without using the
sharp integration. It will not be discussed here (see
\cite{C5}).

{\bf Comment.} (i) Formula (\ref{f32}) obviously holds
for $k=0,$ to be more exact, as $k\to 0$ (we have $1-q^k$ in
the denominator). Indeed, we get
$$1+2\sum_{j=1}^\infty q^{j^2/4}=
\sum_{j=-\infty}^\infty q^{j^2/4}.$$
Hence, the difference between the left-hand and 
the right-hand sides of formula (\ref{f32}) has a zero of 
the {\em second} order at $k=0$.

(ii) For $k=1,$ formula (\ref{f32}) reads as 
$$\sum_{j=0}^\infty q^{(1-j)^2/4}\frac{(1-q^{j+1})^2}{1-q}
=\sum_{j=-\infty}^\infty q^{j^2/4},$$
and can be checked by a simple calculation.

(iii) It holds for $k=-1/2$. Indeed, in this case
the left-hand side is $2q^{1/16}$ and the right-hand side 
can be transformed using Jacobi's triple product
formula as follows:
\begin{align}
&q^{1/16}\prod_{j=1}^\infty (1-q^{j/2})(1+q^{j/2})
(1+q^{j/2-1/2})\prod_{j=1}^\infty
\frac{1-q^{j-1/2}}{1-q^j}\notag\\
&=q^{1/16}\prod_{j=1}^\infty
(1+q^{j/2-1/2})(1-q^{j-1/2})\notag\\
&=2q^{1/16}\prod_{j=1}^\infty
\frac{(1-q^{j})(1-q^{j-1/2})}{(1-q^{j/2})}=
2q^{1/16}.\notag
\end{align}

Let us denote the 
right-hand side of formula (\ref{f31}) by
$\Pi_k$. We have the following properties:

(a) $\Pi_k=-q^{k+1}\Pi_{k+1},\ \ $ (b) $\Pi_{k+2\omega}=q^{\omega^2+k\omega}\Pi_k,$

(c) the zeros of $\Pi_k$ are  
$\{ 0, \omega, 1/2+\omega \}\,\mod\, \BZ +2\omega \BZ.$

\noindent In (c), all zeros are simple. 

These properties determine
$\Pi_k$ uniquely.
The only property which is not immediate is (c). It can be proved
using Jacobi's triple product or deduced from (a),(b). 
Indeed, the latter give that $\Pi_k$ 
has three zeros inside the 
parallelogram of periods with the sum  $1/2.$ 
Two of them are obvious: $\{ 0,\omega\}.$

\subsection{Sharp shift--formula}
Let $\Psi_k^\#\equal \Psi_k^\#(q^{x^2})$ be 
the left-hand side of
formula (\ref{f31}). It is clear that

($b'$) $\Psi_{k+2\omega}^\#=q^{\omega^2+k\omega}\Psi_k^\#,$

($c'$) $\Psi_k^\#$ has zeros at $0$ and $\omega.$

Thus if we prove 
$(a')$  $\Psi_k^\#=-q^{k+1}\Psi_{k+1}^\#,$ 
then $\Pi_k=c\Psi_k^\#$ for a constant $c,$ which 
has to be $1$ (use the normalization).
So we need a variant
of the shift--formula (\ref{f28}) for the sharp integration.

\begin{lemma}\label{lsharpshift}
\begin{equation}\label{fsharpshift}
\Psi_k^\#(f)=(1-q^{k+1})\Psi_{k+1}^\#(f)+q^{k+3/2}
\Psi_{k+2}^\#(S^2(f)),
\end{equation}
provided we have the existence of the integrals.
Here  the shift operator is
$S(f)(x)=\frac{f(x-1/2)-f(x+1/2)}{q^x-q^{-x}},$
$f(x)$ is an even function continuous on
the sharp integration path $C$ and, moreover, 
analytic in the domain
$$
\{x \mid  -1-\delta <\Rea x<1+\delta,\, 
-\epsilon< \Ima x < \epsilon \}\, \hbox{\ for\ } \delta>0.
$$
\end{lemma}

{\bf Proof} is a straightforward 
adjustment of that in the imaginary case.
The analyticity is necessary to ensure the invariance of
the integration with respect to the shifts by $\pm 1.$ 
\sq\medskip

Using the formula $S(q^{x^2})=q^{1/4}q^{x^2}$
(note the change of sign, compared with the previous section),
$$\Psi_k^\#=(1-q^{k+1})\Psi_{k+1}^\#+q^{k+2}\Psi_{k+2}^\#.$$
So the function $\phi_k=\Psi_k^\#+q^{k+1}\Psi_{k+1}^\#$ is periodic:
$\phi_{k+1}=\phi_k$. On the other hand, ($b'$) results in
$\phi_{k+2\omega}=q^{\omega^2+k\omega}\phi_k$. Combining with
the $1$-periodicity, we conclude that
$\phi_k$ has only one zero in the parallelogram of periods.
However we already know that $\phi_k$ has zero of order $2$
at $k=0$. Hence 
\begin{equation}\label{fpsisharp}
\phi_k=0\Rightarrow  \Psi_k^\#=-q^{k+1}\Psi_{k+1}^\#
\Rightarrow \Psi_k^\#=\Pi_k,
\end{equation}
and Theorem \ref{t31} is proved. \sq\medskip

\subsection{Roots of unity} We almost completed 
the first part of our program:

\noindent{\em imaginary integration $\Rightarrow$ 
sharp integration\hfill}

{\em \hfill $\Rightarrow$ Jackson summation 
$\Rightarrow$ Gaussian sums.}

\noindent Switching in (\ref{f32}) to the roots of unity,
we come to the following theorem.

\begin{theorem}\label{t32}
Let $q^{1/4}$ be a primitive $4N$-th root
of unity and  $k$ be an integer such that $0<k\le N/2.$ Then
\begin{equation}\label{f33}
\sum_{j=0}^{N-2k}q^{\frac{(k-j)^2}{4}}\frac{1-q^{j+k}}{1-q^k}
\prod_{l=1}^j\frac{1-q^{l+2k-1}}{1-q^l}=
\prod_{j=1}^k\frac{1}{1-q^j}\sum_{j=0}^{2N-1}
q^{\frac{(k-j)^2}{4}}.
\end{equation}
\end{theorem}

{\bf Proof.} By Galois theory, it is sufficient
to pick any primitive $q^{1/4},$ so we may assume that
$q^{1/4}=e^{\pi \imath/2N}$. Let us substitute
$\widehat q^{\,1/4}\,=e^{\pi \imath/2N-\pi/a}$ for $q$ 
in formula (\ref{f32}).

\begin{lemma}\label{l32}
(Siegel) Let $\, j_m=h+2Nm\, $
where $h, m\in \BZ$. Then asymptotically
as $a\to +\infty,$
$$\sum_{m=0}^\infty\widehat q^{\,(j_m-k)^2/4}\approx 
\frac{\sqrt{a}}{2N}\,q^{(h-k)^2/4}.
$$
\end{lemma}

{\bf Proof.} First,
$$\sum_{m=-\infty}^\infty\widehat q^{\,(j_m-k)^2/4}\,=q^{(h-k)^2/4}
\sum_{m=-\infty}^\infty e^{-\pi (h+2Nm-k)^2/4a}.$$
Second, the sum
$$
\Sigma=\sum_{m=-\infty}^\infty e^{-\pi (h+2Nm-k)^2/4a}=
\sum_{m=-\infty}^\infty e^{-\pi N^2(m+(h-k)/2N)^2/a}
$$
approximates the integral
$$
S=\int_{-\infty}^\infty e^{-\pi N^2(x+(h-k)/2N)^2/a}dx
$$
with the difference $S-\Sigma$  bounded as
$a\to +\infty.$
Third, 
$$
S\approx \int_{-\infty}^\infty e^{-\pi N^2x^2/a}dx=
\frac{\sqrt{a}}{N}
$$
asymptotically. \sq\medskip

Now we compare the asymptotics of the left-hand side and 
the right-hand side of formula (\ref{f32}) as
$\widehat q\to q,$  equivalently, $a\to +\infty.$ We represent
either side of this formula as the sum of $2N$ subsums
corresponding to all possible values $j\modd 2N.$ 

Due to the lemma, the right-hand side approaches
the right-hand side of formula (\ref{f33}) 
times $\sqrt{a}/N.$ 
To manage the left-hand side of formula
(\ref{f32}), we note that the product $\prod_{l=1}^j
\frac{1-q^{l+2k-1}}{1-q^l}$ is a periodic function of $j$
of period $N$ because $q^N=1.$ Moreover, this product vanishes
for $N-2k+1\le j<N$. Therefore the
left-hand side of formula (\ref{f32}) tends to
the left-hand side of formula (\ref{f33}) times
$\sqrt{a}/N$. Theorem \ref{t32} is proved. 
\sq\medskip

\subsection{Gaussian sums} We come to the following corollary.

\begin{corollary}\label{c31}
Taking $q^{1/4}=e^{\pi \imath/2N},$
$$\sum_{j=0}^{2N-1}q^{j^2/4}=(1+\imath)\sqrt{N}.$$
\end{corollary}

{\bf Proof.} Letting $k=[N/2]$ in formula (\ref{f33}),
$$
\sum_{j=0}^{2N-1}q^{j^2/4}=\left\{ \begin{array}{ccl}
q^{n^2/4}\prod_{j=1}^n(1-q^j)&\mbox{if}&N=2n,\\
q^{n^2/4}(1+q^{(2n-1)/4})\prod_{j=1}^n(1-q^j)
&\mbox{if}&N=2n+1.\end{array}\right.
$$

Only even $N=2n$ will be considered.
Setting $\Pi=\prod_{j=1}^n (1-q^{j}),$
we get $\Pi\bar{\Pi}=2N,$
since it is the value
of $(X^{N}-1)(X+1)(X-1)^{-1}$ at $X=1.$
Here the bar is the complex conjugation.
On the other hand, 
$$
\arg(1-e^{\imath\phi})=\phi/2-\pi/2 \hbox{\ for\ angles\ }
0< \phi< 2\pi,\hbox{\ and\ }
$$ 
$$\arg\Pi = \frac{\pi} {N} \frac{n(n+1)}{2}-
\frac{\pi n}{2} =
\frac{\pi (1-n)}{4}.\hbox{\sq}
$$

There is the following reduction of Theorem \ref{t32},
important from the viewpoint of applications to
the Gaussian sums. It has a direct relation to 
the "little" double affine Hecke algebra (see \cite{C7}).

\begin{theorem}\label{t32p}
Let $n=[N/2]$. Then for $0<k\le n,$ 
\begin{align}\label{fgausseven}
\sum_{j=0}^{n-k}q^{j^2-kj}\frac{1-q^{2j+k}}{1-q^k}
\prod_{l=1}^{2j}\frac{1-q^{l+2k-1}}{1-q^l}=\prod_{j=1}^k
\frac{1}{1-q^j}\sum_{j=0}^{N-1}q^{j^2-kj}.
\end{align}
\end{theorem}

{\bf Proof.} Considering the variant of (\ref{f33}) for  
$q^{1/4}\mapsto -q^{1/4},$ 
check that (\ref{fgausseven}) equals either the
half-sum of these two formulas for even $k$
or the half-difference for odd $k.$ 
\sq\medskip

\begin{corollary}\label{c32}
For even $N=2n$ and $k=n,$
\begin{equation}\label{f34}
\sum_{j=0}^{N-1}(-1)^jq^{j^2}=\prod_{j=1}^n(1-q^j),
\end{equation}
which can be rewritten as 
\begin{equation}\label{f35}
\sum_{j=0}^{n-1}(-1)^jq^{j^2}=\prod_{j=1}^{n-1}(1-q^j),
\end{equation}
due to the substitution $j\mapsto j+n.$
\end{corollary}

We are going to use these formulas for $q=\exp(\pi\imath m/n)$
to calculate the corresponding Legendre symbol
in terms of the integer parts $[mj/(2n)].$
Let $n>0$ be any integer and $m>0$ an odd integer, assuming
that they are relatively prime: $(m,n)=1.$ 
We set 
\begin{align}
&\bigl\{\frac {m}{n}\bigr\} 
\equal \imath^{(n-1)(m-1)/2}
(-1)^{\sum_{j=1}^{n-1}\,[\,\frac{mj}{2n}\,]}\label{flegang},\\
&G(m,n)\equal \sum_{j=1}^n
e^{\pi \imath\, j^2\,\frac{m}{n}\,+\pi \imath jm}.\label{fgaussmn}
\end{align}

Formula (\ref{fgaussmn}) is the
classical definition of the generalized Gaussian sum
(see e.g., \cite{Cha}). As we will see in the next theorem,
the first definition extends
the Legendre symbol $\bigl(\frac{m}{n}\bigr),$ 
The latter is $\pm 1$ as $m$ is a quadratic residue 
(non-residue) modulo $n,$ 
where  $n$ is an odd prime number or its power.

\begin{theorem}\label{t33}
(a) For odd $m$ coprime with $n,$
\begin{equation}\label{f36}
G(m,n)=\sqrt{n}\left( \frac{1+\imath}{\sqrt{2}}
\right)^{1-n}\, \bigl\{\frac {m}{n}\bigr\} \,.
\end{equation}

(b) Taking $n=p^a$ for odd prime $p$
and odd $a,\ $ 
$\bigl\lbrace \frac {m}{n}\bigr\rbrace =
\bigl(\frac {m}{n}\bigr).$ 
\end{theorem}

{\bf Proof.} 
The $G(m,n)$ is the sum from 
(\ref{f35}) for $q=e^{\pi \imath/n}.$
The product on the right hand-side 
of this equality was
calculated in Corollary \ref{c31} as $m=1.$ 
Indeed, $G(1,n)$ is given by (\ref{f36}).
 
(a) The equality $|G(m,n)|=\sqrt{n}$ is immediate:
$$G(m,n)\overline{G(m,n)}=\frac{1}{2}\prod_{j=1}^{N-1}(1-q^j)
=\frac{N}{2}=n.$$
Here $G(1,n)$ is sufficient to examine (apply the Galois
automorphisms). 
The formula for the argument of $G(m,n)$ is direct 
from formula (\ref{f35}).
Using $\arg(1-e^{i\phi})=\phi/2-\pi/2$ for $0<\phi <2\pi\,:$
$$
\arg(G(m,n))=\sum_{j=1}^{n-1}(\frac{\pi mj}{2n}-
\left[\frac{mj}{2n}\right]\pi -\frac{\pi}{2})=$$
$$=\frac{\pi m}{2n}\frac{n(n-1)}{2}-(n-1)\frac{\pi}{2}-\pi
\sum_{j=1}^{n-1}\left[\frac{mj}{2n}\right]=$$
$$=\frac{\pi}{4}(m-1)(n-1)-\frac{\pi}{4}(n-1)+
\sum_{j=1}^{n-1}\left[\frac{mj}{2n}\right]\modd 2\pi.$$

(b) For $m=1,$ the coincidence is evident. 
Generally, 
$$
G(m,n)=\pm G(1,n)=\bigl\lbrace \frac {m}{n}\bigr\rbrace G(1,n). 
$$
Using the left-hand side of (\ref{f35}),
the sign is plus if
$m$ is a quadratic residue modulo $n$ and
is constant on all non-residues $m.$
We use that $\BZ_{n}^*$ is cyclic. 
It cannot always be plus because it would 
give the invariance of $G(1,n)$ under
the  Galois automorphisms 
$q\mapsto q^m$ and would result in 
$G(1,n)\in \BQ.$
Therefore the sign is minus on the non-residues.  
\sq\medskip

\subsection{Etingof's theorem}
So far we considered the imaginary integration,
the Jackson summation, and its variant at roots of unity.
However the most natural choice is of course
$\int_{i \eps +\BR}q^{x^2}\delta_k(x)dx$.
The calculation of the latter integral was 
performed by P.~Etingof. 
We are grateful for the permission to include his
note in the paper. Here we use $i$ (instead of $\imath$)
for the imaginary unit.

Recall that $q=\exp{-1/a}$, $\omega=2\pi ia.$
Let $k$ be a positive
real number, $\epsilon>0$ a small positive number.
We will use the function
$\Pi_k,$  the right-hand side of (\ref{f31}), and
its properties, including information about the zeros.
Let
$$
\Psi_k=\frac{1}{g_k}\int_{i\epsilon+\BR}q^{x^2}\delta_k(x)dx.
$$
It is clear that $\Psi_k$ is well defined.

Our goal is to calculate $\Psi_k,$ which will be done in 
Theorem \ref{mainth7} below. 
To formulate the theorem, we need some definitions.

First, the theta function is
$$
\theta(k,\tau)=\sum_{n\in \BZ}(-1)^n
e^{2\pi ink+\pi in(n-1)\tau},\ \Ima(\tau)>0.
$$
It is a periodic entire function with period $1$,
which satisfies the equation
$$
\theta(k+\tau,\tau)=-e^{-2\pi ik}\theta(k,\tau).
$$
It is defined by this equation uniquely up to scaling.
The zeros of $\theta(k,\tau)$ with respect to $k$ are
$m+n\tau$, $m,n\in \BZ$ (all of them are simple).

Define the following (degenerate) Appell function:
$$
A(k,\tau)=2\pi i\sum_{n\in \BZ\setminus 0}
\frac{e^{2\pi in^2\tau+2\pi
    ink}}{e^{2\pi in\tau}-1}.
$$
It is not expressed via theta functions.

\begin{theorem}\label{t41}
The function $\Psi_k$ extends to
an entire function, and one has
$$
\Psi_k=2\sqrt{\pi
    a}\cdot e^{-\pi i k}q^{-k(k+1)/2}
\frac{\theta(k,\omega)}{\theta'(0,\omega)}\times$$
$$\times \biggl(\frac{\theta'}{\theta}(k,\omega)-\frac{\theta'}
{\theta}(1/2,\omega)-A(k,\omega)+A(1/2,\omega)\biggr).
$$
\end{theorem}
We present the proof as a chain of lemmas.

\begin{lemma}\label{l41}
The function $\Psi_k$ has an analytic continuation
to the region $\Rea(k)>0$, $|\Ima(k)|<|\omega|$.
\end{lemma}
{\bf Proof.} The contour of integration can be replaced by
$$
(-\infty+i\epsilon,i\epsilon]\cup [i\epsilon,\frac{1}{2}
\omega-i\epsilon]\cup [\frac{1}{2}\omega-i\epsilon,
\frac{1}{2}\omega-i\epsilon+\infty),
$$
or by
$$
(\frac{1}{2}\omega-i\epsilon+\infty,\frac{1}{2}\omega-
i\epsilon]\cup [\frac{1}{2}\omega-i\epsilon,i\epsilon]
\cup
[i\epsilon,\infty+i\epsilon).
$$
This implies the statement. \sq\medskip

\begin{lemma}\label{l42}
$$
\Psi_{k+\omega}-\Psi_k=-2i\Psi_k^{\#}=-2i\Pi_k.
$$
In particular, $\Psi_k$ extends to a holomorphic function
in the half-plane $\Rea(k)>0$.
\end{lemma}
{\bf Proof.} This follows directly from Cauchy's residue
formula. \sq\medskip

\begin{lemma}\label{l43}
One has
$$
\Psi_k=-q^{k+1}\Psi_{k+1}.
$$
\end{lemma}
{\bf Proof.} This is proved by using the shift--formula
in the same way as the relation (\ref{fpsisharp}).
Actually the sharp integration is a variant of the real
integration so they result
in the same multiplier. However their behavior in the
imaginary direction is different. 
\sq\medskip

\begin{lemma}\label{l44}
$\Psi_k$ extends to an entire function.
\end{lemma}
{\bf Proof.} This is immediate from Lemma \ref{l43}. \sq\medskip

Let us now define the entire function
$$
F_k=q^{k(k+1)/2}e^{\pi ik}\Psi_k.
$$

\begin{lemma}\label{l45}
$F_k$ is periodic with the period $1$, and
\begin{equation}\label{f41}
e^{2\pi ik}F_{k+\omega}+F_k=
2ie^{\pi ik}q^{k(k+1)/2}\Pi_k,
\end{equation}
where $\Pi_k$ is the right-hand side of (\ref{f31}).
\end{lemma}
{\bf Proof.} The real periodicity
is direct from Lemma \ref{l43}. 
Equation (\ref{f41}) is checked following the same
lines as formula (\ref{ffkgk}).
\sq

Now, we are going to construct at least one
holomorphic solution of equation
(\ref{f41}) and then correct it by adding a solution of the
homogeneous equation, i.e., a multiple of $\theta(k,\omega)$.

To construct such a solution,
let us express the right-hand side of (\ref{f41}) in terms
of the theta function.

\begin{lemma}\label{l46}
One has
$$
ie^{\pi ik}q^{k(k+1)/2}\Pi_k=C\,\theta(k+\omega+1/2,2\omega)\,
\frac{\theta(k,\omega)}{\theta'(0,\omega)},
$$
where $C$ is a constant.
\end{lemma}

{\bf Proof} is a combination of the
translational properties  of
$\Pi_k$ and the information about its zeros.\sq\medskip

Now consider the following two equations:
\begin{equation}\label{f42}
G_{k+\omega}-G_k=1,\hbox{\ and\ }
\end{equation}
\begin{equation}\label{f43}
H_{k+\omega}-H_k=\theta(k+\omega+\frac{1}{2},2\omega)-1.
\end{equation}

It is clear
from Lemma \ref{l46} that if $G_k$ solves (\ref{f42}) and $H_k$
solves (\ref{f43}) then
$-2C\frac{\theta(k,\omega)}{\theta'(0,\omega)}(G_k+H_k)$
solves (\ref{f41}).
Concerning (\ref{f42}), it is satisfied by 
$$
G_k=-\frac{1}{2\pi i}\frac{\theta'}{\theta}(k,\omega).
$$
On the other hand,
$$
H_k=\frac{1}{2\pi i}A(k,\omega)
$$
solves (\ref{f43}). This implies that
$$
F_k=-\frac{C}{\pi i}\frac{\theta(k,\omega)}{\theta'(0,\omega)}
(A(k,\omega)-\frac{\theta'}{\theta}(k,\omega)-\beta)
$$
for constants $C,\beta$ (note that the right-hand
side is an entire function).
So we need to find these constants.

\begin{lemma}\label{l47}
The function $F_k$ satisfies the conditions
$F_0=2\sqrt{\pi a}$, $F_{1/2}=0$.
\end{lemma}

{\bf Proof.} The first statement
is clear (use the Gauss integral);
the second one follows from the fact that $g_k$ has a pole
at $k=-1/2$ and $F_k$ is $1$-periodic.
\sq\medskip

Now we see that
$\beta=A(1/2,\omega)-\frac{\theta'}{\theta}(1/2,\omega)$,
and $C=2\pi i\sqrt{\pi a}$. This completes
the proof of the theorem.
\medskip

{\bf Comment.}
The function $A(k,\omega)$ has the following
connection to Appell functions (see e.g., \cite{Po}
and references therein). The Appell function
$\kappa(u,k,\tau)$
is defined as a unique holomorphic, $1$-periodic in $k$
solution of the equation
$$
\kappa(u,k+\tau,\tau)=e^{2\pi iu}\kappa(u,k,\tau)+
\theta(k+\frac{1+\tau}{2},\tau).
$$
It has the Fourier series expansion
$$
\kappa(u,k,\tau)=\sum_{n\in \BZ}\frac{e^{\pi in^2\tau+2\pi
    ink}}{e^{2\pi in\tau}-e^{2\pi iu}}.
$$
The Appell function has a pole at $u=0,$ so one can
introduce its regular part:
$$
\kappa_0(k,\tau)=\sum_{n\in \BZ\setminus 0}
\frac{e^{\pi in^2\tau+2\pi
    ink}}{e^{2\pi in\tau}-1}.
$$
Then
$$
A(k,\omega)=2\pi i(\kappa_0(k,2\omega)+
\kappa_0(k+\omega,2\omega)).
$$

\bigskip
\setcounter{equation}{0}
\section{Nonsymmetric Hankel transform}
We go to the second level of our program: 
{\it to connect
the Hankel transform and the Fourier transform
on $\BZ_N$ in one theory.} 

In the first place, we will interpret the integral
$$\int_{\BR}e^{-x^2}|x|^{2k}dx=\Gamma (k+1/2),\; \Rea k>-1/2,$$
as the structural constant of the classical
Hankel transform, and prove Theorem \ref{t12}. Our approach
is different from the usual treatment of similar integrals 
in classical works on Bessel functions.

\subsection{Operator approach}
Recall that $\phk_\lambda (x)\equal 
\phk (\lambda x)$ is given as follows:
$$\phk (t)=\sum_{n=0}^\infty \frac{t^{2n}\Gamma (k+1/2)}{n!
\Gamma (k+n+1/2)},
\; k\not \in -1/2-\BZ_+.$$
It is even in both $x$ and $\lambda\,:$ 
$\phk(t)=\phk(-t).$ The connection with the classical
Bessel function is established in (\ref{fbesselj}).
 
Our first step is to reprove the classical formula
(see e.g., \cite{Luk}) for the Gauss integral in the presence of
the Bessel functions, using the operator interpretation
of the latter.

\begin{theorem}\label{t51}
For arbitrary complex $\lambda,\mu$ and $\Rea(k)>-1/2,$
\begin{equation}\label{f51}
\int_{\BR}\phk_\lambda (x) \phk_\mu (x)e^{-x^2}|x|^{2k}dx=
\Gamma (k+1/2)
\phk_\lambda (\mu)e^{\lambda^2+\mu^2}.
\end{equation}
\end{theorem}

{\bf Proof} is based on the theory of the differential operator 
$$L=\frac{d^2}{dx^2}+\frac{2k}{x}\frac{d}{dx}.
$$

\begin{lemma}\label{l51}
(i) The function $\phk_\lambda (x)$ is a unique  even solution
of the eigenvalue problem $L\phi(x)=
4\lambda^2\phi (x)$ with the normalization $\phi(0)=1.$

(ii) The operator $L$ is selfadjoint with
respect to the scalar product $\langle f,g\rangle =
\int_{\BR}f(x)g(x)|x|^{2k}dx.$
\end{lemma}

{\bf Proof.} Both statements are
straightforward. Concerning (ii), we check that
$L=|x|^{-k}\circ H\circ |x|^k,$ where $H=\frac{d^2}{dx^2}+
\frac{k(1-k)}{x^2},$ and use that
$H$ is selfadjoint for the scalar
product $\int_{\BR}f(x)g(x)dx.$ \sq\medskip

\begin{corollary}\label{c51}
Asymptotically as $|x|\to \infty,$
$$
\phk_\lambda (x)\sim C(\lambda)(e^{2\lambda x}+
e^{-2\lambda x})|x|^{-k}
$$
for a constant $C(\lambda)$.
\end{corollary}

{\bf Proof.} The eigenfunctions of $H$ with the eigenvalue 
$4\lambda^2$ are asymptotically 
$e^{2\lambda x}$ and $e^{-2\lambda x}$.\sq\medskip

Let $f(x)$ be an even function of $x\in \BR$
such that $f(x)e^{cx}\to 0$ as $x\to \infty$ for any
$c.$ The {\em Hankel transform} is defined as follows:
\begin{equation}\label{fhankel}
\bF(f)(\lambda )\equal \frac{1}{\Gamma (k+1/2)}\int_{\BR}
f(x)\phk_\lambda (x)|x|^{2k}dx,\; \Rea k>-1/2.
\end{equation}
Let $\bF^{op}$ be the corresponding transform
on the operators: $\bF^{op}(A)=\bF\,\bf A\, \bF^{-1}.$ 

\begin{lemma}\label{l53}
(a) $\bF^{op} (L)=4\lambda^2$,

(b) $\bF^{op} (4x^2)=L_\lambda\equal \frac{d^2}{d\lambda^2}+
\frac{2k}{\lambda}\frac{d}{d\lambda}$,

(c) $\bF^{op} (4x\frac{d}{dx})=-4\lambda
\frac{d}{d\lambda}-4-8k$.
\end{lemma}

{\bf Proof.} Formulas (a) and (b) hold
because the operator $L$ is selfadjoint and 
$\phk_\lambda (x)$ is its eigenfunction. 
Formula (c) formally follows from (a), (b), and the identity
$$[L,x^2]=4x\frac{d}{dx}+2+4k.\hbox{\ \sq}$$ 

We set $\phk_{\pm,\lambda}= \phk_\lambda (x)\, e^{\pm x^2},$
and
$$
L^+= e^{x^2}\circ L\circ e^{-x^2},\
L^-= e^{-x^2}\circ L\circ e^{x^2}.
$$ 
Explicitly,
$$L^+=L+4x\frac{d}{dx}+2+4k+4x^2,\; L^-=L-4x
\frac{d}{dx}-2-4k+4x^2.$$
Using Lemma \ref{l53},
$$\bF^{op} (L^-)=L^+.$$
It is immediate from 
Lemma \ref{l51}, (i)
that $L^\pm \phk_{\pm,\lambda}=4\lambda^2 \phk_{\pm,\lambda}$
and
$$\bF(\phk_{-,\mu})=C_\mu e^{\mu^2}\phk_{+,\mu}(\lambda),$$
where $C_\mu$ does not depend on $\lambda.$ Therefore
$$\int_{\BR}\phk_\lambda (x)\phk_\mu (x)e^{-x^2}|x|^{2k}dx=
\Gamma(k+1/2)C_\mu \phk_\mu (\lambda)e^{\lambda^2+\mu^2}.$$
The left-hand side of this equality is invariant 
under the change $\lambda \leftrightarrow \mu,$  
so is the right-hand side. Hence $C_{\mu}=C_0=1.$ 
The theorem is proved. \sq\medskip

\subsection{Nonsymmetric theory} 
The Hankel transform sends even functions 
to even functions (and is zero when applied to odd
functions), by construction.
In this section we consider its nonsymmetric version.
We denote the reflection $f(x)\mapsto f(-x)$ by $s.$
\medskip

{\bf Key Definition} (C.~Dunkl) $D =\frac{d}{dx}-\frac{k}{x}(s-1).$
\sq\medskip

The operator $\cD$ is obviously odd, i.e.,
$s\cD s=-\cD.$
The restriction of the even operator
$\cD^2$ to the space of even functions
coincides with $L=\frac{d^2}{dx^2}+\frac{2k}{x}\frac{d}{dx}:$
$$\cD^2(f)=\cD (f')=f''+\frac{2k}{x}f'\hbox{\ for\ even\ } f(x).$$

The following two lemmas are counterparts of Lemmas
\ref{l51} and \ref{l53}. 

\begin{lemma}\label{l51p}
(i) If $\lambda \ne 0$ or $\lambda=0$ and
$k\not \in -1/2-\BZ_+,$  then the eigenvalue problem
$\cD \psi = 2\lambda \psi$
has a unique analytic at $0$ solution $\psk_\lambda(x)$
with the normalization $\psk_\lambda (0)=1.$ 
Moreover, $\psk_\lambda(x)=\psk (\lambda x)$ and
$$
\psk(t)=\phk(t)+\frac{(\phk(t))'}{2} \, \hbox{\ for\ \,}
\phk(t)=\frac{\psk(t)+\psk(-t)}{2}.
$$
\noindent
If $\lambda =0$ and $k=-1/2-n,\; n\in \BZ_+,$ 
then $\psi =1+Cx^{2n+1}$ for $C\in \BC.$

(ii) Let $\cD^*$ be the adjoint of
$\cD$ with respect to the scalar product $\int_{\BR}
f(x)g(x)|x|^{2k}dx$. Then $\cD^*=-\cD$.
\end{lemma}
{\bf Comment.} For any $k,$ there is an extra solution
$\psi =x/|x|^{2k+1}$ of the equation $\cD \psi =0,$
but it is analytic at $0$ only as $k\in -\frac{1}{2}-\BZ_+.$
\medskip

{\bf Proof.} (i) We set $\psi =\psi_0+\psi_1$
where $\psi_0$ is even and $\psi_1$ is odd. Then 
$\cD \psi =2\lambda \psi$ is equivalent to the system of
equations
$$\psi_0'=2\lambda\,\psi_1,\ 
\psi_1'+\frac{2k}{x}\psi_1=2\lambda\psi_0. 
$$
We get that $\psi_0$ satisfies Lemma \ref{l51}, (i).

(ii) The operator $|x|^k\circ \cD \circ |x|^{-k}=\frac{d}{dx}-\frac{k}{x}s$ is
anti-selfadjoint with respect to
the scalar product $\int_{\BR}f(x)g(x)dx$. 
\sq\medskip

Following Dunkl,
we introduce the {\em nonsymmetric Hankel transform}:
$$\bF(f)(\lambda)=\frac{1}{\Gamma (k+1/2)}
\int_{\BR}f(x)\psk_\lambda (x)|x|^{2k}dx;$$
$\bF^{op}$ is its action on operators.

\begin{lemma}\label{l53p}
(a) $\bF^{op} (\cD )=-2\lambda$,

(b) $\bF^{op} (s)=s_\lambda,\ 
s_\lambda f(\lambda)\mapsto f(-\lambda),$

(c) $\bF^{op} (2x)=\cD_\lambda\equal \frac{d}{d\lambda}-
\frac{k}{\lambda}(s_\lambda-1)$. \sq\medskip
\end{lemma}

The main result of this section is the following theorem.
\begin{theorem}\label{t52}
(Master formula) 
\begin{equation}\label{f52}
\int_{\BR}\psk_\lambda (x)\psk_\mu (x)e^{-x^2}|x|^{2k}dx=
\Gamma (k+1/2)\psk_\lambda(\mu)e^{\lambda^2+\mu^2}.
\end{equation}
\end{theorem}

{\bf Proof} is parallel to the proof of 
Theorem \ref{t51}. However the necessary calculations
are simpler than in the symmetric case.
Introduce functions $\psk_{\pm,\lambda}(x)=\psk_\lambda(x)
e^{\pm x^2}$ and operators $\cD^+=e^{x^2}\circ \cD \circ
e^{-x^2}$, $\cD^-=e^{-x^2}\circ \cD \circ e^{x^2}$. Easy
calculations show that 
$$\cD^\pm \psk_{\pm,\lambda}=
2\lambda \psk_{\pm,\lambda},\ \cD^\pm=\cD \mp 2x,\
\bF^{op} (\cD^-)=\cD^+_\lambda.
$$
Hence 
\begin{equation}\label{f53}
\bF(\psk_{-,\mu})=C_\mu e^{\mu^2}\psk_{+,\mu}(\lambda),
 \end{equation}
where $C_\mu$ does not depend on
$\lambda$. Finally, $C_\mu=C_0=1$ thanks to the
symmetry $\lambda \leftrightarrow \mu$. \sq\medskip

{\bf Comment.} Theorem \ref{t52} is equivalent to 
Theorem \ref{t51}, which is a special feature of the
one-dimensional setup. Generally, the nonsymmetric formula
results in the symmetric one but not the other way round.

In the first place, 
Theorem \ref{t52} implies Theorem \ref{t51}, since $$\phk_\lambda(x)=(\psk_\lambda)_0(x)=
(\psk_\lambda(x)+\psk_{-\lambda}(x))/2,
\hbox{\ and\ for\ } \mu.
$$

To deduce Theorem \ref{t52} from Theorem \ref{t51} we 
need to show that
$$\int_{\BR}(\psk_\lambda)_1(x)\,(\psk_\mu)_1(x)\,e^{-x^2}
|x|^{2k}dx=\Gamma (k+1/2)(\psk_\lambda)_1(\mu)\,
e^{\lambda^2+\mu^2},$$
where $(\psk_\lambda)_1(x)$ is the odd component of
$(\psk_\lambda)(x).$ We may ignore the "cross-terms"
$(\psk_\lambda)_1(x)(\psk_\mu)_0(x)$ since they
are odd and their integrals are zero. 
To calculate the integral above, we can
use the same Theorem \ref{t51} because of the following
{\em shift--formula}:
$$\frac{1}{x}(\psk_\lambda)_1=\frac{2\lambda}{1+2k}
\phi^{(k+1)}_\lambda.
$$

This equivalence somewhat clarifies why the nonsymmetric 
Hankel transform did not appear (as far as we know) in
classical works on Bessel functions. It adds nothing new
to the symmetric (even) one.
\sq\medskip

Using the operator $e^{-\cD^2/4},$ we can rewrite
formula (\ref{f52}) in the following form:
\begin{equation}\label{f54}
\int_{\BR}\psk_\lambda(x)e^{-\cD^2/4}(f(x))e^{-x^2}
|x|^{2k}dx=\Gamma(k+1/2)f(\lambda)e^{\lambda^2},
\end{equation}
where $f(x)$ is a function from a suitable completion of
the space linearly generated by the functions $\psk_\mu(x).$ 
Indeed,
Theorem \ref{t52} shows that (\ref{f54}) holds for $f(x)=
\psk_\mu(x)$. Formula (\ref{f54}) leads to the following
entirely algebraic definition of
the Hankel transform: 
$$\bF=e^{x^2}\circ
e^{\cD^2/4}\circ e^{x^2}.$$
We will use it later to introduce the 
truncated Hankel transform.

Since $\cD$ is nilpotent,  (\ref{f54}) results in the
following important corollaries.

\begin{corollary}\label{c52}
(cf. Corollary \ref{c12}) The
nonsymmetric Hankel transform restricts to a map $\bF:
\BC [x]e^{-x^2}\to \BC [\lambda]e^{\lambda^2}$ and, moreover,
preserves the filtration by degrees of polynomials.\sq\medskip
\end{corollary}

We set $\ \cB_x=\sum_\mu \BC\,\psk_\mu(x)e^{-x^2},\ $
$\cB_\lambda=\sum_\mu \BC\,\psk_\mu(\lambda)e^{+\lambda^2}.$

\begin{corollary}\label{c53} (Inversion)
Introducing the imaginary Hankel transform
$$\bF_{im}(g)(x)\equal \frac{1}{i\Gamma(k+1/2)}\int_{i\BR}
g(\lambda)\psk_x(-\lambda)|\lambda|^{2k}d\lambda,$$
\noindent $\bF\circ \bF_{im}=id$ and  $\bF_{im}\circ \bF=id$
in the spaces $\cB_x,\cB_\lambda,$ or
$\, \BC [x]e^{-x^2},$ $\,\BC [\lambda]e^{\lambda^2},$
or their suitable completions.
\end{corollary}

{\bf Proof.} This follows from the formula
\begin{equation}\label{fhankim}
\frac{1}{i}\int_{i\BR}\psk_\lambda(-x)\psk_\mu(x)
e^{x^2}|x|^{2k}dx=\Gamma(k+1/2)\psk_\lambda(\mu)
e^{-\lambda^2-\mu^2},
\end{equation}
which results from (\ref{f52}) upon $x\mapsto ix, \lambda
\mapsto i\lambda, \mu \mapsto -i\mu.$ 
The passage to $\BC[x]e^{-x^2}$ and $\BC [\lambda]e^{\lambda^2}$
is either by means of (\ref{f54}) or
via the completion.
\sq\medskip

\begin{corollary}\label{c54}
(Plancherel formula) For functions $f(x), g(x)$ in
$\cB_x\,$ or in $\,\BC[x]e^{-x^2},\,$ we set
$\widehat f=\bF(f),\; \widehat g=\bF(g).$ Then
\begin{equation}\label{f55}
\int_{\BR}f(x)g(x)|x|^{2k}dx=\frac{1}{i}\int_{i\BR}
\widehat f(-\lambda)\widehat g(\lambda)|\lambda|^{2k}d\lambda.
\end{equation}
\end{corollary}

{\bf Proof.} The left-hand side
of (\ref{f55}) can be readily
calculated for
$f(x)=\psk_{\mu_1}(x)e^{-x^2}$ and 
$g(x)=\psk_{\mu_2}(x)e^{-x^2}$
thanks to formula (\ref{f53}). We substitute $\sqrt{2}x\mapsto x$
and use that $\psk_\mu(x)$ 
depends on the product $\mu x.$ Similarly, we calculate 
the right-hand side using (\ref{fhankim}).
Completing or using (\ref{f54}), 
we switch to the space $\BC[x]e^{-x^2}.$ 
\sq\medskip

Concerning completions, 
when $k\in \BR$ the last corollary produces
the nonsymmetric Hankel transform on the space 
$L^2(\BR,|x|^{2k})$ because it is the 
$L^2$-completion of the space
$\BC[x]e^{-x^2}$ with respect to the scalar
product $\int_{\BR}f(x)\overline{g(x)}|x|^{2k}dx.$

To be more exact, we first use this corollary for real-valued
functions. The image of the corresponding $L^2$
will be the $\BR$-subspace of 
$L^2(\imath\BR,|\lambda|^{2k})$ formed by
functions satisfying $\overline{f(\lambda)}=f(-\lambda).$
Then we extend the coefficients from $\BR$ to $\BC.$

\subsection{Double H double prime} We denote it by
$\HH''$ ($\HH'$ is reserved for the trigonometric
limit which leads to the Harish-Chandra theory). It is
$$\HH''=\langle \cD, x, s\rangle /\{ s\cD s=-\cD, sxs=-x,
[\cD, x]=1+2ks, s^2=1\} .$$
The algebra $\HH''$ acts on the vector space $\BC [x]$ via
the Dunkl operator:
$$s(f(x))=f(-x),\; x(f(x))=xf(x),\; \cD(f)=
(\frac{d}{dx}-\frac{k}{x}(s-1))f.$$
We will call it the {\em polynomial representation}.

\begin{theorem}\label{t53}
(a) The $\HH''$-module $\BC[x]$ is faithful.

(b) (PBW-property) The elements $x^n\cD^ms^\eps,\; n,m\in
\BZ_+,\; \eps=0,1,$ form a basis of $\HH''$.
\end{theorem}

{\bf Proof.} It is clear that any element of $\HH''$ can be
expressed in the following form:
$$H=\sum_{n,m,\eps}c_{n,m,\eps}x^n\cD^ms^\eps, \;
c_{n,m,\eps}\in \BC.$$
If  $c_{n,m,\eps}\ne 0$ at least once, then
the image of $H$ is nonzero in
the space of either even or odd functions. 
This gives both (a) and (b). \sq\medskip

Because the polynomial representation is always
faithful we will identify $\cD$ considered as a 
generator of the double H double prime with 
its image, the Dunkl operator.

The algebra $\HH''$ has an automorphism $\omega$ defined by
\begin{equation}\label{fomegap}
\omega (\cD)=-2x,\; \omega (2x)=\cD, \; \omega(s)=s.
\end{equation}
It can be represented as follows:
\begin{equation}\label{fbraid}
\omega=e^{x^2}\circ e^{\cD^2/4}\circ^{x^2}=
e^{\cD^2/4}\circ e^{x^2}\circ e^{\cD^2/4}.
\end{equation}
Here we extend $\HH''$ by adding $e^{x^2},\  e^{\cD^2/4},$
and treat the latter as inner automorphisms in the resulting
greater algebra $\widehat{\HH}''.$
It is simple to see that both preserve $\HH''.$
The automorphism $\omega$ is nothing else but
an algebraic version of
the Hankel transform $\bF.$ 
 
The coincidence of the two representations of $\omega$ 
can be deduced from the defining relations or (much simpler)
checked in the polynomial representation. It defines the
action of the {\em projective} $PSL(2,\BZ)$ (due to Steinberg)
on $\HH''.$ 

The polynomial representation is always faithful but not
always irreducible.

\begin{theorem}\label{t54}
The $\HH''$-module $\BC[x]$ is irreducible if
and only if $k\not \in -1/2-\BZ_+$.
\end{theorem}

{\bf Proof.} If $k\not \in -1/2-\BZ_+,$ then the equation
$\cD \psi=0$ has a unique solution $\psi=1$ in $\BC[x]$. Since
the operator $\cD$ is nilpotent in $\BC[x],$ any submodule should
contain $1$ and therefore the whole $\BC[x]$. If $k=-1/2-n$ then
$x^{2n+1}$ generates a nontrivial submodule of $\BC[x]$ since
$\cD (x^{2n+1})=0$.
 \sq\medskip

We see that $\HH''$ can have finite dimensional
irreducible  representations
for some $k.$ It is not difficult to describe them all.
Generally, the theory of finite dimensional representations
of double Hecke algebras associated with root
systems is far from being complete.

\subsection{Finite dimensional representations}
\newcommand{\Pc}{{\BC[x]}}
\newcommand{\zit}{{\BZ}}
\def\dx{ \frac{d}{ dx} }
\def\rit{\BR} 
\def\nit{\BN} 
\def\qit{\BQ} 
\def\cit{\BC} 
\def\Dc{\cD} 
We will use that 
\begin{align}\label{hcom}
 & [h,x]=x,\ [h,\cD]=-\cD\hbox{\ \ for\ }
h=(x\cD+\cD x)/2.
\end{align}

Note that $h$
is $x\dx+k+1/2$ in the polynomial
representation. Since it is faithful,
(\ref{hcom})  is the claim that  
$x,\Dc$ are homogeneous operators 
of degrees $\pm 1,$ which is obvious.

These relations are the defining relations of $osp(2|1),$
but we will not rely on the theory of this super Lie algebra. 
For our purpose, a reduction to $sl_2$ is sufficient.

Namely, we will use that the elements $e=x^2,$  
$f=-\cD^2/4,$ and $h$ 
satisfy the defining relations
of $sl_2(\cit).$ Indeed, $[e,f]=h$ because
$$
[\cD^2,x^2]=[\cD^2,x]x+x[\cD^2,x]=
2\cD x+x (2\cD),
$$
and the relations $[h,e]=2e,[h,f]=-2f$ readily result 
from (\ref{hcom}).
The Casimir operator
$C=h^2-2h+4ef$ becomes  
$$
(h^2-2h)(1)=h(h-2)(1)=(k+1/2)(k-3/2)\hbox {\ in\ }\Pc.
$$

The module $\BC[x]$ is the Verma module with
the $h$-lowest weight $1/2+k$ over $U(osp(2|1)).$ 
It is the direct sum of the
two Verma modules with $h$-lowest weights
$1/2+k$ and $3/2+k,$ formed by even and odd functions respectively,
with respect to the action of $U(sl_2).$

Let $k=-n-\frac12$ for 
$n\in\zit_+.$ Then  $V_{2n+1}=\BC[x]/(x^{2n+1})$ is
an irreducible representation of $\HH''$. 
The elements of $V_{2n+1}$ can be identified with polynomials
of degree smaller than  $2n+1$.

\begin{theorem} 
\label{tfdimpr}
Finite dimensional representations of 
$\HH''$ exist only as $k=-n-1/2$ or $k=n+1/2$ for
$n\in \zit_+.$ 
Given such $k,$ the algebra $\HH''$ has a unique
finite dimensional irreducible representation
up to isomorphisms. It is either $V_{2n+1}$ for negative $k,$ 
or for its image under the $\HH''$-automorphism: 
\begin{equation}\label{fkaut}
x\mapsto x, \cD \mapsto \cD, \ s\mapsto -s, k\mapsto -k,
\end{equation}
in the case of positive $k.$
\end{theorem}
{\bf Proof.} 
Let $V$ be an irreducible finite dimensional representation
of $\HH''.$ 
Then the subspaces $V^0,V^1$ of $V$ 
formed respectively by $s$-invariant and $s$-anti-invariant vectors
are preserved by $e,$  $f,$ and $h.$
Indeed, $sxs=-x, s\cD s=-\cD,$ and $sh=hs.$
Note that $s$ leaves all $h$-eigenspaces invariant.
For instance, it commutes with the $sl_2$-action.
One gets
$$
\cD x =h+k+1/2,\ x\cD =h-k-1/2 \hbox{\ \ in\ } V^0,
$$
and the other way round in $V^1.$

Let us check that $\pm k\in -1/2- \zit_+.$
All $h$-eigenvalues in $V$ are integers 
thanks to the general theory of finite dimensional 
representations of $sl_2(\cit)$.

We pick a highest vector $v$, i.e., a
nonzero $h$-eigenvector $v\in V$ 
with the maximal possible eigenvalue $m.$
Using the automorphism (\ref{fkaut}), we will assume now
and later that it belongs to $V^0.$
Then $m\in \zit_+$ (the theory of $sl_2$) and 
$x(v)=0$ because the latter is an $h$-eigenvector with
the eigenvalue $m+1.$ Hence $\cD x (v)=0,$ 
$m+k+1/2=0,$ and $k=-1/2-m.$   

Let $U^0$ be a nonzero irreducible 
$sl_2(\cit)$-submodule of $V^0.$ The spectrum of $h$ in $U^0$ is 
$\{\,-n,-n+2,\ldots,n-2,n\,\}$ for an integer $n\ge 0.$
Let $v_l\neq 0$ be an $h$-eigenvector with the eigenvalue 
$l.$ If $e(v)=0,$ then $v=cv_{n}$
for a constant $c,$ and if $f(v)=0$ then
$v=cv_{-n}.$  

Let us check that $\cD x (v_{n})=0,$ 
$x\cD (v_{-n})=0,$ and  
$$
\cD x(v_l)\neq 0\hbox{\ for\ } l\neq n, 
\ \  x\cD(v_l)\neq 0 \hbox{\ for\ } l\neq -n.
$$ 

Both operators, $\cD x$ and $x\cD,$ obviously
preserve $U^0:$ 
$$
\cD x(v_l)=(l+k+1/2)v_l,\ x\cD (v_l)=(l-k-1/2)v_l.  
$$
Hence, 
$$
\cD^2 x^2(v_l)=
((\cD x)^2+(1-2k)(\cD x))(v_l)=
(l+k+1/2)(l-k+3/2)v_l.
$$
Setting $l=n,$ we get that $(n+k+1/2)(n-k+3/2)=0$
and $k=-1/2-n,$ because 
$k<0$ and  $n-k+3/2>0.$  Thus $\cD x (v_n)=0.$
The case of $x\cD$ is analogous.

The next claim is that $x(v_n)=0,\ \cD(v_{-n})=0.$ Indeed, 
$x(v')=0$ and $\cD(v')=0$ for  $v'=x(v_n).$ Therefore 
$$
0=[\cD,x](v')=
(1+2ks)(v')=(1-1+n)v'=nv'.
$$
This means that either $v'=0$ or $n=0.$ In the latter case,
$v'$ is proportional to $v_0$ and therefore $v'=x(v_0)=0$
as well. 
Similarly, $\cD(v_{-n})=0.$ 

Now we use the formula 
$$
\cD (x^2(v_l))=x(2+x\cD)(v_l)=
(2+l-k-1/2)x(v_l)=(2+l+n)x(v_l),
$$
and get that $x(v_l)\in \cD(U^0)$ for any $-n\le l\le n.$
Hence $U=U^0+\cD(U^0)$ 
is $x$-invariant. It is obviously $\cD$-invariant
and $s$-invariant ($\Leftarrow$ $\cD(V^0)\subset V^1$).
Also the sum is direct.

Finally, $U$ is an $\HH''\,$-module and 
has to coincide with $V$ because the latter was assumed to be
irreducible. The above formulas make
the $\HH''$- isomorphism $U\simeq V_{2n+1}$ explicit: 
the $h$-eigenvectors $x^i(v_{-n})\in U$ 
go to the monomials $x^i\in V_{2n+1}.$
\sq\medskip

\subsection{Truncated Hankel transform}
Let us consider the representation $V_{2n+1}$ closely.
We endow it with the following two scalar products:
$$\langle f,g\rangle_+\equal 
\Res(f(x)g(x)x^{-2n-1})$$
and $$\langle f,g\rangle_-\equal 
\Res(f(x)g(-x)x^{-2n-1}),
$$
analogous to the standard scalar products considered above.
Here by $\Res,$ we mean the coefficient of $x^{-1}.$
One has:
$\langle sf,g\rangle_\pm=\langle f, sg\rangle_\pm,$ and
\begin{equation}\label{fscalar}
\langle x f,g\rangle_\pm=\pm\langle f, x g\rangle_\pm,\
\langle \cD f,g\rangle_\pm=\mp\langle f,\cD g\rangle_\pm.
\end{equation}

Note that if $f(x)=\sum_{i=0}^{2n}a_ix^i$ then
$$\langle f,f\rangle_\pm=
\sum_{i=0}^{2n}(\pm)^ia_ia_{2n-i},$$ 
so these
scalar products are far from being positive definite.

In contrast to $\Pc,$ the Gaussians are 
well defined in $V_{2n+1}:$
$$e^{\pm x^2}=\sum_{m=0}^{2n}(\pm x^2)^m/m!. $$

We may introduce the truncated Hankel transform on this space
by the old formula
$\bF_+=e^{x^2}e^{\cD^2/4}e^{x^2}.$
The multiplication by $e^{x^2}$ and 
$e^{\cD^2/4}$ are well defined on the space
$V_{2n+1}$ because $x$ and $\cD$ are
nilpotent. Similarly, $\bF_-(f(x))=\bF_+(f(-x)).$
Another (equivalent) approach requires the truncated
nonsymmetric Bessel functions, which will not be discussed
here. See \cite{CM}. 

\begin{theorem}\label{t55} The truncated inversion reads as 
$$
\bF_-\circ\bF_+\ =\ (-1)^n\,\hbox{id}\ =\ \bF_+\circ\bF_-.
$$
Let $\widehat f=\bF_+(f)$ and $\widehat g=\bF_+(g)$.
We have the Plancherel formula
$$\langle f,g\rangle_+=
(-1)^n\langle \widehat f,\widehat g\rangle_-.$$
\end{theorem}

{\bf Proof.} Using (\ref{fscalar}), we conclude that $F_+,$ 
inducing 
the automorphism $\omega$ on $\HH'',$ sends the anti-involution
corresponding to $\langle \cdot ,\cdot \rangle_+$ to that
of $\langle \cdot ,\cdot \rangle_-.$ We need to examine the
generators $s, x,\cD.$ The representation 
$V_{2n+1}$ is irreducible, 
so the scalar products $\langle f,g\rangle_+$
and $\langle \hat f,\hat g\rangle_-$ have to be proportional.

To find the proportionality coefficient, let us
calculate the Hankel transform explicitly. First, $\bF_+(1)$
has to be proportional to $x^{2n}$ since $1$ is a unique
eigenvector of $\cD$ and $\omega$ sends $\cD$ to $-2x.$
Hence $\bF_+(x^m)$ is proportional to $\cD^m\bF_+(1),$
that is $x^{2n-m}$. Second,
\begin{align}
&\bF_+(e^{-x^2})=(e^{x^2}e^{\cD^2/4}e^{x^2})e^{-x^2}=e^{x^2},
\hbox{\ therefore\ }
\notag\\
&\bF_+(x^{2m})=(-1)^m\frac{m!}{(n-m)!}x^{2n-2m},
\hbox{\ and\ }\notag\\
&\bF_+(x^{2m+1})=(\cD/2)\bF_+(x^{2m})=(-1)^m\frac{m!}
{(n-m-1)!}x^{2n-2m-1}\notag.
\end{align}  
In particular $\bF_+(1)=\frac{x^{2n}}{n!}$ and
$\bF_+(x^{2n})=(-1)^nn!$. We get the inversion formula.
The Plancherel formula holds since
$$\langle 1,x^{2n}\rangle_+=1 \hbox{\ and\ }
\langle \bF_+(1),\bF_+(x^{2n})\rangle_-=(-1)^n\langle x^{2n},
1\rangle_-=(-1)^n.\hbox{\ \sq}
$$ 

\bigskip
\setcounter{equation}{0}
\section{DAHA and Macdonald's polynomials}

Recall Theorem \ref{t21}:
\begin{equation}\label{f61}
\frac{1}{i}\int_{i\BR}q^{-x^2}\delta_k(x)dx=2\sqrt{\pi a}
\prod_{j=0}^\infty \frac{1-q^{k+j}}{1-q^{2k+j}}, \; \Rea k>0,
\end{equation}
where
$$\delta_k(x)=\prod_{j=0}^\infty \frac{(1-q^{j+2x})
(1-q^{j-2x})}{(1-q^{k+j+2x})(1-q^{k+j-2x})}$$
for $q=e^{-1/a}$.

For a Laurent series $f(x)=
\sum_{n=-\infty}^\infty c_nq^{nx}$ in the variable $q^x,$ 
$\lr f(x)\rr \equal c_0$
is called the {\em constant term}. 
We expand $\delta_k(x)$ in a
Taylor series in terms of $q^k.$ It becomes a Laurent series
in terms of $q^x$ with the coefficients 
from $\BC [q^k][[q]].$ Let 
$$\widehat{\gamma_{-}}
\equal \sum_{n=-\infty}^\infty q^{nx+n^2/4}.
$$ 
It is a definition, but the right-hand side does coincide
with $q^{-x^2}$ in the space of distributions on periodic
functions in the variable $x$ with the period $\omega =2\pi ia.$ 
Formula (\ref{f61}) is equivalent to
\begin{equation}\label{f62}
\lr \widehat{\gamma_{-}}\delta_k(x)\rr=2\prod_{j=0}^\infty
\frac{1-q^{k+j}}{1-q^{2k+j}}.
\end{equation}

\subsection{Rogers' polynomials}
We are going to interpret the latter formula as a calculation
of the structural constant of the difference spherical 
Fourier transform in the $A_1$-case.
Generally speaking, the spherical transform is an 
integration with the spherical (or hypergeometric) 
functions. In this setup, 
the spherical functions are Rogers' polynomials
and the integration is the constant term functional.

\begin{definition}
Rogers' polynomials $p_n(x)\in \BC(q,q^k)[q^x+q^{-x}]$
are uniquely defined by the properties:

(a) $p_n(x)=q^{nx}+q^{-nx}+\sum_{|m|<n}c_mq^{mx},\, n>0,\ 
p_0= 1;$

(b) $\lr p_n(x)q^{mx}\delta_k(x)\rr=0$ for all $m$ with $|m|<n$.
\end{definition}

The first two nontrivial  Rogers' polynomials are:
$$p_1=q^x+q^{-x},\ p_2=q^{2x}+q^{-2x}+
\frac{(1-q^k)(1+q)}{1-q^{k+1}}.$$
The formula for $p_1$ is immediate from the following
simple lemma.

\begin{lemma}\label{l61}
Let $p_n(x)=\sum_mc_mq^{mx}$. Then $c_m\ne 0$
only for even $n-m.$ 
\end{lemma}

{\bf Proof.} The Laurent series
$\delta_k(x)$ involves only even powers of $q^x$. \sq\medskip

To calculate $p_2,$ we use that $\delta_k(x)$ 
can be replaced in (b) by 
$$\delta_k^0(x)
\equal \frac{\delta_k(x)}{\lr \delta_k(x)\rr}.
$$ 
We calculated it in Theorem \ref{t23}:
$$\delta_k^0(x)=1+\frac{(1+q)(q^k-1)}{2(1-q^{k+1})}
(q^{2x}+q^{-2x})+\ldots.$$
The constant term of $p_2$ is the first coefficient
of this expansion multiplied by $-2.$

More generally, this argument gives
that the coefficients of all
$p_n(x)$ lie in
$\BQ (q,q^k)$ since so do the coefficients of $\delta_k^0(x).$

Rogers' polynomials $p_n(x)$ play the role of
the Bessel functions in the following theorem

\begin{theorem}\label{mastersym}
\begin{align}\label{f63}
\lr p_n(x)p_m(x)&\widehat{\gamma_{-}}\delta_k^0(x)\rr=\notag\\
&= p_n(m+k/2)
p_m(k/2)q^{(m^2+n^2+2k(m+n))/4}\,
\lr \widehat{\gamma_{-}}\delta_k^0(x)\rr.
\end{align}
\end{theorem}

{\bf Proof} will be given later.\sq\medskip

Obviously (\ref{f63})  
results in the $m\leftrightarrow n\,$-invariance of the product 
$p_n(m+k/2)p_m(k/2).$ 
Vice-versa, it can be readily deduced from this
symmetry, as we will see soon. 

The most natural proof of this theorem goes via
the $q$-counterparts of the nonsymmetric Bessel functions.

\subsection{Nonsymmetric polynomials} Recall formulas
(\ref{formmu}),(\ref{fmudelta}):
$$\delta_k^0(x)=\frac{\mu_k^0(x)+
\mu_k^0(-x)}{2},$$ where
$$\mu_k(x)=\prod_{j=0}^\infty \frac{(1-q^{j+2x})
(1-q^{j+1-2x})}{(1-q^{k+j+2x})(1-q^{k+j+1-2x})},$$
$$\mu_k^0(x)=\mu_k(x)/\lr \mu_k(x)\rr=1+\frac{q^k-1}
{1-q^{k+1}}(q^{2x}+q^{1-2x})+\cdots\ .
$$

We use the following linear order on the set of 
monomials $q^{nx}\,:$
$q^{nx}\succ q^{mx}$ if
either $|m|<|n|$ or $n=-m$ and $n$ is negative
(so $q^{-x}\succ q^x$).

\begin{definition}
Nonsymmetric polynomials $e_n(x)\in
\BC(q,q^k)[q^x, q^{-x}]$ for $n\in \BZ$ are uniquely
defined by the properties:

(a) $e_n(x)=q^{nx}+\mbox{lower terms with respect
to} \succ$,

(b) $\lr e_n(x)q^{-mx}\mu_k^0(x)\rr =0$ if $q^{mx}\prec
q^{nx}$ (note the minus sign!)
\end{definition}

The definition and the beginning of
the expansion of $\mu_k^0$ give that
$$e_0=1,\,e_1=q^x,\,e_{-1}=q^{-x}+
\frac{1-q^k}{1-q^{k+1}}q^x,\, e_2=q^{2x}+
q\frac{1-q^k}{1-q^{k+1}}.
$$
Similar to $p_n,$ each $e_n(x)$ involves only monomials
$q^{mx}$ with even $n-m$  because 
$\mu_k^0(x)$ contains only even powers of $q^x.$ 
Similar to $p_n(x),$ 
$$e_n(x)\in \BQ(q,q^k)[q^x, q^{-x}],\hbox{\ since\ }
\mu_k^0(x)\in \BQ(q,q^k)[[q^x,q^{-x}]].
$$ 

The main technical advantage of the theory of $e$-polynomials
versus the $p$-polynomials is that we can construct
them using the intertwining operators of 
the double Hecke algebra. The following 
{\em creation operator} is due to Knop and Sahi
in the case of $A_n.$ See  \cite{KnS} and 
\cite{C1,C3} for the general theory.

\begin{proposition}\label{p61}
Denoting $\pi(f(x))=f(1/2-x),$
$$q^{(1-n)/2}e_n=q^x\pi(e_{1-n}(x))\hbox{\ for\ } n>0.
$$
\end{proposition}

{\bf Proof.} It is clear that 
$q^{(n-1)/2}q^x\pi(e_{1-n}(x))=
q^{nx}+\mbox{lower terms}$. 
So it suffices to check that
\begin{equation}\label{f64}
\lr q^x\pi(e_{1-n}(x))q^{-mx}\mu_k^0(x)\rr =0, \; m=n-1, n-2,
\ldots, 1-n.
\end{equation}
The following two properties of $\pi$ are evident: 
$$(a)\ \lr \pi(f)\rr=\lr f\rr,\ \
 (b)\ \pi(\mu_k^0(x))=\mu_k^0(x).
$$ 
Using them,
we see that (\ref{f64}) is true for $m=n-1, \ldots,
2-n$. The only remaining case $m=1-n$ is trivial: 
$\pi(e_{1-n}(x))q^{nx}\mu_k^0(x)$ involves only odd
powers of $q^x$. \sq\medskip

We will deduce (\ref{f63}) from its 
nonsymmetric analogue.
Setting $\eps_n(x)=e_n(x)/e_n(-k/2),$ it reads as 
\begin{equation}\label{f65}
\lr \eps_n(x)\eps_m(x)\widehat{\gamma_{-}}
\mu_k^0(x)\rr=\lr \widehat{\gamma_{-}}\mu_k^0(x)\rr\,
\eps_n(m_\#)\,q^{(m^2+n^2+2k(|m|+|n|))/4},
\end{equation}
where $n_\#\equal {n+\sgn(n)k}/2$ and we set
$\sgn(0)=-1.$ 

It is not surprising that the nonsymmetric
(more general) formula is easier to prove than
the symmetric one.  Recall that
the advantage of the nonsymmetric setting was quite
clear in the theory of the Hankel transform. 

The coefficient of proportionality in formula (\ref{f65})
is a combination of formula (\ref{f62})
and the constant term conjecture (\ref{f26}):
$$\lr \widehat{\gamma_{-}}\mu_k^0(x)\rr=
\prod_{j=0}^\infty \frac{1-q^j}{1-q^{j+k}}.$$

Formula (\ref{f65}) implies that $\eps_n(m_\#)$ is 
invariant under the change $m\leftrightarrow n$. 
This of course can be established directly. Let us do it. 

\subsection{Double affine Hecke algebra} We denote it by
"double H" $\HH.$ 
It depends on the two parameters $q^{1/2}$ and $t^{1/2}.$
The generators are $X^{\pm 1},Y^{\pm 1}, T,$ the relations
$$\left\{ \begin{array}{c}TXT=X^{-1},\;
TY^{-1}T=Y,\\ Y^{-1}X^{-1}YXT^2q^{1/2}=1,\\
(T-t^{1/2})(T+t^{-1/2})=0\end{array}\right\} .$$

If $t=1,$ then $\HH$ is the extension 
of the Weyl algebra by the reflection $S,$ i.e., 
$$\lr X, Y,S\rr/\{Y^{-1}X^{-1}YX=q^{-1/2}, S^2=1, 
SXS=X^{-1}, SYS=Y^{-1}\}.$$
\medskip

\begin{theorem}\label{t61}
(a) (PBW property) The elements $X^nT^\eps
Y^m,\; n,m \in \BZ,\; \eps=0,1,$ form a basis of $\HH$.

(b) Using $s(f(x))=f(-x)$ and
$\pi(f(x))=f(1/2-x),$  the formulas
\begin{equation}\label{f66}
T\mapsto t^{1/2}s+\frac{t^{1/2}-t^{-1/2}}{q^{2x}-1}(s-1),\;
X\mapsto q^x,\; Y\mapsto \pi T
\end{equation}
define a representation of $\HH$ in the space $\BC[q^x, q^{-x}].$
It is faithful for $q$ apart from roots of unity.
\end{theorem}

{\bf Proof.} We will start with the following
lemma, which is simple to check.

\begin{lemma}\label{l62} Setting
$\pi\equal YT^{-1},$ the algebra $\HH\,$
can be alternatively described as follows:
$$\HH=\lr T, X^{\pm 1}, \pi \rr /
\left\{ \begin{array}{c}TXT=X^{-1},
\; \pi^2=1,\; \pi X\pi^{-1}=q^{1/2}X^{-1},\\
(T-t^{1/2})(T+t^{-1/2})=0\end{array}\right\} .$$
\end{lemma}

The lemma gives that the formulas
from  part (b) of Theorem \ref{t61} define
a representation of $\HH\,.$ Indeed, the
formula for $T$ is well-known in the theory 
of the affine Hecke algebra of type $A_1.$
This operator does satisfy $TXT=X^{-1}.$
So we need to check only the relations involving
$Y$ or, equivalently, $\pi.$ The relations with $\pi$ 
are simple.  

(a) Any element of $\HH\,$ is a 
linear combination of  $X^nT^\eps Y^m.$ 
These monomials are linearly independent in $\HH\,$
if they are independent for at least one pair of special
values of $q,t.$ Let us take $q$ which is not a root
of unity. Then the images ot these monomials in 
$\BC[q^x, q^{-x}]$ are linearly independent. 
It is immediate for $t^{1/2}=1,$ i.e., for generic $t.$
It is also true for arbitrary $t,$ and not difficult to
check.

We get (a), and the remaining part of (b). \sq\medskip

{\bf Comment.}
In fact, the theorem can be reformulated as follows:
$$\BC[q^x, q^{-x}]=\Ind_{<T,Y>}^{\HH} (\BC),$$ 
where the subalgebra $\lr T, Y\rr\subset \HH\,$ acts on 
the one-dimensional
space $\BC$ via the character $T\mapsto t^{1/2},$
$Y\mapsto t^{1/2}.$ Indeed, this relation
includes the PBW-theorem
and readily results in the formulas from part (b).
Such reformulation does not include only one statement:
the fact that the polynomial representation is faithful
for generic $q.$
\sq\medskip

We introduce the {\em conjugation}
$f\mapsto \bar f$ on the polynomial representation:
$$\overline{q^x}=q^{-x},\,
\overline{q^{1/2}}=q^{-1/2},\,\overline{t^{1/2}}=t^{-1/2},
$$
and define the scalar product: 
$$\lr f, g\rr\equal 
\lr f\bar g\mu_k^0(x)\rr,\ \hbox{\ where\ } q^k=t.
$$

Here it is necessary to change the field of coefficients
$\BC$ to the ring $\BC[q^{\pm 1/2}, t^{\pm 1/2}]$
or its field of rationals. We will do the same with $\HH\,.$
From now on 
{\em $q^{1/2}$ and $t^{1/2}$ will be
considered as formal parameters,}
unless stated otherwise.
Note that the operators $T$ and $\pi$ from
(\ref{f66}) preserve the space 
$\BC[q^{\pm 1/2},t^{\pm 1/2},q^{\pm x}].$ 

\begin{theorem}\label{t62}
The operators $T, X, Y, \pi, q, t$ are
unitary with respect to the scalar product $\lr \; , \rr.$
\end{theorem}

{\bf Proof.} It is evident for
$X, q, t$. It is also clear that $\pi^*=\pi=\pi^{-1}.$
For $T,$ we use
the formula $\mu_k^0(-x)=\frac{1-q^{2x+k}}{q^k-q^{2x}}
\mu_k^0(x).$ It gives that the adjoint of $s$ is 
$(1-q^{2x}t)(t-q^{2x})s.$ Concerning $Y,$ use $Y=\pi T.$ 
\sq\medskip

Since the coefficints of polynomials $e_n(x)$
have denominators, $\HH\,$ and
the polynomial representation will be considered over
the field $\BC(q^{1/2},t^{1/2})$ in what follows. 
We will continue using the notation $\BC[q^{\pm x}].$

In fact it is sufficient to take $\BQ$ instead of $\BC.$ 
Moreover, the coefficients of $e_n$ are expressed in 
terms of $q,t.$ However the action of $T$ and $\pi$ requires
the square roots of $q,t.$ 

\begin{theorem}\label{t63}
Considering the polynomials $e_n(x)$ as
elements of the $\HH\, $-module
$\BC[q^{\pm x}],$ 
$$Ye_n(x)=q^{-n_\#}\,e_n(x),\ n_\#\equal(n+\sgn(n)k)/2,\;
\sgn(0)=-1,\, n\in \BZ.$$
\end{theorem}

{\bf Proof.} The operator $Y$ is unitary and 
preserves the filtration induced by $\prec,$
which is easy to check. 
Note that $T$ does not preserve this filtration. It readily
gives that $e_n(x)$ are eigenvectors of $Y$. Use the
definition.
Examining the leading term $q^{nx}$ of $e_n(x),$
one gets the corresponding eigenvalue.
\sq\medskip

\begin{theorem}\label{t64}
Setting $\eps_n(x)\equal e_n(x)/e_n(-k/2),$
$$\eps_n(m_\#)=\eps_m(n_\#).
$$ 
\noindent Equivalently,
the pairing $\{f\, ,\, g\}\equal f(Y^{-1})(g)(-k/2)$
is symmetric on $\BC[q^{\pm x}],$ 
where $f(Y^{-1})=f(q^x\mapsto Y^{-1}).$
\end{theorem}

{\bf Proof.} First, we introduce the anti-involution
$\phi$ of the algebra $\HH\,$, 
by setting
$$\phi(X)=Y^{-1},\ \phi(Y)=
X^{-1},\ \phi(T)=T,\ \phi(q)=q,\ \phi(t)=t.
$$
Later it will be associated with 
the Fourier transform on the space of the 
generated functions.

Second, we need the {\em evaluation map} 
$:\HH \to
\BC(q^{1/2},t^{1/2})$ defined by 
$$\{ X^nT^\eps Y^m\}
=t^{-n/2}t^{\eps/2}t^{m/2}.$$ 
It is obviously $\phi$-invariant.
Explicitly,
$$
\{ H\} \hbox{\ is\ the\ evaluation\ of\ }
H(1)\in \BC[q^{\pm x}]
$$ 
for $H\in \HH.$ 
The evaluation of a polynomial
in $q^{\pm x}$ is its value at the point $x=-k/2.$

Third, the scalar product  $\{ A, B\}=\{
\phi(A)B\}$ on $\HH$ is symmetric, i.e.,
$\{ A, B\}=
\{ B, A\}.$ It follows from 
$\{ \phi(H)\}=\{ H\}$ for  $H\in \HH.$
It coincides with $\{ f\, ,\, g\}$ on Laurent polynomials
upon the substitution $X=q^x.$

Finally, 
\begin{align}
\{ e_n(X), e_m(X)\}=\{ e_n(Y^{-1})e_m(X)\}
=e_n(m_\#)e_m(-k/2),
\end{align}
where we used that $Y^{-1}e_n(x)=q^{n_\#}e_n(x)$ in
the polynomial representation.  \sq\medskip

\subsection {Back to Rogers' polynomials}
The relation of the Macdonald polynomials
to the double affine Hecke algebras
was the first obvious confirmation of their
importance.
These algebras were designed
for a somewhat different purpose: to connect the 
Knizhnik-Zamolodchikov equation with the 
quantum many-body (eigenvalue) problem in the $q$-case.
Paper \cite{C11} is devoted to it.
It was an important step in the
new theory of generalized $q$-hypergeometric functions.
However the application to the celebrated
Macdonald constant term conjecture and 
the norm-conjecture was the first recognized success of 
the double Hecke algebras.

In our case, symmetric means even.
We will constantly use $X=q^x$ instead of $x.$
Then the reflection $s$ sends $X\mapsto X^{-1}.$ 
Concerning operators, we call them 
symmetric if they commute with $s.$
Symmetric operators on the space
$\BC[X^{\pm 1}]$ preserve the subspace $\BC[X+X^{-1}]$
of symmetric polynomials. For the operators we 
consider it is necessary and sufficient. 

Recall that 
$p_0=1$ and the other $p_n$ are defined 
from the relations
$$
p_n(x)=X+X^{-1}+\mbox{lower terms},\ \lr p_n(x)
p_m(x)\delta_k^0(x)\rr =\delta_{mn}C_n,
$$ 
where $\delta_k^0(x)
=(\mu_k^0(x)+\mu_k^0(-x))/2$ and $C_n$ are certain constants.

\begin{theorem}\label{c72}
(a) The operator $Y+Y^{-1}$ commutes
with $T$ and is symmetric. Setting  
$\ \varpi(f)(x)\equal f(x+1/2),\ $ 
$\cL\equal Y+Y^{-1}\,\mid\,_{\BC[X+X^{-1}]},$
$$\cL=\frac{t^{1/2}X-t^{-1/2}
X^{-1}}{X-X^{-1}}\varpi+\frac{t^{1/2}X^{-1}-t^{-1/2}X}
{X-X^{-1}}\varpi^{-1}.$$

(b) For any $n>0,$ 
\begin{align}
&p_n=(1+t^{1/2}T)e_n=(1+s)\left(\frac{t-X^2}
{1-X^2}e_n\right),\label{fpen} \\
&\cL p_n=(q^{n/2}t^{1/2}+q^{-n/2}t^{-1/2})p_n.\label{flpn}
\end{align}
\end{theorem}

{\bf Proof.} (a) The commutativity of $T$ and $Y+Y^{-1}$
immediately results from the defining relations. Generally,
it is due to Bernstein and Zelevinsky.
The symmetric polynomials $f\in \BC[X+X^{-1}]$ are exactly 
those satisfying $Tf=t^{1/2}f$. Therefore $Y+Y^{-1}$ is
symmetric. The calculation of the restriction of
$Y+Y^{-1}$ to $\BC[X+X^{-1}]$ is simple, especially if
one uses that it is symmetric.

(b) Let us use (\ref{fpen})
to introduce the
polynomials $\tilde p_n$. We must check that they
coincide with $p_n.$ First, 
$$
\cL \tilde p_n=(q^{n/2}t^{1/2}+q^{-n/2}t^{-1/2})\tilde p_n,
$$
since $T$ commutes with $Y+Y^{-1}$. Second, 
$(Y+Y^{-1})^*=Y+Y^{-1},$ where $\ast$ is 
with respect to the scalar product $\lr f\bar
g\mu_k^0(x)\rr$ because the operator $Y$ is unitary.
It readily results in $\cL^*=\cL,$
where $\cL^*$ is now defined with respect to the scalar product
$\lr f\bar g\delta_k^0(x)\rr,$
which is the symmetrization of the above pairing with $\mu.$ 

Note that here we can drop the conjugation
and use $\lr fg\delta_k^0(x)\rr,$ because $\cL$ is 
bar-invariant. Respectively, $\cL$ has $\ast$-invariant 
eigenfunctions.

For generic $q$ and $t,$ the
eigenvalues of $\cL$ in the space of even polynomials
are pairwise distinct and therefore the
eigenvectors are pairwise orthogonal. We get the desired
coincidence with $p.$ \sq\medskip

{\bf Comment.} When $k=1,$ we have $\delta_1(x)=
(1-q^{2x})(1-q^{-2x})$ and it is easy to show that
$p_n^{(k=1)}(x)=\frac{q^{(n+1)x}-q^{-(n+1)x}}
{q^x-q^{-x}}$. Thus $p_n^{(k=1)}(x)$ is the
character of the irreducible representation of $sl_2$
of dimension $n+1.$ Is there any reasonable
interpretation of the coefficients of polynomilas 
$e_n^{(k=1)},$ and the coefficients of general $p_n,\, e_n$
considered as rational functions or series in terms of $q,t\ $? 
\sq\medskip

\subsection{Conjugated polynomials}
Let us examine the action of the bar-involution
$\bar{X}= X^{-1},\, q^{1/2}\mapsto q^{-1/2},\,
t^{1/2}\mapsto t^{-1/2}$
on the $\eps$-polynomials.
It preserves $p_n$ and $p_n/p_n(-k/2)$ thanks to
formula (\ref{flpn}) and the relation $\cL^*=\cL.$
One can also see it using that $\delta_k^0$ is bar-invariant.

We will need the following extension to $\HH\,$
of the Kazhdan-Lusztig involution (cf. \cite{KL1}):
\begin{equation}\label{feta}
\eta(T)=T^{-1},\, \eta(\pi)=\pi,\, \eta(X)=X^{-1},\
q^{\frac{1}{2}}\mapsto q^{-\frac{1}{2}},\, 
t^{\frac{1}{2}}\mapsto t^{-\frac{1}{2}}.
\end{equation}
Introducing $T_0=\pi T\pi,$ we get $\eta(T_0)=T_0^{-1}$
so it coincides with the Kahdan-Lusztig involution
on the affine Hecke algebra genearated by $T_1=T,\ T_0.$

The importance of $\eta$ in
the theory of double affine Hecke algebras
is due to the following proposition.

\begin{proposition} \label{petabar}
For arbitrary
$f\in \BC(q^{1/2},t^{1/2})[q^{\pm x}],$
$$ 
H(\bar f)=\overline{\eta(H)f},\ H\in \HH\,. 
$$
For $m\in \BZ,$
the polynomial $\overline{e_m}(x)$
is an eigenvector of $\eta (Y)$ with the eigenvalue
$q^{m_\#}.$
\end{proposition}
{\bf Proof.} It suffices to check that $\overline{T}$
is $T^{-1},$ which is straightforward.
\sq\medskip

\begin{proposition}\label{p91}
For all $m\in \BZ,$

(a)\ $ \overline{\eps_m}(x)=t^{-1/2}T(\eps_m(x)),$

(b)\ $ \overline{\eps_m}(x)=t^{-1/2}X^{-1}\eps_{1-m}(x),$

(c)\ $ XT(\eps_m(x))=\eps_{1-m}(x).$
\end{proposition}

{\bf Proof.} 
(c) Using Proposition \ref{p61}, 
\begin{align}
&XT(e_m)=X\pi Y(e_m)=q^{-m_\#}X\pi(e_m)\notag\\
&=q^{-m_\#+m/2}e_{1-m}=q^{-k/2}e_{1-m}=t^{-1/2}e_{1-m}\notag
\end{align}
for $m>0.$ Thus (c) results from  
$e_{1-m}(-k/2)/e_m(-k/2)=t^{1/2},$ which is a simple
corollary of the evaluation formula in the next section.

To check it,
one can also proceed as follows. 
In the relation $X\pi(\eps_m)=C\eps_{1-m},$ 
we need to find $C.$
The right-hand side is $C$ at $-k/2.$ So we
need to evaluate the left-hand side at $-k/2.$ It suffices
to know $\eps_m(1/2+k/2).$ Using the duality:
$$\eps_m(1/2+k/2)=\eps_1(m_\#)=q^{m_\#}q^{k/2}.$$

(b) Thanks to the previous proposition,
$\overline{\eps_m}(x)$
is an eigenvector of $\eta (Y)$ with eigenvalue
$q^{m_\#}$. Since $\eta (Y)=YT^{-2}=\pi
Y^{-1}\pi,$ we get that $\overline{\eps_m}(x)$ is
proportional to $\pi \eps_m(x)$ and to
$X^{-1}\eps_{1-m}(x).$ Thus
$$\overline{\eps_m}(x)=CX^{-1}\eps_{1-m}(x) \hbox{\ for\ 
a\ constant\ } C,
$$ 
which has to be $t^{-1/2}$ due to $\eps_n(-k/2)=1.$

(a) It is a combination of (b) and (c).
\sq\medskip

\bigskip
\setcounter{equation}{0}
\section{Four corollaries}
Let us emphasize the main points of the previous section.
So far we have five key definitions, and five key properties
of the double Hecke algebra and the $e$-polynomials.

\subsection{Basic definitions}

(i) Double affine Hecke algebra:
$$\HH=\lr X,Y,T\rr /~\left\{ \begin{array}{cc}
(T-t^{1/2})(T+t^{-1/2})=0,&
TXT=X^{-1},\\
Y^{-1}X^{-1}YXT^2q^{1/2}=1,&TY^{-1}T=Y\end{array}\right\}.
$$

(ii) There is also an alternative description in terms of the
generators $T, X, \pi =YT^{-1}$:
$$\HH=\lr X, T, \pi \rr /~\left\{ \begin{array}{cc}
TXT=X^{-1},& \pi X\pi^{-1}=q^{1/2}X^{-1},\\ \pi^2=1,
&(T-t^{1/2})(T+t^{-1/2})=0\end{array}\right\}.
$$

(iii) The $\mu$-function $\mu_k^0(x)$  
is uniquely defined by
the properties 
$$\mu_k^0(x+1/2)=\frac{q^{2x+k}-1}{q^{2x}-q^k}\mu^0_k(x),
\hbox{\ and\ }\lr \mu^0_k(x)\rr =1,
$$
where $\lr \, \rr$ is the constant term. 

(iv) The formal conjugation is the automorphism of
$\BC[X^{\pm 1}]$ defined by $\bar X=X^{-1}, \bar q=q^{-1},
\bar t=t^{-1}$. We will constantly identify 
$X$ with $q^x.$ 

(v) The linear order $\prec$ on the set of monomials $X^n$: 
$$
X^n\prec X^m \hbox{\ if\ either\ }|n|<|m| \hbox{\ or\ }
n=-m>0.
$$  

Let us summarize what has been proved.

(a) (PBW property) The elements
$X^nT^\eps Y^m$, $n,m\in \BZ$, $\eps=0,1,$ form a 
basis of $\HH$.

(b) (Polynomial representation) For 
$s(X)=X^{-1},\, \varpi(X)=q^{1/2}X,$ 
the formulas
$$T\mapsto t^{1/2}s+\frac{t^{1/2}-t^{-1/2}}{q^{2x}-1}(s-1),\;
\pi \mapsto s\varpi,\; X\mapsto X,\; Y\mapsto s\varpi T$$
define a representation of $\HH$ in the space 
$\cP\equal \BC[X^{\pm 1}].$

(c) (Unitary structure) The operators $X, q, Y, \pi, T$ 
are unitary with respect to the scalar product.
$\lr f,g\rr = \lr f\bar g\mu_k^0(x)\rr$ on the space $\cP.$

(d) (Nonsymmetric polynomials) They can be defined 
as follows: $e_n(x)=X^n\hbox{\, mod\,} \{X^m\prec X^n\},$
$$
Ye_n=q^{-n_\#}e_n,\ n_\#=(n+\sgn(n)k)/2,\, 0_\#=-k/2.
$$

(e) (Duality) 
$
\eps_n(m_\#)=\eps_m(n_\#),\hbox{\ where\ }
\eps_n(x)\equal e_n(x)/e_n(-k/2).
$

Let us discuss applications.
Our main instruments will be the duality and
the anti-involution
$\phi:\HH \to \HH$ defined on the generators by $\phi(X)=Y^{-1}$,
$\phi(Y)=X^{-1}$, $\phi(T)=T$, $\phi(q)=q$, $\phi(t)=t$.

\subsection{Creation operators} 
In this section we will consider $q$ and $t=q^k$ as
numbers (not as formal variables). We will assume that $q$
is not a root of unity. The coefficients of the polynomials
$e_n$ are rational functions of $q$ and $q^k$, so the
polynomials $e_n(x)$ are not well defined for some particular
values of $q$ and $k$. The theory of $e$-polynomials is
(relatively) simple because they can be produced using the 
{\em intertwining operators}. 
The latter play the role of the
creation (raising) operators in Lie theory.

\begin{corollary}\label{c71}
The polynomials $e_n(x),\, e_{1-n}(x)$ for $n>1$  are well
defined for $k\not \in -\{ [n/2], \ldots, n-1\}.$
Explicitly, they can be obtained using the 
following operations:

(A) Introducing  $\Pi\equal X\pi,$ 
$$e_{1-n}=q^{-n/2}\Pi e_n\hbox{\ for\ } n\in \BZ.
$$

(B) Assuming that $q^{2n_\#}\ne 1,$
$$e_{-n}=q^{k/2}(T+\frac{t^{1/2}-t^{-1/2}}{q^{2n_\#}-1})e_n.
$$
\end{corollary}

{\bf Proof.} Recall that $Y$ preserves the filtration on
$\cP=\BC[X^{\pm 1}]$ induced by $\prec,$ and also subspaces
$\BC[X^{\pm 2}],\,$  $X\BC[X^{\pm 2}]$. The eigenvector of
$Y$ with an eigenvalue $q^{-n_\#}$ is clearly well defined
if $q^{-n_\#}\ne q^{-m_\#}$ for all $m$ such that $X^m\prec
X^n$ and $n-m$ is even. Since $q$ is not a root of unity, this
condition is equivalent to $n_\#\ne m_\#$. It always holds
when $n,m>0$ or $n,m\le 0.$ If
$n>0,\, m\le 0,$ then the conditions
$(n+k)/2\ne (m-k)/2$ or $k\ne (n-m)/2$ are sufficient.

(A) Let us check that $\phi(\Pi)=\Pi$:
$$\phi(\Pi)=\phi(X\pi)=\phi(X\pi^{-1})=\phi(XTY^{-1})=
XTY^{-1}=\Pi.$$
Using  $\Pi X\Pi^{-1}=q^{1/2}X^{-1},$
we get that $\Pi Y^{-1}
\Pi^{-1}=q^{1/2}Y.$ So we only need to calculate the 
coefficient of proportionality.
Compare the proof of Proposition \ref{p61}, which is
based directly on the definition.

(B) The relation $TXT=X^{-1}$ implies
$$\left(T+\frac{t^{1/2}-t^{-1/2}}{X^2-1}\right)X=X^{-1}
\left(T+\frac{t^{1/2}-t^{-1/2}}{X^2-1}\right).$$
Applying $\phi$ we come to
$$\left(T+\frac{t^{1/2}-t^{-1/2}}{Y^{-2}-1}\right)Y=
Y^{-1}\left(T+\frac{t^{1/2}-t^{-1/2}}{Y^{-2}-1}\right),$$
which gives the desired up to a coefficient of 
proportionality. We get the latter coefficient
from the consideration
of the action of $T$ on the leading term of $e_n.$ 
\sq\medskip

{\bf Comments.}
(i) We will later see that the polynomials $e_n,e_{1-n}, n>1,$ 
are not well defined
for $k\in -\{ [n/2],\ldots, n-1\}.$ 
The polynomial $e_n$ is well defined if and only
if $e_{1-n}$ is well defined thanks to (A).

(ii) The corollary gives an inductive procedure for calculating
the nonsymmetric polynomials, which will be used a great deal
for the classification of the finite dimensional 
representations:
$$1=e_0\stackrel{A_0}{\longrightarrow}e_1\stackrel{B_1}
{\longrightarrow}e_{-1}\stackrel{A_{-1}}{\longrightarrow}
e_2\stackrel{B_2}{\longrightarrow}e_{-2}\stackrel{A_{-2}}
{\longrightarrow}\ldots.
$$

\subsection{Standard identities} Recall that 
$\eps_n(x)=e_n(x)/e_n(-k/2)$. We continue using
$X$ instead of $q^x$ in the formulas. 

\begin{corollary}\label{c73}
(a) (Pieri rules) Setting $\nu=1$ for $m\le 0$ and
$\nu=-1$ otherwise,
\begin{align}
&X^{-1}\eps_m=\frac{t^{1/2+\nu}q^{-m+1}-t^{-1/2}}
{t^\nu q^{-m+1}-1}\eps_{m-1}-
\frac{t^{1/2}-t^{-1/2}}{t^\nu q^{-m+1}-1}\eps_{1-m},
\label{fpierin}\\
&X\eps_m=\frac{t^{-1/2+\nu}q^{-m}-t^{1/2}}
{t^\nu q^{-m}-1}\eps_{m+1}-
\frac{t^{-1/2}-t^{1/2}}{t^\nu q^{-m}-1}\eps_{1-m}.
\label{fpierip}
\end{align}

(b) (evaluation formulas) We have the equalities
\begin{equation}\label{fevaluat}
e_m(-k/2)=t^{-|m|/2}\prod_{0<j<|m|'}\frac{1-q^jt^2}{1-q^jt},
\end{equation}
where $|m|'=m$ if $m>0$ and $|m|'=1-m$ if $m\le 0$.

(c) (norm formulas) For $m,n\in \BZ,$
\begin{align}
&\lr \eps_m\, ,\, \eps_n\rr=
\lr \eps_m\overline{\eps_m}\mu_k^0(x)\rr=
\prod_{0<j<|m|'}\frac{1-q^j}{t^{-1}-q^jt},\label{fnormep}\\
&\lr e_m\, ,\, e_n\rr=
\prod_{0<j<|m|'}\frac{(1-q^j)(1-q^j t^2)}
{(1-q^j t)(1-q^{j}t)}.\label{fnorme}
\end{align}
\end{corollary}

{\bf Proof.} (a) By Theorem \ref{t63}, $Y\eps_n=q^{-n_\#}\eps_n.$ 
Evaluating this equality at points
$X=q^{m_\#},$ we get $(Y\eps_n)(m_\#)=q^{-n_\#}\eps_n(m_\#).$
Let $m\le 0.$
Using the formula for the action of $Y$ in $\cP,$ 
$$t^{1/2}\eps_n((m-1)_\#)+\frac{t^{1/2}-t^{-1/2}}
{q^{1-2m_\#}-1}(\eps_n((m-1)_\#)-\eps_n((1-m)_\#))=$$
$$=q^{-n_\#}\eps_n(m_\#).
$$
Now we apply Theorem \ref{t64}: 
$$t^{1/2}\eps_{m-1}(n_\#)+\frac{t^{1/2}-t^{-1/2}}
{q^{1-m}t-1}(\eps_{m-1}(n_\#)-\eps_{1-m}(n_\#))=
q^{-n_\#}\eps_m(n_\#).$$
This gives formula (\ref{fpierin})) 
at points $x=n_\#.$
Since there are infinitely many such points we get the first
formula in (a) for $m\le 0.$ Positive $m$ are considered
in the same way. The relations
\begin{equation}\label{fpieps}
\eps_m=q^{(m-1)/2}\Pi \eps_{1-m},
\end{equation}
connect the second formula with the first.

(b) Assume that $m\le 0.$ 
The leading terms on the 
left-hand side and
the right-hand sides of (\ref{fpierin}) for $\nu=1$ are 
$(1/e_m(0_\#))X^{m-1}$ and respectively, 
$$\frac{(t^{3/2}q^{-m+1}-t^{-1/2})}
{e_{m-1}(0_\#)(tq^{-m+1}-1)}X^{m-1}.
$$
This implies (b). The case $m>0$ is analoguos.

(c) Let $m\le 0.$ Using (\ref{fpierin}):
$$\lr \eps_m, \eps_m\rr=\lr X^{-1}\eps_m,
X^{-1}\eps_m\rr=$$
$$=\frac{t^{3/2}q^{-m+1}-t^{-1/2}}{tq^{-m+1}-1}
\cdot \frac{t^{-3/2}q^{m-1}-t^{1/2}}{t^{-1}q^{m-1}-1}
\lr \eps_{m-1}, \eps_{m-1}\rr +$$
$$+\frac{t^{1/2}-t^{-1/2}}{tq^{-m+1}-1}\cdot
\frac{t^{-1/2}-t^{1/2}}{t^{-1}q^{m-1}-1}\lr
\eps_{1-m},\eps_{1-m}\rr.$$
Involving $\Pi,$ which is unitary  by Corollary \ref{c71} (b), 
\begin{align}\label{flepepr}
&\lr \eps_{1-m},\eps_{1-m}\rr=
\lr \eps_m,\eps_m\rr\hbox{\ since\ }
\eps_m=q^{(m-1)/2}\Pi \eps_{1-m},\\
&\lr \eps_m, \eps_m\rr=\frac{t^{-1}-tq^{-m+1}}
{1-q^{-m+1}}\lr \eps_{m-1}, \eps_{m-1}\rr\notag
\end{align}
and we get (c) by induction. In the case
$m>0,$ one uses (\ref{flepepr}).
\sq\medskip

{\bf Comments.} (i) The existence of the three-term relation 
in the form (a) can be seen directly
from the orthogonality of polynomials
$e_n(x)$. The approach based on the duality is simpler and
readily gives the exact coefficients. 


(ii) Formula (b) shows that the polynomial $e_n$ is not
well defined for values of $k$ from Corollary \ref{c71} (a),
because otherwise $e_n(-k/2)$ would be defined. So we have
the complete list of singular $k$ for each $e_n.$ 

(iii) Formula (b) as $k=1$ 
is a polynomial in the variable $q$
with positive integral coefficients for either positive
odd $m$ or negative even $m$. These coefficients 
have some combinatorial interpretation. 
Note that $p_n^{(k=1)}(1/2)$ is the so-called $q$-dimension
of the representation of $sl_2$ of dimension $n+1.$ 
There is a partial interpretation
of $e_n(-k/2)$ and $p_n(k/2)$ considered as
a series in terms of $q,t^{1/2}.$ However it is for
$A_1$ only and there is no known connection with the 
representation theory.

\subsection{Change $k\mapsto k+1$} 
We are going to discuss the shift--formula in
the presence of the polynomials $p_n^{(k)}(x).$ 
The $p$-polynomials appear because of the following
corollary.

\begin{corollary} \label{pshiftp}
Let $\cA=(t^{1/2}X-t^{-1/2}X^{-1})$
and $\cB=(t^{1/2}Y^{-1}-t^{-1/2}Y)$. The operator
$S\equal t^{-1/2}\cA^{-1}\cB$ preserves
the space $\BC[X+X^{-1}].$  Its restriction 
to this space is
$\frac{\varpi^{-1}-\varpi}{X-X^{-1}}$
for $\ \varpi(f)(x)\equal f(x+1/2)\ ,$ and
\begin{equation}\label{fsonp}
S(p_n^{(k)}(x))=
(q^{-n/2}-q^{n/2})p_{n-1}^{(k+1)}\hbox{\ for\ } n>0.
\end{equation}
\end{corollary}

{\bf Proof.}
It is straightforward to check the following equality
in $\HH$:
$$(T+t^{-1/2})(t^{1/2}X-t^{-1/2}X^{-1})=(t^{1/2}X^{-1}-
t^{-1/2}X)(T-t^{1/2}).$$
As a corollary, we get that the eigenspace of $T$ in 
$\cP$
with the eigenvalue $-t^{-1/2}$ is $\cA \BC[X+X^{-1}]$. Applying
the anti-involution $\phi$ to this equality we get
$$(t^{1/2}Y^{-1}-t^{-1/2}Y)(T+t^{-1/2})=
(T-t^{1/2})(t^{1/2}Y-t^{-1/2}Y^{-1}).$$
Therefore $\cB$ maps $\BC[X+X^{-1}]$ to the
$T$-eigenspace of $\cP$ with eigenvalue
$-t^{-1/2},$ that is $\cA \BC[X+X^{-1}].$ This
shows that $S$ is well defined on
$\BC[X+X^{-1}]$ and, moreover,
preserves this space. The formula for the restriction
of $S$ to $\BC[X+X^{-1}]$ is simple. Cf. Theorem \ref{c72}.

The formula for $S$ (or directly the
definition) givs that the left-hand side and
the right-hand side of \ref{fsonp} coincide.
Therefore it suffices to check
$$\lr S(p_n^{(k)}(x))g(x)\delta_{k+1}^0(x)\rr =0
$$
for any $g(x)\in \BC[q^x+q^{-x}]$ of degree $m<n-1.$
Let us prove that
$$
\lr S(p_n^{(k)}(x))g(x)\mu_{k+1}^0(x)\rr=0,
$$ 
for any $g(x)\in \BC[q^{\pm x}]$ of degree $m<n-1.$ 
Symmetrizing the latter for even $g$, we get the former. 

First,
$$\mu^0_{k+1}(x)=\frac{1-q^{k+1}}{(1-q^{2k+1})
(1+q^{k+1})}(1-q^{k+2x})(1-q^{k+1-2x})\mu_k^0(x)=$$
$$=\cA (x)\tilde \cA (x)\mu_k^0(x),$$
where $\cA(x)=t^{1/2}X-t^{-1/2}X^{-1}=q^{(k+2x)/2}-
q^{-(k+2x)/2}$ and $\tilde\cA$ is a polynomial
in $q^x$ of first degree. Second,
$$\lr S(p_n^{(k)}(x))g(x)\mu_{k+1}^0(x)\rr=
\lr \cB p_n^{(k)}(x)\tilde \cA g(x)\mu_k^0(x)\rr.$$
Here $\tilde \cA g(x)$ is a polynomial of degree
$m'<n$. Third, $p_n^{(k)}(x)$ is a linear
combination of $e_n^{(k)}(x)$ and $e_{-n}^{(k)}(x)$
by Theorem \ref{c72} (b). The same holds for
$\cB p_n^{(k)}(x).$
We get the desired result.
\sq\medskip

\subsection{Shift--formula}
Let us discuss an analytic interpretation of this corollary.
We need to choose the space of functions $\ccF$ and
the integration $\varrho.$  
The following cases can be considered:

(a) $\ccF$ is the space of even Laurent series in terms of $q^{x}$
with the constant term functional taken as $\varrho;$ 

(b) $\ccF$ is  the space of even analytic functions on the strip
$$
-1/2-\delta\le \Rea x\le 1/2+\delta 
\hbox{\ for\ } \delta>0,\ 
\hbox{\ and\ } \varrho=\int_{i\BR};
$$ 

(c) $\ccF$ is the space of continuous even functions on
the real line,  and $\varrho=\int_{\BR}$;

(d) $\ccF$ is a space of even functions which are
defined and continuous on the two lines
$z=\pm\eps \imath+\BR$ for small $\eps>0,$ i.e.,
on $C\cup \{-C\}$ for the 
the sharp integration path $C$ from (\ref{fpsish}),
and, moreover, are analytic in the rectangle
$$\{x\, \mid \, -\eps<\Ima x<\eps,\ 
-1/2-\delta\le \Rea x\le 1/2+\delta\}.$$ 
In this case, 
$$\varrho=2\int_C=\int_{\widehat{C}},\hbox{\ where\ }
\widehat{C}=C\cup\{-\overline{C}\}
$$
for the complex conjugation (preserving $C$ but changing
its orientation). 

In all cases, the integration $\varrho$ is 
invariant under the change of variable $x\mapsto -x,$
because the integration paths are $s$-invariant. 
Note that, topologically, $\widehat{C}$ 
is the cross of the diagonals at the origin both directed 
upwards. On even functions, $\frac{1}{2}\int_{\widehat{C}}$ 
coincides with the sharp integration from (\ref{fpsish}).  
It is zero on odd functions as well as for
the other integrations.

The analytic properties of $f$ are
necessary to ensure the invariance of $\varrho$
under the change of variables $x\mapsto x\pm 1/2.$
The corresponding domains of analyticity allow  deforming
the contours after the shift by $\pm 1/2.$ The real case
is the most relaxed: no domains are necessary. However 
the $q$-Gauss integrals are the most difficult to calculate
in this case.

\begin{corollary}\label{c74}
Let $\ccF, \varrho$ be one of the 
pairs above. Then for any $f\in \ccF$ such
that either the left-hand or 
the right-hand side below is well defined, 
the other is well defined too and 
$$\varrho(S(f)p_{n-1}^{(k+1)}(x)\delta_{k+1}(x))=
q^k(q^{(-k-n)/2}-
q^{(k+n)/2})\varrho(fp_n^{(k)}(x)\delta_k(x)).$$
\end{corollary}
{\bf Proof} will be  given below.
\sq\medskip

It readily results in the main property of the $q$-Mellin
transform, namely, Theorem \ref{t25}.
  
For $g_k=\prod_{j=0}^\infty \frac{1-q^{k+j}}{1-q^{2k+j}},$
we defined the $q$-Mellin transform as follows: 
$$\Psi_k(f(x))=
\frac{1}{g_k}\varrho(f(x)\delta_k(x)).$$ 
Let us deduce from the corollary that 
\begin{equation}\label{f74}
\Psi_k(f(x))=(1-q^{k+1})\Psi_{k+1}(f(x))+q^{k+3/2}
\Psi_{k+2}(S^2(f(x))).
\end{equation}

First of all,
$$\delta_{k+1}(x)=(1-tX^2)(1-tX^{-2})\delta_k(x)=$$
$$=\left(1+t^2-t\left(p_2^{(k)}(x)-\frac{(1-t)(1+q)}
{1-tq}\right)\right)\delta_k(x)=$$
$$=\left(\frac{(1+t)(1-t^2q)}{1-tq}-tp_2^{(k)}(x)\right)
\delta_k(x).
$$
Using Corollary \ref{c74}  twice,
\begin{align}
&\varrho(f(x)\delta_{k+1}(x))\notag\\
&=\frac{(1+t)(1-t^2q)}{1-tq}\varrho(f(x)\delta_k(x))-
t\varrho(f(x)p_2^{(k)}(x)\delta_k(x))=\notag\\
&=\frac{(1+t)(1-t^2q)}{1-tq}\varrho(f(x)\delta_k(x))-\notag\\
&-\frac{t^{-1}q^{-1}}{(tq-t^{-1}q^{-1})(tq^{3/2}-
t^{-1}q^{-3/2})}\varrho(f(x)\delta_{k+2}(x)),\notag
\end{align}
which is equivalent to formula (\ref{f74}).

\subsection{Proof of the shift--formula}
First, we  switch from $\delta$ to $\mu:$
$$\delta_{k+1}(x)=(1-q^{-2x})(1-q^{k+2x})\mu_k(x).$$
Let us rewrite it in terms of $\cA,$ $q^x=X,$ and
$q^{k/2}=t^{1/2}:$
$$\cA^2\mu_k(x)=-t^{-1/2}\frac{\cA}{X-X^{-1}}
\delta_{k+1}(x).
$$
Since $f(x)$ is an even function, 
\begin{align}\label{f71}
&\varrho(\cA^2f(x)\mu_k(x))=\varrho(-t^{-1/2}\frac{\cA}
{X-X^{-1}}f(x)\delta_{k+1}(x))=\notag\\
&=-\frac{1+t^{-1}}{2}\varrho(f(x)\delta_{k+1}(x)).
\end{align}
When replacing $\frac{\cA}
{X-X^{-1}}f(x)$ by its symmetrization,
we used the
invariance of $\varrho$ and $\delta_{k+1}(x)$ under the
change $x\mapsto -x.$
 
Second, we take an even function $g(x)$ 
of the same type as $f$ such that $\overline{g(x)}=g(x).$
For instance, it can be a Laurent polynomial or
convergent Laurent series
with the coefficients invariant under $q\mapsto q^{-1}$
and $t\mapsto t^{-1}.$ It is not restrictive because
later $g(x)$ will be one of
$p_n^{(k)},$ which are bar-invariant. One gets
\begin{align*}\label{f72}
&\frac{1+t^{-1}}{2}\varrho(S(f(x))S(g(x))\delta_{k+1}(x))\\
&\stackrel{(1)}{=}-\varrho(\cA^2S(f(x))S(g(x))\mu_k(x))
\stackrel{(2)}{=}\varrho(\cA^2S(f(x))\overline{ S(g(x))}\mu_k(x))\\
&\stackrel{(3)}{=}-\varrho(\cB(f(x))\overline{\cB(g(x))}\mu_k(x))
\stackrel{(4)}{=}\varrho(\cB^2(f(x))\overline{g(x)}\mu_k(x)),
\end{align*}
where the first equality follows from (\ref{f71}), the
second follows from $\overline{S}=-S$, and the
third follows from $\bar \cA=-\cA.$  Only the fourth
equality requires some comment. In the algebraic
variant, i.e., for the constant term integration (a), 
we can simply use that $Y$ is unitary. See
Theorem \ref{t62}. The same argument is applied for
the other three integrations since
$\varrho$ is invariant under  $x\mapsto -x$ and 
$x\mapsto x\pm 1/2.$  Note that the latter symmetry is
necessary when collecting $\cB$ together.

Third, we apply the $t$-symmetrizer
$P\equal (1+t)^{-1}(1+t^{1/2}T)$. It projects
$\BC[X^{\pm 1}]$ onto $\BC[X+X^{-1}].$
One has
\begin{equation}\label{fpbp}
P\cB^2P=-\cB \cB'P\hbox{\ where\ } \cB'=t^{1/2}Y-
t^{-1/2}Y^{-1}.
\end{equation}
Indeed, $\cB \cB'=-Y^2-Y^{-2}+t+t^{-1}$ 
commutes with $T$ and, in
particular, is even. So it suffices to check that
$P\cB (\cB+\cB') (f)=0$ for any $f\in \BC[X+X^{-1}]$.
However $\cB+\cB'=(t^{1/2}-t^{-1/2})(Y+Y^{-1})$ commutes
with $T.$ Hence $f'=(\cB+\cB') (f)\in
\BC[X+X^{-1}]$ and  $T\cB (f')=-t^{-1/2}f'.$
Cf. the proof of Corollary \ref{pshiftp} above. 
Formula (\ref{fpbp}) is checked.

Combining the previous formulas, we get
$$1/2(1+t^{-1})\varrho(S(f(x))S(g(x))\delta_{k+1}(x))=
\varrho(\cB^2(f(x))\overline{g(x)}
\mu_k(x))=$$ $$=\varrho(\cB^2(P(f(x)))
\overline{P(g(x))}\mu_k(x))=\varrho(P\cB^2P(f(x))
\overline{g(x)}\mu_k(x))=$$ $$=-\varrho(\cB \cB'(f(x))
\overline{g(x)}\mu_k(x))=-\varrho(
f(x)\overline{\cB\cB' (g(x))}\mu_k(x))=$$ $$=
-1/2(1+t)\varrho(f(x)\overline{\cB\cB' (g(x))}
\delta_k(x)).
$$
Thus the relation
\begin{equation}\label{f73}
\varrho(S(f(x))S(g(x))\delta_{k+1}(x))=
-t\varrho(f(x)\overline{\cB\cB'(g(x))}\delta_k(x))
\end{equation}
is proved, provided that $f(x), g(x)\in \BC[q^x+q^{-x}]$
and $\overline{g(x)}=g(x)$.

Now we simply take $g(x)=p_n^{(k)}(x)$ in (\ref{f73}) 
for  $n\ge 1.$ Using Theorem \ref{c72}, 
\begin{align*}
&\cB\cB'(p_n^{(k)}(x))=(-q^nt-q^{-n}t^{-1}+t+
t^{-1})p_n^{(k)}(x) \\
&=(tq^{n/2}-t^{-1}q^{-n/2})
(q^{-n/2}-q^{n/2})p_n^{(k)}(x).
\end{align*} 
Using Corollary \ref{pshiftp}, 
$$\varrho(S(f)p_{n-1}^{(k+1)}\delta_{k+1})=
-t(tq^{n/2}-t^{-1}q^{-n/2})\varrho(f(x)p_n^{(k)}
\delta_k).$$
Corollary \ref{c74} is proved. 

\bigskip
\setcounter{equation}{0}
\section{Difference Fourier transforms}
Recall the notation. 
We set $0_\#=-k/2$ and
$n_\#=(n+\sgn(n)k)/2$ for integers $n\neq 0,$
where $t=q^k.$
For generic $q, t,$ 
the Laurent polynomials $\eps_n(x)$ 
are uniquely defined by the properties
$Y\eps_n(x)=q^{-n_\#}\eps_n(x)$ and $\eps_n(0_\#)=1$.
We permanently use $X=q^x.$

The Gaussian is $\widehat{\gamma_-}=
\sum_{n=-\infty}^\infty q^{n^2/4}X^n.$  We will also use
$\gamma,$ a solution of
the corresponding difference equation treated as a formal
symbol. When Jackson's summation is considered, $\gamma=q^{x^2}.$

The scalar product $\lr f,g\rr=\lr f\bar{g}\mu_k^0\rr$
is given in terms of
$$\mu_k^0(x)=1+\frac{q^k-1}{1-q^{k+1}}(q^{2x}+q^{1-2x})
+\ldots,
$$
the constant term functional $\lr \cdot \rr,$ 
and the involution $f\mapsto \bar f$ of the space
$\cP=\BC(q^{1/2},t^{1/2})[X^{\pm 1}]:$  
$\ \bar X=X^{-1}$,
$\bar q^{\,1/2}=q^{-1/2}$, $\bar{t}^{\,1/2}=t^{-1/2}$. 
Our aim is the following theorem.

\begin{theorem}\label{masternon}
(Master formula) For arbitrary $m,n\in \BZ,$
\begin{align}\label{f81}
&\lr \eps_n(x)\eps_m(x)\widehat{\gamma_{-}}\mu_k^0(x)\rr
=q^{\frac{m^2+n^2+2k(|m|+|n|)}{4}}\eps_m(n_\#)\lr \widehat
{\gamma_{-}}\mu_k^0(x)\rr,\\
\label{f82}
&\lr \eps_n(x)\overline{\eps_m}(x)\widehat{\gamma_{-}}
\mu_k^0(x)\rr =q^{\frac{m^2+n^2+2k(|m|+|n|)}{4}}
\overline{\eps_m}(n_\#)\lr \widehat{\gamma_{-}}
\mu_k^0(x)\rr.
\end{align}
\end{theorem}

Before proving the theorem,
we need to establish several general facts.

\subsection{Functional representation} 
For $\BZ_\#\equal\{ n_\#|n\in
\BZ \} \subset \BC,$ let $\hat \cF$ be the space of
functions on $\BZ_\#$ and $\cF \subset \hat \cF$
the subspace of functions with compact support.
We have an evident {\em discretization map}
$\chi: \BC[X^{\pm 1}]\to \hat \cF$.

\begin{theorem}\label{t81}
The space $\hat \cF$ has a natural
structure of $\HH$-module making the map
$\chi: \cP\to \hat \cF$ an $\HH\,$-homomorphism;
$\cF\subset \hat \cF$ is an $\HH\,$-submodule.
\end{theorem}

{\bf Proof.} Let $f\in \hat \cF$. We define
$(Xf)(n_\#)=q^{n_\#}f(n_\#)$, $(\pi f)(n_\#)=
f(1/2-n_\#).$ Note that the set $\BZ_\#$ is invariant
under  $\pi: x\mapsto 1/2-x.$ It is not
$s$-invariant. However the operator  
\begin{equation}\label{ftsharp}
(Tf)(n_\#)=
\frac{t^{1/2}q^{2n_\#}-t^{-1/2}}{q^{2n_\#}-1}
f(-n_\#)-\frac{t^{1/2}-t^{-1/2}}{q^{2n_\#}-1}f(n_\#).
\end{equation}
is well defined because $-n_\# \not\in \BZ_\#$ 
for $n=0,$ precisely when first term of (\ref{ftsharp})
vanishes. \sq\medskip

\begin{definition}
(a) The Fourier transform $\bS$
acts from $\cP$ to $\cF$ and is given by the formula
$$\bS(f)(n_\#)\equal \lr f\eps_n(x)\mu_k^0(x)\rr .$$

(b) Its  antilinear counterpart, conjugating $q$ and $t,$
is introduced as follows:
$$\bE(f)(n_\#)\equal \lr \bar f\eps_n(x)\mu_k^0(x)\rr .\hbox{\ \sq}
$$
\end{definition}

The fact that functions $\bS(f)$ and $\bE(f)$ have
compact supports readily follows from orthogonality
relations for the polynomials $\eps_n(x);$ any
polynomial is their linear combination.

\begin{theorem}\label{mainth7}
(a) The formulas
$$\ep (X)=Y,\; \ep (Y)=X,\; \ep (T)=T^{-1},\;\ep(q)=q^{-1};
\ep(t)=t^{-1}
$$
can be extended to an antilinear automorphism of $\HH$. 
The formulas
$$\sigma (X)=Y^{-1},\; \sigma (T)=T,\; \sigma (Y)=
q^{-1/2}Y^{-1}XY=XT^2$$
define a linear automorphism of $\HH\,,$
fixing $q,t.$
The connection with the involution $\eta$ from
(\ref{feta}) is as follows: $\eta=\ep \sigma=\sigma^{-1}\ep.$ 

(b) For any $H\in \HH$, $f\in \cP,$ 
$$\bS(H(f))=\sigma(H)(\bS(f)),\; \; \bE(H(f))=
\ep(H)(\bE(f)).$$
\end{theorem}

{\bf Proof.} Claim (a) is straightforward.
Let us check (b) for the operator $\bE.$
It holds for $Y$ because of the definition of $\eps_n,$
and therefore for $X$ thanks to the duality 
(Theorem \ref{t62}). Hence it holds for 
$T^2=q^{-1/2}X^{-1}Y^{-1}XY$ and for 
$T=(t^{1/2}-t^{-1/2})^{-1}(T^2-1).$
Here we use that $X,T,Y$ are unitary
for $\lr\,f\,,\,g\,\rr.$ The special case $t=1$
is not a problem because it is sufficient to
check (b) for generic $t.$ 

The automorphism corresponding to the operator $\bS$
is $\ep\eta=\sigma$ thanks to Proposition \ref{petabar}.
\sq\medskip

We already defined automorphisms
$\ep, \sigma, \eta=\ep \sigma$ of the algebra
$\HH\,$ and the anti-involution $\phi.$ 
Completing the preparation to
proving the Master formula, let us introduce two more
automorphisms.

\begin{proposition}\label{p81}
(a) The formulas
$$\tau_+: X\mapsto X,\; T\mapsto T,\; Y\mapsto q^{-1/4}XY,$$
$$\tau_-: Y\mapsto Y,\; T\mapsto T,\; X\mapsto q^{1/4}YX$$
can be extended to linear automorphisms of $\HH\,$.

(b) We have the following identities in $\Aut (\HH)$:
$$\tau_+\tau_-^{-1}\tau_+=\sigma =
\tau_-^{-1}\tau_+\tau_-^{-1},\ \sigma^2=T^{-1}(\cdot)T,$$
$$\sigma \tau_+\sigma^{-1}=\tau_-^{-1},\ \ep \tau_+\ep=
\tau_-,\ \phi \tau_+\phi=\tau_-.$$
\end{proposition}

{\bf Proof.} Straightforward. \sq\medskip

{\bf Comment.} It is immediate from (b) that
$\sigma^2$ commutes with $\tau_+$ and $\tau_-.$
This also follows from $\sigma^2=T(\cdot)T^{-1}.$
The automorphism $\tau_-$ 
is inner in the representation $\cP.$  This is not
true for $\tau_+$. \sq\medskip

The automorphism $\tau_+$ is directly related to
the Gaussian $\widehat{\gamma_{-}}.$ 
Namely, 
$$H\widehat{\gamma_{-}} =\widehat{\gamma_{-}}\tau_+(H)
$$
in the polynomial representation.
Actually what we need here is an
abstract function $\gamma$ satisfying the relation:
$$ \gamma(x+1/2)=q^{1/4}q^x\gamma(x),$$
which readily results in
$$
 \gamma H \gamma^{-1}=\tau_+(H)
$$
in any functional space.

\subsection{Proof of the Master formula}
We start with $\bS.$
Treating  $\gamma^{-1}$ as $\widehat{\gamma_{-}}$
and $\gamma$ as $q^{x^2},$ 
we put 
$$\eps_m^-(x)\equal \widehat{\gamma_{-}}\eps_m(x),\ 
\eps_m^+(n_\#)\equal (\gamma \eps_m)(n_\#).
$$ 
Note that they belong to different spaces.
Then (\ref{f81}) is
equivalent to the formula
$$\bS (\eps_m^-)(n_\#)=q^{m_\#^2-0_\#^2}\eps_m(n_\#)
q^{n_\#^2}\lr \widehat{\gamma_{-}}\mu_k^0(x)\rr$$
or, equivalently, to
$$\bS (\eps_m^-)=q^{m_\#^2-0_\#^2}\eps_m^+
\lr \widehat{\gamma_{-}}\mu_k^0(x)\rr .$$

Actually it is sufficent to establish that
\begin{equation}\label{f84}
\bS(\eps_m^-)=C_m\eps_m^+
\end{equation}
for a constant $C_m,$ since the left-hand side
of (\ref{f81}) is $m\leftrightarrow n$-symmetric
and this constant can be only $1.$

Introducing the operators
\begin{equation}\label{fyplmin}
Y^+\equal \gamma Y\gamma^{-1},\ \
Y^-\equal \gamma^{-1}Y\gamma,
\end{equation}
we claim that $\sigma(Y^-)=Y^+$. Indeed,
$Y^-=\tau_+^{-1}(Y)=q^{1/4}X^{-1}Y$ and
$$
\sigma (Y^-)=\sigma(q^{1/4}X^{-1}Y)=q^{-1/4}YY^{-1}XY=
q^{-1/4}XY=\tau_+(Y).
$$ 

Obviously $Y^{\pm}(\eps_m^\pm)=q^{-m_\#}\eps_m^\pm.$
For generic  $q$ and $k,$ the eigenvalues are distinct,
the $\HH\,$-module  $\cP\widehat{\gamma_{-}}$
is irreducible, and $\{\eps_m^-\}$ is its basis.
Similarly, all eigenvectors of $Y^+$ in $\cF\gamma$
are simple. This gives (\ref{f84}) and therefore (\ref{f81}).
Now we can make $q,t$ arbitrary, provided that (\ref{f81})
is well defined. 

Switching to (\ref{f82}), its
left-hand side is not $m\leftrightarrow n$-symmetric 
anymore. So we have to proceed in a slightly different way.

First, we set
$\bar \eps_m^{(-)}\equal 
\overline{\eps_m}(x)\widehat{\gamma_{-}},\ \ $ 
$\bar \eps_m^{(+)}(n_\#)\equal 
\overline{\eps_m}(n_\#)\gamma(n_\#),\ $ 
and come to the relation
\begin{equation}\label{f85}
\bS(\bar \eps_m^{(-)})(n_\#)=C_m\overline{\eps_m}(n_\#)
q^{n_\#^2}.
\end{equation}

Indeed, for $\widetilde Y\equal \eta(Y),$ 
$$
\widetilde Y(\overline{\eps_n})=q^{n_\#}
\overline{\eps_n}
$$  
and $\bar \eps_m^{(\pm)}$ are eigenvectors of
the operators
$$
\widetilde Y^+\equal \gamma
\widetilde Y\gamma^{-1}=\tau_+(\widetilde Y),\  \widetilde Y^-
\equal \gamma^{-1}\widetilde Y\gamma=\tau_+^{-1}(\widetilde Y).
$$
Then we observe that $\sigma(\widetilde Y^-)=\widetilde Y^+,$
which results from the relation
$\sigma=\tau_+\tau_-^{-1}\tau_+:$
\begin{align}
&\sigma(\widetilde Y^-)=\sigma\tau_+^{-1}\eta(Y)=
\sigma\tau_+^{-1}\ep\sigma(Y)=
\tau_+(\tau_+^{-1}\sigma\tau_+^{-1})\ep\sigma(Y)\notag\\
&=\tau_+\tau_-^{-1}\ep\sigma(Y)=\tau_+\ep\sigma\tau_-^{-1}(Y)=
\tau_+\ep\sigma(Y)=\widetilde Y^+\notag.
\end{align}
Formula (\ref{f85}) is checked. 

Second, we establish that
\begin{equation}\label{f86}
\bE (\eps_n^-)(m_\#)=D_n\overline{\eps_n}(m_\#)
q^{m_\#^2}.
\end{equation}
This formula follows from the
identity $\ep (Y^-)=\tau_+\eta(Y)$
for the operator $Y^{-}$ from (\ref{fyplmin}).
Let us prove it:
\begin{align}
&\ep (Y^-)=\ep\tau_+(Y)=
\tau_+\eta(\sigma^{-1}\ep\tau_+^{-1}\ep\tau_+)(Y)\notag\\
&=\tau_+\eta(\sigma^{-1}\tau_-^{-1}\tau_+)(Y)
=\tau_+\eta\tau_-(Y)=\tau_+\eta(Y).\notag
\end{align}

Third, formulas (\ref{f85}) and (\ref{f86}) result in
$C_m=q^{m_\#^2}C$, $D_n=q^{n_\#^2}C,$ where $C$
does not depend on $m$ and $n,$ which concludes
the proof of formula (\ref{f82}). \sq\medskip

\subsection{Topological interpretation}
We use that the relations of $\HH\,$
are mainly of group nature and introduce the {\em group}
$$
\cB_q=\lr T, X, Y, q^{1/4}\rr / \left \lr 
\begin{array}{c}TXT=X^{-1},
TY^{-1}T=Y,\\ Y^{-1}X^{-1}YXT^2q^{1/2}=1
\end{array}\right \rr.
$$
Now $T,X,Y,q^{1/4}$ are treated as group generators,
$q^{1/4}$ is central.
The double affine Hecke algebra $\HH\,$ is the
quotient of the group algebra of $\cB_q$ by the quadratic
Hecke relation. It is easy to see that the change of
variables $q^{1/4}T\mapsto T$, $q^{-1/4}X\mapsto X$,
$q^{1/4}Y\mapsto Y$ defines an isomorphism
$\cB_q\cong \cB_1\times \BZ,$ where the generator of
$\BZ$ is $q^{1/4}.$ 

In this section, we give 
a topological interpretation
of the group
$$\cB_1=\lr T, X, Y\rr / \lr TXT=X^{-1}, TY^{-1}T=Y,
Y^{-1}X^{-1}YXT^2=1\rr.$$

Let $E$ be an elliptic curve over $\BC$,
i.e., $E=\BC/\Lambda$ where $\Lambda=\BZ+\BZ\imath.$
Topologically, the lattice can be arbitrary. 
Let $o\in E$ be the zero point, and
$-1$ the automorphism $x\mapsto -x$ of $E$.
We are going to calculate the fundamental group of
the space $(E\setminus o)/\pm 1=\BP_{\BC}^1\setminus o.$ 
Since this space is contractible, its usual
fundamental group is trivial. We can take the quotient after
removing all (four) ramification points of $-1.$ However
it would enlarge the fundamental group dramatically.
  
So we need to understand
this space in a more refined way. 
We take the base point $p=-\eps-\eps\imath\in \BC$ 
for small $\eps>0.$

\begin{proposition}\label{t82}
We have an isomorphism $\cB_1\cong
\pi_1^{orb}((E\setminus o)/\pm 1)$ where
$\pi_1^{orb}(\cdot)$ is the orbifold
fundamental group, which will be defined in
the process of proving.
\end{proposition}

{\bf Proof.} The projection map $E\setminus o\to
(E\setminus o)/\pm 1=\BP_{\BC}^1\setminus o$ 
has three branching points, which
come from the nonzero points of order 2 on $E.$ 
So by definition, 
$\pi_1^{orb}((E\setminus o)/\pm 1)$ is 
generated by three {\em involutions} $A, B, C,$ namely, 
the clockwise loops from $p$ around the branching
points in $\BP_{\BC}^1.$ 
There are no other relations. We
claim that the assignment $A=XT$, $B=T^{-1}Y$, $C=XTY$
defines a homomorphism 
$$\pi_1^{orb}((E\setminus o)/\pm 1)\to \cB_1.$$ 
Indeed, $A$ and $B$ are obviously involutive.
Concerning $C,$ the image of its square is
$$XTYXTY=T^{-1}X^{-1}YXTY=T^{-1}YT^{-1}Y=1.
$$ 
This homomorphism is an
isomorphism: $ACB=T$, $ABCA=X$
and $AC=Y$. \sq\medskip

This approach can be hardly generalized to arbitrary root
systems. The following (equivalent) constructions can.
The definition of the fundamental group is modified as follows.
See \cite{C12} and paper \cite{Io}.  

We switch from $E$ to its universal cover $\BC$ 
and define the {\em paths} 
as curves $\gamma\in \BC\setminus \Lambda$ 
from $p$ to $\widehat{w}(p),$ where
$\widehat{w}\in \widehat{W}=\{\pm 1\}\lsmash \Lambda.$ 

The {\em composition} of the paths is via $\widehat{W}:$
we add the image of the second path under $\widehat{w}$
to the first path if the latter ends at $\widehat{w}(p).$ 
The corresponding variant of Proposition \ref{t82} 
reads as follows.

\begin{proposition}\label{t82a}
The fundamental group of the above paths modulo homotopy
is isomorphic to $\cB_1,$ when 
$T$ is the half-turn, i.e., the clockwise half-circle from
$p$ to $s(p),$ $X,Y$ are $1$ and $\imath$ 
considered as vectors from  $p.$ \sq\medskip
\end{proposition}

Actually this definition is close to the calculation
of the fundamental group of $\{E\times E\setminus $\,diagonal$\}$,
divided by the transposition of the components. See 
\cite{Bi}. However there is no exact coincidence.
Let us also mention the relation to the elliptic braid group due to
v.d.Lek, although he removes all points of second order and his
group is in a sense bigger.

The topological interpretation is the best way
to understand why the group $SL_2(\BZ)$ acts on $\cB_1$
projectively.

Its elements act on $\BC$ natuarally,
by real linear transformations.
On $E,$ they commute with $-1,$
preserve $o$, and permute three other points of second order.
The position of the base point may be changed,
so we need to connect its image with the base point by
a path.
We will always do it in a small neighborhood of zero.
This makes the corresponding automorphism of $\cB_1$
unique up to powers of $T^2.$ All such automorphisms
fix $T,$ because they preserve zero and the orientation.  

Thus we constructed a homomorphism 
$$\alpha: SL_2(\BZ)\to 
\Aut_T(\cB_1)/T^{2\BZ},
$$ 
where $\Aut_T(\cB_1)$ is the
group of automorphisms of $\cB_1$
fixing $T.$ 
The elements from $T^{2\BZ}=\{T^{2n}\}$ are
considered as inner automorphisms. 

Let  $\tau_+$,
$\tau_-$ be the $\alpha$-images
of the matrices
$\left(\begin{array}{cc}1&1\\0&1\end{array}
\right)$ and $\left(\begin{array}{cc}1&0\\1&1
\end{array}\right).$ 
Then $\sigma=\tau_+\tau_-^{-1}\tau_+$
corresponds to $\left(\begin{array}{cc}0&-1\\1&0
\end{array}\right),$
and $\sigma^2$ has to be the conjugation by 
$T^{2l-1}$ for some $l.$ Similarly, 
$$ \tau_+\tau_-^{-1}\tau_+=T^{2m} \tau_-^{-1}\tau_+\tau_-^{-1}.$$
Using rescaling $\tau_{\pm}\mapsto T^{2m_{\pm}}\tau_\pm$
for $m_-+m_+=m,$ we can eliminate $T^{2m},$ and make 
$0\le l\le 5.$ Generally, $\,l\,$mod $6\,$ is the invariant of the
action, due to Steinberg.
Taking the "simplest" pullbacks for $\tau_{\pm},$
we easily check that $l=0$
and calculate the images of the generators under
$\tau_\pm$ and $\sigma.$    
We arrive at the relations from Proposition \ref{p81}.
\smallskip

{\bf Noncommutative Kodaira-Spencer map.}
Generalizing, let $\e$ be an algebraic, or
complex analytic, or symplectic, or real analytic
manifold, or similar. It may be noncompact and singular.
We assume that there is a continuous family of {\em topological
isomorphisms} $\e\to\e_t$ for manifolds $\e_t$ as $0\le t\le 1,$
and that $\e_1$ is isomorphic to $\e_0=\e.$
The path $\{\e_t\}$ in the moduli space $\m$ of $\e$
induces an outer automorphism $\eps$ of the fundamental
group $\pi_1(\e,\star)$ defined as above. Namely, we take
the image of $\gamma\in \pi_1(\e,\star)$ in 
$\pi_1(\e_1,\star_1)$ for the image $\star_1$ of the
base point $\star\in \e$ and conjugate it by the path
from $\star_1$ to $\star.$ We obtain that
{\em the fundamental group $\pi_1(\m)$} (whatever it is)
{\em acts in $\pi_1(\e)$ by outer automorphisms modulo inner
automorphisms.} 

The above considerations correspond to the case
when a group $G$ acts in $\e$ preserving
a submanifold $D.$ Then  $\pi_1(\m)$ acts in 
$\pi_1^{orb}((\e\setminus D)/G)$ by outer automorphisms.

Another variant
is with the Galois group taken instead of $\pi_1(\m)$
assuming that $\e$ is an algebraic variety over a
field which is not algebraically closed.      

The action of $\pi_1(\m)$ on an {\em individual}
$\pi_1(\e)$  generalizes, in a way, the celebrated
{\dfont Kodaira-Spencer map} and 
is of obvious fundamental importance.
\index{Kodaira-Spencer map}
However calculating the fundamental
groups of algebraic (or similar) varieties, generally
speaking, it is difficult. The main known examples are
concerning the products of algebraic curves and related 
configuration spaces. 
Not much can be extracted from the claim above
without an explicit description of the fundamental
groups.
\medskip

\subsection{Plancherel formulas}
Recall that we have two representations of the double
affine Hecke algebra $\HH\,,$ the polynomial
representation $\cP=\BC[X^{\pm 1}]$ and the representation
$\cF$ in the space of finitely supported functions 
on the set $\BZ_\#$. We constructed two Fourier transforms
$\bS, \bE:\cP\to \cF$ defined by
$$\bS(f)(n_\#)=\lr f\eps_n(x)\mu_k^0(x)\rr_0,\; \;
\bE(f)(n_\#)=\lr \overline{f}\eps_n(x)\mu_k^0(x)\rr_0$$
where $\lr \cdot \rr_0$ denotes the constant term.

The space $\cP$ is equiped with the scalar product
$\lr f, g\rr_0\equal \lr f\bar g\mu_k^0(x)\rr_0.$ Its
$\cF$-counterpart is as follows: 
\begin{align}\label{ffroduct}
&\lr f,g \rr_1\equal \lr f\bar g\mu_k^1\rr_1, 
\hbox{\ where\ }\bar{g}(n_\#),
=\overline{g(n_\#)},\\
&\lr f\rr_1\equal \sum_{n\in \BZ}
f(n_\#),\ \ \mu_1^k(z)\equal \mu_k(z)/\mu_k(0_\#).\notag 
\end{align}
A simple calculation for $n\in \BZ_{>0}$ gives that
\begin{equation}\label{fmuone}
\mu_k^1(n_\#)=\mu_k^1((1-n)_\#)=q^{-k(n-1)}
\prod_{j=1}^{n-1}\frac{1-q^{2k+j}}{1-q^j}.
\end{equation}
Note that $(1-n)_\#=1/2-n_\#$ and hence
$\pi(\mu_k^1)=\mu_k^1$. Also 
$\overline{\mu_k^1}=\mu_k^1$ which makes
the form $\lr \cdot ,\cdot \rr_1$ symmetric.

\begin{theorem}\label{t91}
(Plancherel formula I). For any
$f, g\in \cP,$ 
\begin{equation}\label{f91}
\lr f, g\rr_0=\lr \bS(f), \bS(g)\rr_1=\lr \bE(f),
\bE(g)\rr_1.
\end{equation}
\end{theorem}

{\bf Proof.} Let $H\mapsto H^\star$ denote taking the
adjoint with respect to the scalar
product $\lr \cdot ,\cdot \rr_0$. From Theorem \ref{t62},
$$
X^\star=X^{-1},\, Y^\star=Y^{-1}\, T^\star=T^{-1}\,
q^\star=q^{-1},t^\star=t^{-1}.
$$
It is an anti-automorphism of $\HH\,.$

The theorem follows formally from the following
statements.

(i) For automorphisms  $\sigma$ and $\ep$ from
Theorem \ref{mainth7},
$\sigma \star\sigma^{-1}=\star$, $\ep \star\ep^{-1}=\star.$
Use that $\sigma$ and $\ep$ are homomorpisms of the
group $\cB_q.$ Any group automorphisms
commute with the inversion $g\mapsto g^{-1}$.

(ii) The representations $\cP$ and $\cF$
of the algebra $\HH$ are irreducible. It is true
for generic $q,t.$ We can assume this
because it sufices to check (\ref{f91}) for generic $q,t.$ 

(iii) Normalization: (\ref{t91}) holds for $f=g=1$. \sq\medskip

We introduce the characteristic and delta-functions
 $\chi_m^\#,\, \delta_m^\# \in \cF$ by the formulas 
\begin{equation}\label{fdeldef}
\chi_m^\#(n_\#)=\delta_{mn},\ \,\delta_m^\# (n_\#)=
\delta_{mn}/\mu_k^1(n_\#),
\end{equation}
where $\delta_{mn}$ is the Kronecker delta.

\begin{theorem}\label{t92}
(a) $\bE (\eps_m(x))=\delta_m^\#$,\qquad (b) $\bS
(\eps_m(x))=t^{1/2}T^{-1}(\delta_m^\#).$
\end{theorem}

{\bf Proof.} (a) By the definition of $\bE$ and
the norm formulas (\ref{fnormep}), 
$$\bE(\eps_m(x))=\chi_m^\# \lr \eps_m(x),
\eps_m(x)\rr_0=\delta_m^\#.
$$
Instead of using the norm formulas we can apply the
Plancherel theorem:
$$
\lr \eps_m(x), \eps_m(x)\rr_0=
\lr \chi_m^\#, \chi_m^\#\rr_1 \lr \eps_m(x),
\eps_m(x)\rr_0^2.
$$ 
Since  $\lr \chi_m^\#, \chi_m^\#\rr_1=\mu_k^1(m_\#),$  
we get that $\lr \eps_m(x), \eps_m(x)\rr_0=$
$\mu_k^1(m_\#)^{-1}.$ 

(b) From (a), $\bS(\overline{\eps_m}(x))
=\delta_m^\#$. Then we use part (a) of Proposition \ref{p91}. 
\sq\medskip

\begin{corollary}\label{c92}
Let $\widehat{\eps}(m_\#;x)\equal \eps_m(x).$ 
For any $g\in \cF,$ 
\begin{equation}\label{f92}
\bS^{-1}(g)=t^{-1/2}\lr g\,T^{(1)}(\widehat{\eps}(m_\#;x)) 
\mu_k^1\rr_1,
\end{equation}
where $T^{(1)}$ acts on the first argument of $\widehat{\eps}$
via (\ref{ftsharp}).
\end{corollary}

{\bf Proof.} Using that
 $\bS(\eps_m(x))=t^{1/2}T^{-1}\delta_m^\#$ due to Theorem \ref{t92},
\begin{align}
&\bS(\,t^{-1/2}\lr g\,T^{(1)}(\widehat{\eps}(m_\#;x))\, )= 
t^{-1/2}\lr g\,T(t^{1/2}T^{-1}\delta_m^\# )\mu_k^1\rr_1\notag\\
&=\lr g\,\delta_m^\# \mu_k^1\rr_1=\lr g\,\chi_m^\#\rr_1=g(m_\#).
\hbox{\ \sq }
\end{align} 

We are going to drop the bar-conjugation in the scalar
products:
$$\lr\!\lr f,g \rr\!\rr_0\equal 
\lr fg\mu_k^0(x)\rr_0
\hbox{\ on \ } \cP,\ \
\lr\!\lr f,g \rr\!\rr _1\equal \lr fg\mu_k^1\rr_1 
\hbox{\ on \ } \cF.
$$
Theorem \ref{t92} and Theorem \ref{t81}.
result in the following corollary.

\begin{corollary}\label{c91}
For arbitrary $n,m\in \BZ,$
\begin{align}\label{ftdelta}
&\lr \eps_n(x)\eps_m(x)\mu_k^0(x)\rr_0=
t^{1/2}T^{-1}\delta_n^\# (m_\#),\hbox{\ where\ }\\
&T(\delta_n^\#)=
\frac{t^{1/2}q^{2n_\#}-t^{-1/2}}{q^{2n_\#}-1}\delta_{-n}^\#
-\frac{t^{1/2}-t^{-1/2}}{q^{2n_\#}-1}\delta_n^\#. \notag
\end{align}
In particular, the inequality
$\lr \eps_n(x)\eps_m(x)\mu_k^0(x)\rr_0\ne 0$ implies
$m=\pm n.$ \sq\medskip
\end{corollary}

\begin{theorem}\label{t93}
(Plancherel formula II) For any
$f,g\in \cP,$ we have
\begin{equation}\label{f93}
t^{-1/2} \lr\!\lr \bS(f), T\bS(g)\rr\!\rr_1=
\lr\!\lr f,g\rr\!\rr_0=t^{-1/2}\lr\!\lr\bE(f),
T\bE(g)\rr\!\rr_1.
\end{equation}
\end{theorem}

{\bf Proof.} Let $\diamond$ 
denote the anti-involution corresponding to the scalar
product $\lr\!\lr f,g \rr\!\rr_0.$ It equals $\star\circ\eta:$
$$
X^\diamond =X,\, T^\diamond =T,\, \pi^\diamond =\pi,\,
Y^\diamond =TYT^{-1}.
$$ 
The same anti-involution serves the scalar product
$\lr\!\lr f,g \rr\!\rr_1.$

Following Theorem \ref{t91} (see (i,ii,iii)), we
need to check that 
$$\ \sigma \diamond \sigma^{-1}=
\ep \diamond \ep^{-1}=T^{-1}\diamond T.
$$
This is straightforward. \sq\medskip

Preparing for the inversion theorems,
let us introduce "conjugations" of the transforms
$\bS, \bE:$ 
\begin{align}\label{fnewtrans}
&\bar \bS:f\mapsto \lr f\,\overline{\eps}_n(x)\mu_k^0(x)\rr_0,\\
&\bar \bE: f\mapsto \lr \overline{f\,\eps}_n(x)
\mu_k^0(x)\rr_0.\notag
\end{align}
Theorem \ref{mainth7} states
that for any $H\in \HH$, $f\in \cP,$ 
$$\bS(H(f))=\sigma(H)(\bS(f)),\; \; \bE(H(f))=
\ep(H)(\bE(f)).$$
It is not difficult to find the automorphisms
corresponding to $\bar \bS$ and $\bar \bE$:
$$\bar \bS(H(f))=\sigma^{-1}(H)(\bar \bS(f)),\; \;
\bar \bE(H(f))=\eta \ep(H)(\bar \bE(f)).$$

Following Theorem \ref{masternon},
we get one more Master formula:
\begin{align}\label{f94}
&\lr \overline{\eps_n}(x)\overline{\eps_m}(x)
\widehat{\gamma_-}\mu_k^0(x)\rr\\
&=q^{m_\#^2+n_\#^2-20_\#^2}\,
t^{-1/2}\,(T\overline{\eps_m}(x))(n_\#)\lr 
\widehat{\gamma_-}\mu_k^0(x)\rr.\notag
\end{align}

Algebraically, it is equvalent to
the formula 
$$\sigma^{-1}\tau_+^{-1}\eta (Y)\ =\ \tau_+ (T(\eta(Y)T^{-1}),
$$
resulting from $\sigma^{-2}=T(\cdot )T^{-1}$
and $\eta\tau_{\pm}\eta=\tau_{\pm}^{-1}:$
$$\sigma^{-1}\tau_+^{-1}\eta=
\tau_+\tau_+^{-1}(\sigma^{-2}\sigma)\tau_+^{-1}\eta=
\tau_+\sigma^{-2}(\tau_+^{-1}(\sigma)\tau_+^{-1})\eta$$
$$=\tau_+\sigma^{-2}\tau_-^{-1}\eta=
\tau_+\sigma^{-2}\tau_-^{-1}\eta=\tau_+\sigma^{-2}\eta\tau_-.
$$

\subsection{Inverse transforms} 
Let us summarize what has been done.
We defined four transforms acting from $\cP$ to $\cF:$
$$\bS_0(f)=\lr f\,\eps_n(x)\mu_k^0(x)\rr_0,
\; \; \bE_0(f)=\lr \overline{f}\,\eps_n(x)\mu_k^0(x)\rr_0,$$
$$\bar{\bS}_0(f)=\lr f\,\overline{\eps}_n(x)\mu_k^0(x)\rr_0,\; \;
\bar{\bE}_0(f)=
\lr \overline{f}\,\overline{\eps}_n(x)\mu_k^0(x)\rr_0.
$$
Let us introduce their counterparts 
acting in the opposit direction, from the space $\hat \cF$ of all 
functions on $\BZ_\#$ to $\cP.$ Replacing  
$$
\lr \cdot ,\cdot \rr_0\leadsto \lr \cdot , \cdot \rr_1,\ 
\mu_k^0\leadsto \mu_k^1,
$$
$$\bS_1(f)=\lr f\,\eps_n(x)\mu_k^1\rr_1,\; \; \bE_1(f)=
\lr \overline{f}\,\eps_n(x)\mu_k^1\rr_1,$$ 
$$\bar \bS_1(f)=
\lr f\,\overline{\eps}_n(x)\mu_k^1\rr_1,\; \;
\bar \bE_1(f)=\lr \overline{f}\,\overline{\eps}_n(x)\mu_k^1\rr_1.$$

\begin{theorem}\label{t101} (Inversion theorem)
We have the following identities:
$$\bar \bS_1\circ \bS_0=\id =\bS_1 \circ \bar \bS_0,
\; \; \bar \bS_0\circ \bS_1=\id =\bS_0 \circ
\bar \bS_1,$$
$$\bE_1\circ \bE_0=\id =\bar \bE_1 \circ \bar \bE_0,
\; \; \bE_0\circ \bE_1=\id =\bar \bE_0 \circ \bar \bE_1.$$
\end{theorem}

{\bf Proof.} We know that the automorphisms
corresponding to $\bS$, $\bE$, $\bar \bS$, $\bar \bE$ are
$\sigma$, $\ep$, $\sigma^{-1}$, $\eta \sigma,$ respectively.
This gives that all maps in the theorem are
homomorphisms of $\HH$-modules. Now we use that
$\HH$-modules $\cP$, $\cF$ are irreducible and
the inversion formulas are true for $f=1\in \cP$ or
$f=\delta_0^\#\in \cF$. \sq\medskip

We can obtain Jackson-type
Master formulas from the old ones by formal
conjugating. We also replace $\widehat{\gamma_-}$
by $\gamma=q^{x^2}.$ 

For example, formula (\ref{f94}) results in
\begin{equation}\label{f95}
\lr \eps_m(x)\eps_n(x)\gamma \mu_k^1\rr_1=
q^{20_\#^2-m_\#^2-n_\#^2}\,t^{1/2}T^{-1}
(\eps_m(x))(n_\#)\lr \gamma \mu_k^1\rr_1.
\end{equation}

Formulas (\ref{f81}), (\ref{f82}) read
\begin{equation}\label{f96}
\lr \overline{\eps_m}(x)\overline{\eps_n}(x)\gamma
\mu_k^1\rr_1=q^{20_\#^2-m_\#^2-n_\#^2}\,
\overline{\eps_m}(n_\#)\lr \gamma \mu_k^1\rr_1,
\end{equation}
where $20_\#^2-m_\#^2-n_\#^2=
-\frac{m^2+n^2+2k(|m|+|n|)}{4},$ and
\begin{equation}\label{f97}
\lr \eps_m(x)\overline{\eps_n}(x)\gamma \mu_k^1\rr_1=
q^{20_\#^2-m_\#^2-n_\#^2}\,\eps_m(n_\#)\lr \gamma
\mu_k^1\rr_1.
\end{equation}

The proof of these formulas (including the
convergence of the integrals) is a straightforward 
copy of that
in the compact case. Plancherel formulas  can
also be transformed to the Jackson case 
with ease using Theorem \ref{t101}. 
We would like to mention
here papers \cite{KS1,KS2} devoted to a variant
of the Fourier transform $\bS_1,$ including 
its inversion in the analytic setting (in
contrast to our, essentially algebraic, discussion).

{\bf Comments.} 
(i) The limit $q\to \infty$ of the  transform 
$\bar{\bS}_1: \cF \to
\cP$ is the $p$-adic spherical transform
due to Matsumoto \cite{Mat}. Its inverse  $\bS_0$ is also
compatible with the limit.
In papers \cite{O3,O4}, Opdam developed the Matsumoto
theory of "nonsymmetric" spherical functions towards the
theory of nonsymmetric polynomials. The operator $T_{w_0}$
appears there in the inverse transform in a way 
similar to that of the present paper 
(\cite{O3}, Proposition 1.12).

(ii) As $q\to 1,$ the transform $\bS_0$ becomes the
Hankel transform. One can also get the Harish-Chandra
spherical transform under the same limit by switching from
the variable $x$ to the variable $X=q^x.$ That is, we need to
rewrite all formulas in terms of $X$ and then leave $X$
untouched in the limit. Note that the self-duality of the 
Fourier transform and the Gaussian do not survive 
in this limit as well as in the $p$-adic case.

(iii) We can transform the norm formula 
(\ref{fnormep}) and its variant without conjugation
(\ref{ftdelta}) to the Jackson case:
$$\lr \eps_m(x)\overline{\eps_n}(x)\mu_k^1\rr_1=
\mu_k^1(n_\#)^{-1}\delta_{mn},$$
$$\lr \eps_m(x)\eps_n(x)\mu_k^1\rr_1=
t^{1/2}T^{-1}(\delta_n^\#)(m_\#).$$
However these formulas hold only for $\Rea k<<0.$
In the absence of the nonsymmetric polynomials, $\Rea k <0$
is sufficient. It is a particular case of the Aomoto conjecture,
proved by Ian Macdonald. See \cite{Ao,M4}.

\bigskip
\setcounter{equation}{0}
\section{Generic finite dimensional representations}
Finite dimensional representations 
can appear either for special $k,$ namely, for half-integers,
or for special $q$ (roots of unity). Let us start their 
classification with the case of generic $q.$

\subsection{Generic $q,$ special $k$}
We assume that $q$ is not a root of unity.
When $q^{a/b}$ appear in the claims and/or formulas
for $a/b\in \BQ,$ then we will assume in the next theorem
that $q^{1/b}$ is somehow fixed and $q^{a/b}=(q^{1/(b)})^a.$
In particular, a product $q^a t^b$ can be $1$ but may not be
a nontrivial root of unity. 
 
\begin{theorem} \label{tnegk}
(i) An arbitrary irreducible
finite dimensional representation $V$ is a 
quotient of either the polynomial represenation 
$\cP=\BC[X^{\pm 1}]$ or its image under the automorpisms
$\iota,\,\varsigma_y,\,\iota\varsigma_y,\,$ where 
\begin{equation}\label{fiota}
\iota:\, T\mapsto -T,\, X\mapsto X,\, Y\mapsto Y,
\, q^{1/2}\mapsto q^{1/2},\, t^{1/2}\mapsto t^{-1/2},   
\end{equation}
\begin{equation}\label{fvarsy}
\varsigma_y:\, T\mapsto T,\, X\mapsto X,\, Y\mapsto -Y,
\, q^{1/2}\mapsto q^{1/2},\, t^{1/2}\mapsto t^{1/2}.   
\end{equation}
The $\cP$ is $Y$-semisimple if
and only if $k$ is not a negative integer.
It is irreducible if and only if
$k\neq -1/2-n,\, n\in\BZ_+.$

(ii) Let $k=-1/2-n$ for $n\in\BZ_+.$
The polynomials $e_m$ are well defined for all $m$
and form a basis of $\cP.$ They are multiplied by
$(-1)^m$ under the automorphism of $\HH\,$ 
\begin{equation}\label{fvarsx}
\varsigma_x:\, T\mapsto T,\, X\mapsto -X,\, Y\mapsto Y,
\, q^{1/2}\mapsto q^{1/2},\, t^{1/2}\mapsto t^{1/2}.   
\end{equation}
The values $e_m(-k/2)$ are nonzero
and therefore the polynomials 
$$\eps_m=e_m/e_m(-k/2)\hbox{\ exist\ for\ }
M=\{-2n\le m \le 2n+1\}.
$$
The series $\mu_k^0$ and the pairing $\lr f,g\rr$ on $\cP$
are well defined for $k=-1/2-n.$ The scalar squares
of $e_m$ are nonzero precisely at the same set $M.$
The radical $Rad_0$ of the pairing $\lr f,g\rr$ is 
$\oplus_{m\not\in M} \BC e_m.$ 

(iii) Continuing (ii), $Rad_0=(e_{-2n-1})$ as an ideal in $\cP,$
and the $\HH\, $-module
$\cP/Rad_0$ is the greatest 
finite dimensional quotient of $\cP.$ 
It is the direct sum of the two non-isomorphic
$\HH\,$-submodules of dimension $2n+1$
\begin{align}
&V^\pm_{2n+1}\equal\oplus_{m=1}^{2n+1}\BC \eps_m^\pm 
\hbox{\ mod\ }Rad_0,\
\eps_m^\pm=\eps_m\pm \eps_{-2n-1+m},\\
&V^\pm_{2n+1}\cong \BC[X^{\pm 1}]/(\eps^\mp) \hbox{\ for\ } 
\eps^\pm=\eps_{n+1}\pm \eps_{-n}.\label{fvnegk}
\end{align}
These modules are orthogonal to each other with respect to
$\lr\, ,\, \rr,$ and $V^-_{2n+1}$ is isomorphic to
the $\varsigma_x$-image of 
$V^+_{2n+1}.$

(iv) The representation $V^+_{2n+1}$ is the quotient of
$\cP$ by the radical $Rad$ of the pairing
$\{f\, ,\, g\}=f(Y^{-1})(g)(-k/2)$ from Theorem \ref{t64},
which is $V^-_{2n+1}+Rad_0.$ Respectively, $V^-_{2n+1}$
corresponds to the paring 
$$\{f\, ,\, g\}_-\ =\ f(Y^{-1})(g)\mid_ {X\mapsto-t^{-1}}.
$$
The discretization map 
$$
\chi:
f\mapsto f(z),\ z\in\, \bowtie\,' \equal \, 
\{1/4+m/2\,\mid\, -n\le m\le n\}
$$  
identifies $V^+_{2n+1}$ with 
$F_{2n+1}\equal $ Funct $(\bowtie'),$ where
the action of $\HH\ $ is via formulas (\ref{ftsharp}) from
Theorem \ref{t81}. The $\eps^\pm$ from (ii) are proportional to
$$e^{\pm}=e_{n+1}\pm t^{-1/2}e_{-n}=
X^{-n}\prod_{m=-n}^{n}(X\pm q^{1/4+m/2}).$$
\end{theorem}

{\em Proof.} 
Let $V$ be an arbitrary $\HH\, $-module,
and $v$ a $Y$-eigenvector of weight $\lambda:$
$Y(v)=q^{\lambda}v.$
Following Corollary \ref{c71}, we
construct the chain
\begin{equation}\label{fchain}
v=v_0\stackrel{A_0}{\longrightarrow}v_1
\stackrel{B_1}{\longrightarrow}v_{-1}
\stackrel{A_{-1}}{\longrightarrow}v_2
\stackrel{B_{2}}{\longrightarrow}\ldots,
\end{equation}
where the operation $A_m$ is $q^{-m/2}X\pi$ and the
operation $B_m$ is the application of 
$$t^{1/2}(T+
\frac{t^{1/2}-t^{-1/2}}{q^{-2\lambda_m}-1}).$$
They come from the intertwining operators, so we 
will call them intertwiners too.
The normalization is adjusted to the case 
$v=1\in \BC[X^{\pm 1}],$
when $v_m=e_m,$ $\la_m=-n_{\#}.$ 
Generally, $v_m$ is a $Y$-eigenvector 
of weight 
$$\lambda_{m}=-\lambda -m/2\ \hbox{\ for\ } m>0,\ \  
  \lambda_m=  \lambda -m/2\ \hbox {\ for\ } m\le 0.
$$

Because $q$ is not a root of unity, the chain 
is infinite and  contains only invertible intertwiners
$B_m$ ($A_m$ are always invertible) unless for some $m\in \BZ_+,$

 (a) either  $q^{2\lambda}=q^{-m}$ ("non-existence"), 
 
 (b) or $q^{2\lambda}=t^{\pm 1} q^{-m}$
("non-invertibility").
 
We treat $B_m$ as elements of the nonaffine
Hecke algebra $\bH=\BC+\BC T.$ Note that $B_m$ is
simply $s$  in the degenerate case $k=0,$ and there
are no finite dimensional representations in this case.
 
(i) Let $V$ be a finite dimensional irreducible module. 
Then $k\neq 0$ and
at least one of the conditions (a), (b) must hold
for some $m.$
Applying intevertible intertwiners, we can always
make 
$$
\hbox{either\ \ } (a')\ \lambda=0  \hbox{\ \ or\ \ }
(b')\ \lambda= k/2.
$$

In the second case, we may need to use the automorphisms $\iota,$ 
and $\varsigma_y.$
Let us start with $(b').$

The corresponding chain of intertwiners
results in the weights
\begin{align}
\lambda_0=k/2,\, \lambda_1=-1/2-k/2,\, \lambda_{-1}=1/2+k/2,\\ 
\ldots,\, \lambda_m=-m/2-k/2,\, \lambda_{-m}=m/2+k/2,\ldots,
\ m>0.\notag
\end{align}
Here all  
intertwiners exist and are invertible respectively
if and only if 
$0\neq m/2+k/2\neq -k/2$ for $m\in \BZ_+.$ These inequalities
make $V$ infinite dimensional, so 
we can assume that $k\in -\BZ_+/2.$
Recall that the $B$-intertwiners are applied to
the $v_m$ with positive indices $m.$

Let us check that $\tilde{v}\equal 
(T-t^{1/2})(v)=0.$
Indeed, $\tilde{v}$ is a $Y$-eigenvector of weight 
$-k/2.$ Taking $v_0=\tilde{v},$ the corresponding
chain (\ref{fchain}) produces the weights
$$\lambda_0=-k/2,\,\lambda_1= -1/2+k/2,\, 
\la_{-1}=1/2-k/2,\, -1+k/2,\, 1-k/2,\ldots.
$$ 
The existence and
invertibility of the intertwiners are satisfied respectively as  
$0\neq m/2-k/2\neq k/2$ for $m\in \BZ_+.$ Since $k<0,$
the chain is always invertible and dim$\, V=\infty$
unless $\tilde{v}=0.$
Thus, under assumption $(b')$, $k$ has to be in $-\BZ_+/2,
k\neq 0,$  and $V$ is a quotient of $\cP,$ where the 
covering map sends 
$1\mapsto v.$

Actually  case $(b')$ with integral $k$
is also a part of $(a').$ Indeed, if  
$k=-n$ for natural $n,$ then the chain behaves as follows. 
It loses the existence at $m/2=n/2$ and then
the invertibility at $m/2=n.$ Thus
we can reach $\lambda=0$ using the invertible intertwiners
in this case, which leads to $(a').$

Let us prove that $V$ is always infinite dimensional in
case $(a')$, i.e., when we can find $\lambda=0$ in $V.$
Instead of (\ref{fchain}),
we use its variant with the space 
$V_0=\bH v.$ The latter is an irreducible $2$-dimensional
module over the affine Hecke algebra $\langle T,Y\rangle.$
The operator $Y$ is not semisimple there. This module is
the space of all solutions of $(Y-1)^2(w)=0.$
The operations $A,B$ transform $V_0$ to 
$$
V_m\equal \{w\,\mid\,
(Y-q^{\lambda_m})^2 (w)=0\hbox{\ for\ }\lambda_m=-m/2, m\in \BZ_+\}.
$$
All intertwiners here are invertible
unless $\pm k/2=n/2$ for a positive integer $n.$
Note that the invertibility of all intertwiners 
is possible only for infinite dimensional $V.$ 

Let us consider the case $\{(a'),\ k\in \BZ\ \setminus\{0\}\}.$ 
In this case, the chain from $\lambda_0=0$  
loses the invertibility
exactly once at $m=n.$ The dimension of the next 
(after $V_n$) space $V_{-n}$  becomes $1$ and then remains
unchanged. Therefore $V$ is infinite dimensional too.

We conclude that an irreducible $V$ containing a 
$Y$-eigenvector
of weight $0$ is always infinite dimensional and, 
moreover, it is $Y$-non-semisimple precisely as
$k=\pm\, n$ for a positive integer $n.$ 

In the case under consideration, 
we may assume that $k=-n$ using $\iota,$
so it is covered by $(b')$ and $V$ is a quotient of $\cP.$ 
However the polynomial representation is 
irredicible for such $k.$ 
Indeed, let us assume that there is a
submodule $W\subset \cP.$ Then 
it contains at least one $Y$-eigenvector $v.$
Denoting its weight by $\lambda,$ $W$ 
contains the whole space 
$$\{v\mid (Y-q^{\lambda})^2(v)=0\},$$ 
which is $2$-dimensional.
Using the intertwiners, we can make $\lambda=0$ here. So $1$
belongs to $W$ and we get $W=\cP.$  

(ii) Corollary \ref{c71} provides the existence of
$e_m$ for any half-integral $k.$ The statement about
$\mu_k^0$ is direct from Theorem \ref{t23}.
The other claims follow from formulas
(\ref{fevaluat}) and (\ref{fnorme}) from 
Corollary \ref{c73}:
\begin{equation}\label{fevnew}
e_m(-k/2)=t^{-|m|/2}\prod_{0<j<|m|'}\frac{1-q^jt^2}{1-q^jt},
\end{equation}
where $|m|'=m$ if $m>0,\ $ $|m|'=1-m$ if $m\le 0, \ $ and
\begin{align}\
&\lr e_l\, ,\, e_m\rr=\de_{lm}
\prod_{0<j<|m|'}\frac{(1-q^j)(1-q^j t^2)}
{(1-q^j t)(1-q^{j}t)}.\label{fnormenew}
\end{align}

(iii) Let us show  that
$\cP$ has a unique maximal finite dimensional
quotient and calculate it as $k=-1/2-n$ for $n\in \BZ_+.$
We can construct $e_m$ explicitly using the chain of
intertwiners from $v_0=1.$ 
The intertwiners are
well defined. However one of them, namely, 
$$
B_{2n+1}=t^{1/2}(T+\frac{t^{1/2}-t^{-1/2}} {tq^{2n+1}-1})=
t^{1/2}(T+\frac{t^{1/2}-t^{-1/2}} {t^{-1}-1})=
t^{1/2}(T-t^{1/2}).
$$
is not invertible. Nevertheless
$e_{-2n-1}=B_{2n+1}(e_{2n+1})$ because the leading monomial
of the latter is $X^{-2n-1}.$  
We see that this chain does produce
all $e_m$ regardless of the 
non-invertibility of $B_{2n+1}.$
By the way, this consideration makes the above
reference to Corollary \ref{c71} (concerning the existence
of $\{e_m\}$) unnecessary.
See (ii). 

The polynomials
$e_m$ and $e_{-2n-1+m}$ have coinciding weights  for  
$1\le m\le 2n+1,$ namely, $\ -m_\#=-(-2n-1+m)_\#.$ 
We obtain that the image of the space $J$ linearly generated by 
$e_m,e_{1-m}$ for $m>2n+1$ is zero in $V.$ Indeed,
otherwise dim$\, V$ would be infinite. Moreover,
all eigenvectors in $V$ are images of $\{e_m,e_{1-m}\}$
as $1\le m \le 2n+1$ up to proportionality.
The multiplicities of the images can be either all equal to $1$
or all equal to $2.$ If the multiplicities are all $2,$
then $V$ has to be $\cP/J.$ 

Let us check that $J$ is
an $\HH\, $-submodule of $\cP.$ The $Y$\~invariance
of $J$ is obvious. It is $T$\~invariant, 
since the intertwiners are always 
well defined and 
$$e_{-2n-1}=t^{1/2}(T-t^{1/2})(e_{2n+1})\Rightarrow
T(e_{-2n-1})=-t^{-1/2}e_{-2n-1}.
$$ 
Hence it is $\pi$\~invariant.
The $A$-operations give the $X\pi$\~invariance of $J.$
So it is $X^{\pm 1}$\~invariant, which can also be seen from
the Pieri rules (Corollary \ref{c73}).
Obviously $J=(e_{-2n-1})$ as an ideal, because 
the dimension of the quotient  $\cP/J$ is $2(2n+1).$

If the $Y$-multiplicities are all $1,$ then the dimension
of $V$ becomes $2n+1$ and $V$ has to be irreducible. 
We may set $V=\cP/(e)$ 
for a polynomial
$e=X^{n+1}+\ldots +$ const$\cdot X^{-n}$ 
and a nonzero constant.
Here $e$ must be a $Y$-eigenvector. 
Indeed, there is a linear combination of $e, Y(e)$ 
which is lower than $e,$ unless these vectors are proportional.
Therefore $e=e_{n+1}+ c e_{-n}$ (both have coinciding
weights, namely, $-1/4$).

It is easy
to find that $c=\pm t^{-1/2}$ using the $A,B$\~invariance
of $(e),$ i.e., $e=e^{\pm}$ up to proportionality.
Switching to the $\eps$-polynomials, 
$e$ has to be proportional to
$\eps_{n+1}\pm \eps_{-n}.$ 

However it is easier
to obtain the value of $c$ using (iv), which will be
considered next.
Indeed, if we know that the ideal $V^-_{2n+1}=(e^-)$ is a
nontrivial submodule of $\cP/J,$ then its orthogonal 
complement with respect to $\lr\, ,\, \rr$ is obviously 
$V^+_{2n+1}=(e^+).$
Their sum is direct and they are non-isomorphic as 
$\HH\,$-modules because $e^{\pm}$ are non-proportional.
Note that they are transposed by $\varsigma_x$ 
(acting on polynomials) because $\varsigma_x(e^+)=\pm e^-.$
So $c=\pm t^{-1/2}.$ 

\smallskip
By the way, the modules $V^{\pm}_{2n+1}$
are not isomorphic to each other because otherwise
$\cP/J,$ generated by $1$ as an $\HH\,$-module,
would contain a proper submodule generated by 
$\eps_0^+ +\eps_0^-=2\eps_0=2.$ The latter is impossible.
\smallskip

(iv) The radical $Rad$ of the pairing $\{\ ,\ \}$
contains a $Y$-eigenvector $e'$ if and only if
$e'(-k/2)=0.$ This follows directly from  
the definition (see Theorem \ref{t64}). 
Since the $e$-polynomials form a basis
in $\cP$, a given polynomial $f$ belongs to
the radical if and only if $f(m_\#)=0$ for all $m$
such that $e_m(-k/2)\neq 0.$

The quotient $V=\cP/Rad$ is finite dimensional.
Indeed, it contains all 
$$
\eps_{m}-\eps_{-2n-1+m}\hbox{\ for\ } 
1\le m\le 2n+1.
$$ 
Actually one of them is sufficient for
dim$\,V<\infty$ because $Rad$ is an ideal. 

By the way, it gives 
that  $e_m(-k/2)\neq 0\neq e_{1-m}(-k/2)$
for the  $1\le m\le 2n+1$ and this is the complete list 
of such $m,$ without using the explicit formula 
(\ref{fevaluat}).

As an immediate application, we obtain that
the set $\{m_\#\}$ for $e_m$ with nonzero values at $-k/2$ is
precisely $\,\bowtie'$, and $V_{2n+1}^+$ coincides with
$F_{2n+1}=$ Funct $(\bowtie).$ This 
proves the product formula for $e$ from (iii)
and justifies that $e=e_{n+1}-t^{-1/2}e_{-n}.$
\sq

The theorem gives that for each negative half-integer 
$k=-1/2-n,$
there are $4$ irreducible finite dimensional
representations, namely, 
$\varsigma_x^{\pm 1}\varsigma_y^{\pm 1}(V^+_{2n+1})$  
for all possible combinations of the signs. 
This is different from
the case of the rational DAHA, which has only one such
representation (Theorem \ref{tfdimpr}) as $k$ is a
half-integer. The automorphisms
$\varsigma$ do not have counterparts in the rational limit.

We note that a construction connected
with the finite dimensional quotients  of $\cP$ 
appeared in \cite{DS}.
The authors did not consider the double Hecke algebras but
found some finitely supported measures for orthogonal
polynomials, which is directly related to the theorem.
\medskip

\subsection{Additional series}
If it is not supposed that $q^at^b$ do not represent nontrivial
roots of unity for integral $a,b,$ then the {\dfont additional
series}
\index{additional series ($A_1$)}
of finite dimensional representations
will appear. We continue to 
assume that $q$ is not a root of unity.
Parts (i)-(ii) of the theorem below are due to A.~Oblomkov.
A particular "additional" representation 
was considered in the Appendix of \cite{C4}. 

We continue the above consideration, however $t$ now
cannot be written in the form $q^k$ without certain reservations.
In practical terms,
we need to go back to the multiplicative notation. 
The proof of the statement that finite dimensional 
representations  are
quotients of the polynomial representation up to
the products of the automorphisms
$\varsigma_y$ and $\iota$ remains 
unchanged, so we only need to describe the quotients of the
polynomial representation.

Following the previous section,
we consider the chain (\ref{fchain}):

\begin{equation}\label{fchainn}
v_0=1\stackrel{A_0}{\longrightarrow}v_1=X
\stackrel{B_1}{\longrightarrow}v_{-1}
\stackrel{A_{-1}}{\longrightarrow}v_2
\stackrel{B_{2}}{\longrightarrow}\ldots,
\end{equation}
where the operation $A_{-m}$ is $q^{m/2}X\pi$ and the
operation $B_m$ is the application of 
$$t^{1/2}(T+
\frac{t^{1/2}-t^{-1/2}}{tq^m-1}), \ m>0.$$

Because $q$ is not a root of unity, the chain 
is infinite and  contains only invertible intertwiners
$B_m$ ($A_m$ are always invertible) unless for some $m\in \BZ_+,$

 (a) either  $tq^m=1$ ("non-existence"), 
 
 (b) or $tq^m=t^{-1}$
("non-invertibility").
 
There is only one possible degeneration which was 
not covered by the previous theorem, when
\begin{align}\label{addit}
t=-t^{-n/2}\hbox {\ for\ some\ integral\ } n>0.
\end{align}
In this case, we have no problem with (a).
Also, the polynomials $e_m$ always exist.
Indeed, the $B$-intertwiners always create the 
polynomials with the desired leading term, 
even if they are not invertible.

Because $e_{-n}$ is proportional
to $B_{n}(e_n),$ one has:
$$ T(e_{-n})=t^{1/2}e_{-n}\and s(e_{-n})=e_{-n}.$$
Thus $e_{-n}=X^{-n}+$const$X^{-n+2}\ldots\ +X^{n}.$
Only terms with $X^{-n+2l}$ can appear in the decomposition. 

Using the standard arguments, we conclude that
the space 
$$\BC e_{-n}\oplus \BC e_{n+1}\oplus \BC e_{-n-1}\oplus \ldots$$
is an ideal $(e_{-n})$ and a $\HH\,$-submodule.
We come to the following theorem.

\begin{theorem}\label{tnegka}
(i) Assuming that $q$ is not a root of unity, the 
nonzero finite dimensional
representations of $\HH\, $ which are not described in
Theorem \ref{tnegk} exist only for $t=-q^{\pm n/2}$ as
$0\neq n\in \BZ_+.$ Given $q,t,$
such a representation is unique, has dimension $2n,$
and is irreducible.
For $t=-q^{- n/2},$ it is 
isomorphic to $V_{2n}\equal\p/(e_{-n}).$ 
The automorphisms $\varsigma_x$ and $\varsigma_y$ do not
change its isomorphism class; $\iota(V_{2n})$ corresponding
to $t=-q^{n/2}$ is not a quotient of $\p$ for such $t.$

(ii)
The series $\mu_k^0$ and the pairing $\lr f,g\rr$ on $\cP$
are well defined for such $t.$ The values $e_m(t^{-1/2})$
of $e_m$ and their scalar squares from
(\ref{fevnew}) and (\ref{fnormenew}) are nonzero 
precisely for $m=0,1,-1,\cdots,n.$
The radical $Rad_0$ of the pairing $\lr f,g\rr$ 
coincides with $(e_{-n}),$ so its restriction to $V_{2n}$ 
is well defined.

(iii) The discretization map identifies it with
the space
\begin{align}\label{funadd}
F_{2n}\equal
\hbox{Funct}\Big(t^{-1/2}q^{\frac{1-n}{2}}\cdots, t^{-1/2},
t^{1/2}q^{1/2},\cdots, t^{1/2}q^{\frac{n}{2}}\Big).
\end{align}
Its kernel coincides with $(e_{-n})$ and 
is the radical $Rad$ of the pairing
$\{f\, ,\, g\}=$ $f(Y^{-1})(g)(X\mapsto t^{-1/2})$ 
from Theorem \ref{t64}.
\end{theorem}
{\em Proof.} Claim (i) has been mainly checked. The invariance
of $V_{2n}$ with respect to the substitution $X\mapsto -X$
follows from the the structure of $e_{-n}$ described above.
The map $Y\mapsto -Y$ acts in $V_{2n}$ since the $Y$-spectrum
of $V,$ which is the inverse of the set (\ref{funadd})
in the multiplicative notation, that is,
\begin{align*}
\Big\{t^{1/2}q^{\frac{n-1}{2}}\cdots, t^{1/2},
t^{-1/2}q^{-1/2},\cdots, t^{-1/2}q^{\frac{-n}{2}}\Big\},
\end{align*}
is invariant with respect to the multiplication by $-1.$

Claim (ii) is straightforward. Claim (iii) readily folows
from the descrition of the $Y$-spectrum of $V_{2n}.$  
\sq

\subsection{Fourier transform}
To conclude the consideration of the case
of generic $q$, let us describe the 
action of the projective $PSL(2,\BZ)$ on
$V^+_{2n+1}$ for $k=-1/2-n$ and $V_{2n}$ for $t=-q^{-n/2}.$ 

We will beging with the case $k=-1/2-n.$

First of all, note that the $\sigma$\~invariance
(up to proportionality) of an $\HH\,$-module 
$V$ gives that if $Y(v)=q^\lambda v,$ then 
$X(v')=q^{-\lambda} v'$ for some $v'.$
This holds for $V^+_{2n+1}$ and 
$\varsigma_x\varsigma_y(V^+_{2n+1})$
but is not true for the remaining two, which are
transposed by $\sigma.$ That is why only
$V^+_{2n+1}$ will be considered.
Applying $\varsigma_x\varsigma_y$ we can manage the
second "self-dual" irreducible module. It is not
difficult to extend the  results below
to the direct sum of the remaining two modules. We
leave it as an exercise.

It is convenient to switch from $V^+_{2n+1}$
to its "functional realization" $F_{2n+1}$
from Theorem \ref{tnegk} (iv).

This space has the following scalar product:
\begin{align}
&\lr f,g\rr'\equal 
\lr f\bar g\mu^1\rr',\ \bar g(m_\#)
=\overline{g(m_\#)},\ \lr f\rr'=\sum_{m_\#\in \bowtie\,'}
f(m_\#),  \label{fmuprime}\\ 
&\mu^1(m_\#)=\mu^1((1-m)_\#)=q^{-k(m-1)}
\prod_{i=1}^{m-1}\frac{1-q^{2k+i}}{1-q^i}\hbox{\ as\ } m>0.
\notag
\end{align}

Recall that $(1-m)_\#=1/2-m_\#,$ 
$\pi(\mu^1)=\mu',$ and $\pi(\bowtie')=\bowtie'.$ 
Also $\overline{\mu^1}=\mu^1$ and therefore
the form $\lr \cdot ,\cdot \rr'$ is symmetric.
The operators $X,Y,T,q,t$ are unitary with respect
to this scalar product.

We come to the following truncation of 
Theorem \ref{t81}.

\begin{theorem}\label{masterpr}
(Truncated Master formula) 
(i) Let us introduce the  \index{Gaussian $\ga$} 
{\dfont Gaussian} by the
formula $\gamma(m_\#)=q^{m_{\#}^2}.$ Then  
it induces $\tau_+$  upon the action 
$H\mapsto \gamma H\gamma^{-1}$ on $H\in \HH\,.$
Respectively, the automorphisms of $F_{2n+1}$
\begin{align}\label{fbsbe}
&\bS(f)(m_\#)=\lr f\eps_m\mu^1\rr',\
\bE(f)(m_\#) =\lr \bar{f}\eps_m\mu^1\rr'
\end{align}
induce $\sigma$ and $\ep$
on $\HH\, .$   

(ii) For $l_\#,m_\#\in \,\bowtie',$
\begin{align}
&\lr \eps_l\overline{\eps_m} \gamma\mu^1\rr'
=q^{-\frac{m^2+n^2+2k(|m|+|n|)}{4}}\eps_l(m_\#)
\lr\gamma\mu^1\rr',\ \ \lr\gamma\mu^1\rr'=
\label{fmastpr}\\
\label{fmastprr}
&=\sum_{j=n+1}^{2n+1} q^{\frac{(j-k)^2}{4}}
\frac{1-q^{j+k}}{1-q^{k}} \prod_{i=1}^j
\frac{1-q^{i+2k-1}}{1-q^{i}}=  
q^{1/16}\prod_{i=1}^{n}
(1-q^{1/2-i}).
\end{align}
\end{theorem}
\sq

\medskip

Now let us consider the {\dfont additional series}
\index{additional series ($A_1$)}, namely
the modules $V_{2n}$ for $t=-q^{-n/2}, n\in \BN.$
The analysis is 
straightforward, so we will simply formulate
the counterpart of the previous theorem in this case.

We set 
\begin{equation}\label{qeakn}
q=e^a,\ k=-n/2+\pi i a \for a>0,\ n\in\BN.
\end{equation}
The formulas below are algebraic identities
in terms of $q,t$ and their fractional powers.
This special setting is convenient to continue using
the "exponential" notation.

We use the standard $m_{\#}$
and the functional realization $F_{2n}$
from Theorem \ref{tnegka} in the set:
\begin{align}\label{bowminus}
\bowtie'_-\equal\Big\{\frac{1-n-k}{2}\cdots, -\frac{k}{2},
\frac{k+1}{2},\cdots, \frac{n+k}{2}\Big\}.
\end{align}
     
The Gaussian will be $\gamma(m_\#)=q^{m_{\#}^2}.$
For instance,
$$
\gamma(m_\#)=\exp(\frac{(m-n/2)^2}{4a}+
\frac{(m-n/2)\pi i}{2}-
\frac{\pi^2 a}{4}) \for m>0.
$$ 
The $\eps_m$ are the spherical polynomials, normalized by
$\eps_m(0_\#)=1;\ $
$\mu^1=\mu/\mu(0_\#).$

\begin{theorem}\label{masterpra}
(Additional Series)
Provided (\ref{qeakn}), for 
$l_\#,m_\#\in \,\bowtie_-',$
\begin{align*}
&\lr \eps_l\overline{\eps_m} \gamma\mu^1\rr'_-
=q^{-\frac{m^2+n^2+2k(|m|+|n|)}{4}}\eps_l(m_\#)
\lr\gamma\mu^1\rr'_-,\\
&\lr\gamma\mu^1\rr'_-=
q^{\frac{(j-k)^2}{4}}+
\sum_{j=1}^{n} q^{\frac{(j-k)^2}{4}}\
\frac{1-q^{j+k}}{1-q^{k}}
\prod_{i=1}^j
\frac{1-q^{i+2k-1}}{1-q^{i}},
\end{align*} 
where the Gaussian sum is given by the formulas:
\begin{align}
&\lr\gamma\mu^1\rr'_-=
e^{-\frac{\pi^2 a}{4}}\prod_{i=0}^{l-1}
(1+q^{-i})\for m=2l,\ l\in\BN,\and
\label{fmastprev}\\
&\lr\gamma\mu^1\rr'_-=
\sqrt{2}\ e^{\frac{1}{16a}-\frac{\pi^2 a}{4}}\,
\prod_{i=0}^{l-1}
(1+q^{-i-1/2}),\ m=2l+1,l\in\BZ_+.
\label{fmastprodd}
\end{align}
\end{theorem}
\sq

\medskip
\subsection{Roots of unity $q$, generic $k$} 
Let us study the case
which is in a sense opposite to the previous one.
We assume that $q^{1/2}$ is a
primitive $2N$-th root of unity for $N\ge 1$ and
consider generic $k.$ Till the
end of the subsection, $k\not\in\BZ/2.$
We will continue using the symbol $\BC$ for the field of constants,
although it can be made now $\BQ(q^{1/4},t^{1/2}).$
As above, by generic we mean that
{\em all fractional powers of $q$ must be defined in
terms of the "highest"  
primitive root of unity.}

\begin{theorem}\label{trootgenk}
Let $k\not\in \BZ/2.$
\\
(i) The polynomials $e_m(x)$, $\eps_m(x)$ are well defined 
and constitute a basis of $\cP=\BC[X^{\pm 1}].$
For instance, $e_{-N}=X^N+X^{-N}, e_{-2N}=X^{2N}+X^{-2N}.$
The form
$\lr \cdot , \cdot \rr_0$ is also well defined:
$$
\lr \eps_l(x),\eps_m(x)\rr_0=\delta_{lm}(\mu^1(m_\#))^{-1},
$$
where $\mu^1$ is given by formula (\ref{fmuprime}).
\\
(ii) The vectors $\eps_m(x), m\le -N,$ together with $\eps_m(x),
m\ge N+1,$ form a basis of the radical $Rad_0$ of the pairing
$\lr \cdot , \cdot \rr_0$ on the space $\cP.$
The $\HH'\, $-module $V_{2N}\equal \cP/\Rad_0$ $=
\oplus_{N\ge n\ge -N+1}\BC \eps_m(x)$ is irreducible
of dimension $2N$ with the simple $Y$-spectrum:
$$
\Bigl\{\frac{k+N-1}{2},\ldots,\frac{k}{2},-\frac{k+1}{2},
\ldots,-\frac{k+N}{2}\Bigr\}.
$$
As an ideal, $Rad_0=(e_{-N}),$ so $V_{2N}=\cP/(X^N+X^{-N}).$
\\
(iii) The polynomials $\eps_l-\eps_m$ for
$m_\#=l_\# \mod N$ linearly generate the radical $Rad$ of
the pairing 
$\{f\, ,\, g\}=f(Y^{-1})(g)(-k/2)$ from Theorem \ref{t64}.
The quotient $\cP/Rad$ is irreducible
and isomorphic to the $\HH\,$-module
Funct $(\bowtie_N)$ for the set $\,\bowtie_N\equal\BZ_\#
\hbox{mod\,} N$
of cardinality $4N,$ under the action from Theorem \ref{t81}.
As an ideal, 
$$Rad=(X^{2N}+X^{-2N}-t^N-t^{-N}).
$$
\\
(iv) An arbitrary nonzero irreducible quotient of $\cP$ is
either 
$$V^C\equal \cP/(X^{2N}+X^{-2N}+ C) \for C\neq 2,
$$
or $V_{2N}=\cP/(X^N+X^{-N}),$ which is a quotient of $V^2$
by the submodule isomorphic to the image of $V_{2N}$ under 
$\varsigma_y.$\\
The dimension of $V^C$ is $4N$ and it is
$Y$-semisimple with simple weights constituting the set
$-\bowtie_N.$ It is also $X$-semisimple 
(with simple weights) unless $C=-2.$ 
The module $\cP/Rad$ from (iii) equals $V^{C_t}$ as 
$C_t=t^N+t^{-N}=q^{kN}+q^{-kN};$ $C_t\neq \pm 2$ since
$k\not\in \BZ/2.$
\end{theorem}
{\em Proof.} 
(i) The existence of the series $\mu^0$, the $e$-polynomials, 
and the $\eps$-polynomials
readily follows from Theorem \ref{t23}, Corollary \ref{c71},
and (\ref{fevnew}). One can also involve the Pieri rules
from Corollary \ref{c73} or directly
the chain (\ref{fchain}). Let us calculate $e_{-N},e_{-2N}.$
We get $Y(e_{-N})=q^{N/2+k/2}e_{-N}=-t^{1/2}e_{-N},$
where $q^{N/2}=-1$ because $q^{1/2}$ is a primitive root
of unity. Then $e_{-N}=t^{1/2}(T+t^{-1/2})e_{N}$ thanks
to (\ref{fchain}), so $T(e_{-N})=t^{1/2}e_{-N}.$
The latter immediately gives that $e_{-N}$ is $s$\~invariant
(even), and that 
$$\pi(e_{-N})=YT^{-1}(e_{-N})=-e_{-N},\hbox{\ where\ }
\pi(f(x))=f(1/2-x),$$ 
so $e_{-N}(x+1/2)=-e_{-N}.$ Therefore $e_{-N}$ can contain
only the monomials $X^{N}$ and $X^{-N}.$ 
Combining it with the $s$\~invariance
we obtain the formula for $e_{-N}.$ 
The calculation of $e_{-2N}$ is similar.

(ii) The norm-formulas (\ref{fnormenew}) give the description
of $Rad_0$ as a linear space. Since it is an ideal,
$Rad_0=(e_{-N}).$ As an immediate application, we obtain
that $(e_{-N})$ is an $\HH\, $-submodule. 

Note that it is not difficult to check directly
the $\HH\,$\~invariance of $(e_{-N})$
using the Pieri formula 
for $Xe_{-N}$  together with the above formulas for the
action of $Y,T$ on $e_{-N}.$ 

(iii) The calculation of the radical $Rad$ of the form
$\{\, ,\,\}$ is similar to that for $Rad$ from the previous
theorem. The module $\cP/Rad$ (of dimension $4N$)
has to be irreducible because of the following argument.
If it contains a nontrivial submodule $W,$ then using 
the corresponding chain of invertible intertwiners we 
can find there a $Y$-eigenvector
either of weight $k/2$ or of weight $N/2+k/2,$ i.e.,
either $1$ or $e_{-N}$. The former is impossible.
The latter is impossible too because 
the generator $X^{2N}+X^{-2N}-t^{N}-t^{-N}$ of $Rad$
is not divisible by $e_{-N}=X^N+X^{-N}:$
$$
e_{-N}(q^{k/2})=q^{kN/2}+q^{-kN/2}=q^{-kN/2}(q^{kN}+1)=0
\Leftrightarrow k\in 1/2+\BZ.
$$

(iv) The element $X^{2N}+X^{-2N}$ is central in $\HH\,.$
Indeed, it commutes with $\pi$ and with 
$$
T_i+\frac{t^{\frac{1}{2}}-t^{-\frac{1}{2}}}{X_i^2-1}\hbox{\, for\, }
 i=0,1,\, T_1=T,\,X_1=X,\, T_0=\pi T_1\pi,\,X_0=\pi(X),
$$
generating a proper localization of $\HH\,.$ 
Therefore it is central in $\HH\,$ without the localization.

Hence $V^C$ is an $\HH\,$-module. 
Using the chain of intertwiners
(see (iii)), it is irreducible if and only if
 $X^{2N}+X^{-2N}+ C$ is divisible by $X^N+X^{-N},$ 
so for $C=2$ only. Otherwise its $Y$-spectrum is simple.
Concerning the $X$-spectrum,
the polynomial  $X^{2N}+X^{-2N}+ C$ has simple roots
unless $C=\pm 2.$ The action of $X$ is not semisimple
in $V^{-2}$ (use the affine Hecke algebra $\lr T,X\rr$).
\sq

{\bf Comment.} The form $\lr\, ,\,\rr$ making
the generators $X,Y,T,q,t$ unitary and inducing
$\star$ on $\HH\,$ cannot be 
introduced on $V^C$ for $C\neq 2.$ It follows from
formula (6.26) (before Proposition 6.3) in \cite{C12}.
The form $\{\, ,\, \}$  inducing $\phi$ exists only 
on $\cP/Rad$ from (iii). Note that if
$V_{2N}$ had such a form, then $\star\cdot \phi=\ep$
would be an inner automorphism of $V_{2N}$ and 
Spec$_X=-$Spec$_Y,$ which is impossible.

\bigskip
\setcounter{equation}{0}
\section{Classification, Verlinde algebras} 
In this
section we continue to assume that $q^{1/2}$ is a
primitive $2N$-th root of unity, $N\ge 1.$
However now $k$ will be arbitrary (possibly special).
The actuial field of constants is $\BQ(q^{1/4},t^{1/2}),$
although the next theorem is stated over $\BC.$
Recall, that the $Y$-weights $\lambda$
are defined as follows: $Yv=q^\lambda v,$ $v\neq 0.$

The classification below gives the list of all
possible pairs $\{\la,k\}.$ However in process of
proving the next theorem we give a complete description
of the defining parameters of the irreducible representations.

Note that the consideration of the regular case can be
simplified using the structure of the center of DAHA. 
It makes the parametrization $X\leftrightarrow Y$\~ symmetric.
Using the center makes sense for the special cases
too, although the technique of intertwiners seems
more convenient.

There is a variant of the general theory
for odd $N,$ when one picks $q^{1/2}$ in 
primitive $N$-th roots of unity.
This case will not be considered in this section, 
as well as the equivalent case of
the Little DAHA generated by $\lr X^2,T,Y^2\rr\in \HH$
(see below).
\medskip

\subsection{The list}
\begin{definition}
(i) A number $\lambda\in \BC$  is called 
\index{la@$\lambda$ regular ($A_1$)} {\dfont regular} if
the orbit $O_\lambda\equal $
$\pm \lambda+\BZ/2 \mod N$ is simple, i.e., contains
$4N$ elements. Equivalently, $2\lambda\not\in \BZ/2\mod N.$
Otherwise it is called 
\index{laa@$\lambda$ half-singular ($A_1$)} 
{\dfont half\~singular} as 
$2\lambda\in 1/2+ \BZ\mod N,$ and 
\index{laaa@$\lambda$ singular ($A_1$) } 
{\dfont singular}
as $2\lambda\in \BZ\mod N.$

(ii) An irreducible representation $V$
of $\HH\, $ is called $Y$-{\dfont principal}
\index{principal-special ($A_1$)}
if its dimension is $4N$
and $Y$ has a simple spectrum in it. Otherwise it is called
$Y$-{\dfont special}.
\end{definition}

Note that the $O_\lambda$ can consist of $4N$ or $2N$ elements. 

\begin{theorem}
Let $V$ be an irreducible finite dimensional
$\HH\,$-module with
a $Y$-eigen\-vector of weight $\lambda.$
The following list exhausts all possible
pairs $\{\lambda,k\}.$

(A) Let $\lambda$ be regular. Then $V$ is $Y$-principal
if $k\not\in \pm 2\lambda+\BZ.$  Otherwise,
$k\in \pm 2\lambda+\BZ,$ so
$k\not\in \BZ/2,$ and up to $\iota,$ $\varsigma_y,$
and their product
from (\ref{fiota},\ref{fvarsy}), $V$ is one of the
irreducible quotients of $\cP$ described in
Theorem \ref{trootgenk}. The latter are $Y$-principal
unless $V=\cP/(X^N+X^{-N}).$

(B) Let $\lambda$ be half-singular, i.e., 
$2\lambda\in 1/2+\BZ,$ and either
$k\not\in \BZ/2$ or $k\in \BZ.$
Then $V$ is $Y$-semisimple. The dimension can be 
$2N$ (then the $Y$-spectrum is simple) or 
$4N$ (then all weights are of multiplicity $2).$ 


(C) Let $\lambda$ be singular, i.e., 
$2\lambda\in \BZ,$
and $k\not\in \BZ.$ 
Then the dimension of $V$ is
$4N$ and it is not $Y$-semisimple, unless
$t=1.$

(D) Let $\lambda$ be either half-singular under 
$k\in 1/2+\BZ$
or $\lambda $ singular under $k\in \BZ.$
Then $V$ is a quotient of $\cP$ or its image under the 
automorphisms $\iota,\varsigma_y,\iota\varsigma_y.$
\end{theorem}

{\em Proof.} 
Given a $Y$-eigenvector $v=v_0$ of weight $\lambda=\lambda_0,$
we will apply the chain of intertwiners
from (\ref{fchain}):
\begin{equation}\label{fchainnew}
v_0\stackrel{A_0}{\longrightarrow}v_1
\stackrel{B_1}{\longrightarrow}\cdots{\longrightarrow}v_{-m+1}
\stackrel{A_{-m+1}}{\longrightarrow}v_{m}
\stackrel{B_{m}}{\longrightarrow}\cdots v_{2N},
\end{equation}
to construct $v_{1-m},v_m$ for  $1\le m\le 2N.$ 
Recall that
\begin{align*}
&A_{-m+1}=q^{(m-1)/2}X\pi, \ \ 
B_m=t^{1/2}(T+\frac{t^{1/2}-t^{-1/2}}{q^{-2\lambda_m}-1}),\and \\
&\la_1=-\frac{1}{2}-\la,\la_{-1}=\frac{1}{2}\la,\ldots,
\la_m=-\frac{m}{2}-\la=\la_{-m} \for m>0.
\end{align*}

We can also start here 
with the space $V_0=V(\lambda)\subset V$ of 
all $\lambda$-eigenvectors. Then the composition
$$
K\equal B_{2N}A_{1-2N}\ldots B_2A_{-1}B_1A_0
$$ 
preserves $V_0=V(\lambda)$ if $K$ is well defined.
Similarly, the {\em chain} can be applied to the space
$\widehat{V}_0=\widehat{V}(\lambda)$ of the {\em generalized
eigenvectors}, satisfying $(Y-\lambda)^M(u)=0$ for
sufficiently large $M.$ In this case, however, we need to go back
to the formula 
$$
B_m=t^{1/2}(T+\frac{t^{1/2}-t^{-1/2}}{Y^{-2}-1}).
$$

(A) All intertwiners exist and are invertible
if $k/2\not\in O_\lambda.$ We can set 
$K(v_0)=q^{\kappa}v_0$ for arbitrary $\kappa\in \BC$
and uniquely extend it to the action of $\HH\,$ on
$\oplus_{m=1-2N}^{2N} \BC v_m.$ It is a simple exercise.
This representation is principal. If $k/2$ belongs to
the orbit $O_\lambda,$ then we use the intertwiners to find
a $Y$-eigenvector $v_0$ of weight $\pm k/2\mod N/2$
such that the $B$-intertwiner is zero on it.
Indeed, if $B(v_0)\neq 0,$ then we take it as $v_0.$ 
Using
$\iota$ and $\varsigma_y,$ we can assume that this
weight is $k/2.$
This means that $V$ is a quotient of $\cP.$ Now we can
apply Theorem \ref{trootgenk}.

(B) All intertwiners exist and are invertible.
We can make $\lambda=-1/4.$ 
Then $S=q^{-1/4}X\pi$ creates the weight $-1/2-\lambda=-1/2,$
and therefore preserves the space $V_0$ of all eigenvectors
of weight $-1/4.$ 
The operators $K,S$ act in $V_0$ and satisfy 
the relation $S^2=1, SKS=K^{-1}.$ The irreduciblity
of $V_0$ as a $\{S,K\}$-module is necessary and sufficient
for the irreduciblity of $V.$ Similar to (A), we can 
uniquely extend an arbitrary action of $S,K$ on $V_0$ to a
structure of $\HH\,$-module on $V.$ Thus $V_0$ is either
one-dimensional ($K=\pm 1$, $S=\pm 1$) or two-dimensional.
We obtain the desired result.

(C) The intertwiners are invertible provided we have  their existence.
We can find $v_0$ with $\lambda=0.$ Then $\cH_y v_0$
is a two-dimensional irreducible representation of the 
affine Hecke algebra $\cH_Y= \lr T,Y\rr,$
which coincides with 
$$
\widehat{V}_0=\widehat{V}(0)= 
\{u\,\mid\,(Y-1)^2(u)=0\}.
$$
The action of $Y$ here is not semisimple. 
We take invertible 
$$
L\equal A_{1-N}\ldots B_2A_{-1}B_1A_0
$$ 
instead of $K.$ It sends $\lambda=0$ to $-N/2$ and
$\widehat{V}_0$ to $\widehat{V}_N.$
Recall that $q^{N/2}=-1,$ so
the the intertwiner
$B_N$ acting in $\widehat{V}_N$
is singular. We obtain that $(Y+1)^2=0$ on $\widehat{V}_N$
and $T$ preserves it. If we know the action of $T$ in this space, 
then it is sufficient to reconstruct uniquely the $\HH\,$-action
on $V.$ Indeed, we know the action of 
$Y$ and $X\pi,$ which
sends $\widehat{V}_N$ back to $\widehat{V}_{1-N}.$
The $T$-action depends on one parameter because we can 
conjugate $T$ by the matrices in the centralizer of $Y$  
(the $Y$-action is known). Any choice of this parameter
gives the corresponding action of $\HH\, $ on $V.$
Here we need to exclude the case $k=0,$ which require
somewhat different consideration. 

(D) Similar to the second part of (A), we can find
a $Y$-eigenvector $v_0$ of weight $\pm k/2\mod N/2$
such that the $B$-intertwiner is zero on it. Using
$\iota$ and $\varsigma_y$ if necessary, we obtain a surjection
$\cP\to V.$
\sq 
\smallskip

Actually we have proved more than what was stated
in th theorem.
In cases (A-C), we gave the following complete description of 
{\em all} parameters which determine the irreducible
representations.
\smallskip

\noindent
Given a $Y$-weight $\lambda$, we need either the action of
$K$ (A), or $K$ and $S$ (B)
in the corresponding $Y$-eigenspace. 
\smallskip

\noindent
Respectively, we
need the action of $T$ in 
$\{v\, \mid\, (Y+1)^2(v)=0\}$ for (C) as $k\neq 0$, 
and it can be arbitrary, compatible with the given action of $Y.$
\smallskip

\noindent
Case (D) requires an extension of Theorem \ref{trootgenk}
to integral and half-integral $k,$ 
which we are going to discuss now.
\smallskip

\subsection{Special spherical representations}
First of all, the substitution 
$$T\mapsto -T,\ t^{1/2}\mapsto -t^{1/2}$$ 
identify the polynomial
representations for $t^{1/2}$ and $-t^{1/2}.$ 
It does not act on monomials $X^m,$ so the 
nonsymmetric polynomials $e_{m}$ remain
unchanged, as long as they are well defined.
The spherical polynomials $\eps_m$ for even $m$ do not change
either, $\eps_{m}\mapsto -\eps_m$ for odd $m.$ 

Thus it is sufficient
to decompose $\cP$ upon the transformation $k\mapsto k+N,$
and we can assume that 
$\-N/2\le k <N/2.$ We will also
use the outer involutions $\iota,\varsigma$ of $\HH$ from
(\ref{fiota},\ref{fvarsy}). 

\begin{theorem}\label{tsph}
(i) Let $k\in \BZ/2\mod 2N$ and $-N/2\le k < N/2.$
Then the representations
$V^C=\cP/(X^{2N}+X^{-2N}+C)$ of dimension $4N$
from Theorem \ref{trootgenk}
remain irreducible for $C\neq \pm 2.$
However now $V^C$ for such $C$ becomes $Y$-non-semisimple 
for integral $k\neq 0$ and $Y$-semisimple with the $2$-fold 
spectrum otherwise. 

Moreover, $V^{-2}$ is irreducible for
half-integral $k$, as well as the quotient
$V_{2N}=$$\p/(X^N+X^{-N})$ of $V^{2}$ for integral $k.$
Recall that the kernel of the map $V^2\to V_{2N}$ is 
isomorphic to the image of $V_{2N}$ under 
$\varsigma_y$ sending $Y$ to $-Y.$  

(ii)
If $0<2k<N,$ then either $V^{-2}$ for $k\in \BZ$ and
or $V_{2N}$ for $k\in 1/2+\BZ$ has a unique irreducible 
nonzero quotient
$V_{2n}=\cP/(\eps_{-n})$ of dimension $2n$ for $n=N-2k.$ 
Its $Y$-spectrum is simple.
The eigenvectors are the
images of $\eps_{m}$ for $-n+1\le m\le n,$ 
which are all well defined. It coincides with 
$\cP/Rad$ for the radical of the pairing $\{\, ,\, \};$
\begin{align}\label{TeY}
T(e_{-n})\ =\ -t^{-1/2}e_{-n}\ =\ Y(e_{-n}).
\end{align}

These statements can be extended to $k=0$ when $t=1$ 
if we take $\eps_{-n}=X^N-X^{-N}$ in place of $\eps_{-n}.$

(iii) 
Let $k=-1/2-n$ for integral $0\le n< (N-1)/2.$ Then the module 
$V_{2N}=\cP/(X^N+X^{-N})$ has two non-isomorphic irreducible
the $Y$-semisimple quotients, namely, the
representations from (\ref{fvnegk}):
\begin{align*}
&V^\pm_{2n+1}=\cP/(\eps_{n+1}\pm \eps_{-n}),\\
&\hbox{\ dim}(V^\pm_{2n+1})=|2k|=2n+1,\ 
V^-_{2n+1}=\varsigma_x(V^+_{2n+1}).
\notag
\end{align*}
Here $\eps_m$ are well defined when $-2n\le m\le 2n+1.$
The binomials $X^{n+1}\pm X^{-n}$ must be taken in place of 
$\eps_{n+1}\pm \eps_{-n}$ to extend these statements to
the boundary case $N=2n+1.$ 

The kernel $(\eps_{-2n-1})$ of the map from 
$V_{2N}$ to the direct sum of  $V^+_{2n+1}$ and 
$V^-_{2n+1}$ has dimension
$2N-|4k|$ and is isomorphic to $V_{2N-4|k|}$ from (ii) under
the involution $\iota$ sending $k$ to $-k,$
$T\mapsto -T.$ The vector $e=\eps_{-2n-1}$ satisfies
$T(e)\ =\ -t^{-1/2}e,\ $ $Y(e)= t^{-1/2}e.$
 
(iv) 
In the last case $k\in -1-\BZ_+,$ the module
$V^{-2}$ has a unique irreducible nonzero
quotient $V_{2N+4|k|}$ (of dimension $2N+4|k|$).
It equals $\cP/Rad$ for the radical of the pairing 
$\{\, ,\, \},$ and is isomorphic to $\cP/(e)$ for
$e=\eps_{N}-\eps_{-N-2|k|}$ satisfying relations
(\ref{TeY}).
Here the polynomials $\eps_{N}, \eps_{-N-2|k|}$
is well defined. 

It is also isomoprphic to
the kernel
of the map $V^{-2}\to V_{2N-4|k|}$ from (ii) under
the outer involution $\iota\varsigma_y.$ 

(v) The polynomials $\eps_m$ for $\{m=-2|k|,2|k|+1,\cdots,-N+1,N\}$
exist and their images generate the $Y$-semisimple part 
of $V_{2N+4|k|}$ (of dimension $2N-4|k|$). The corresponding
$Y$-weights are 
$$
\{\la=\frac{|k|}{2},\frac{-|k|-1}{2},\frac{|k|+1}{2},\cdots
\frac{N-1-|k|}{2},\frac{|k|-N}{2}\}.
$$
The rest of $V_{2N+4|k|}$
is the direct sum of $2|k|$ Jordan $2$-blocks
of the total dimension $4|k|$
corresponding to the remaining $2|k|$ weights in the orbit
$O_\la=$
$\{\la=0,-1/2,1/2,-1,\cdots,$$(N-1)/2,-N/2\}.$
Namely, there are two segments of non-semisimple $\la:$
$$\{-\frac{|k|}{2},\frac{|k|-1}{2},\cdots,-\frac{1}{2},0\},\ 
\{\frac{N-k}{2},\frac{|k|-N-1}{2},
\cdots,\frac{N-1}{2},-\frac{N}{2}\}.
$$
\end{theorem}
{\em Proof.} (i) We use that the $X$-spectrum is
simple in $V^C$ and the intertwiners are always
invertible. This readily gives the irreducibility.
The Jordan $Y$-blocks will appear if and only if the $Y$-spectrum
contains $0$ (with a reservation about the degenerate case
$k=0,$ which must be considered separately).
 
The space 
$\widehat{V}(0)=\{v\,\mid\, (Y-1)^2v=0\}$ is
an $\cH_Y$-module which is not  
$Y$-semisimple. So are all other generalized 
$Y$-eigenspaces. This happens precisely for integral $k\neq 0.$   
This argument has been used quite a few times. 
We obtain (i).
\smallskip

Concerning the existence of the $e$-polynomials,
the invertibility of the 
intertwiners is sufficient but not necessary.
The following lemma gives the general construction.
In fact, it has been already used before.
It will be applied a couple of times, especially
in part (iv).

\begin{lemma}\label{epolyn}
Let $V_0=\C,$ $V_1=\C X,\cdots,$  
$$
\widehat{V}_{-m}=B_m\widehat{V}_m,\
\widehat{V}_{m-1}=A_{-m}\widehat{V}_{-m},\cdots,
$$
where $m>0$, $A_{-m}=q^{m/2}X\pi,$
$B_m$ is the restriction
of the $B$-intertwiner 
$B=t^{1/2}(T+\frac{t^{1/2}-t^{-1/2}}{Y^{-2}-1})$ 
to $\widehat{V}_m$ provided that  $q^{2\la_m}\neq 1$
for $\la_m=-m/2-k/2.$ If $q^{2\la_m}=1,$ then we set
$B_m=t^{1/2}T,$ $\widehat{V}_{-m}=\widehat{V}_m+
T\widehat{V}_m.$

The space $\widehat{V}_{\pm m}$ is always one or
two-dimensional, and $(Y-q^{\pm\la_m})^2$ equals $0$ there. 
It contains
a unique $e$-polynomial, however of
degree smaller than $|m|$ if the space is two-dimensional.

If dim$\widehat{V}_{m}=1$
then the corresponding chain of intertwiners
$$A_{1-m}B_{m-1}\cdots B_1A_0(1)$$ 

\noindent
produces the polynomial
$e_{m}\in \widehat{V}_{m}.$ Also, 
$e_{-m}=B_{m}e_{m}$ for such $m>0.$
When dim$\widehat{V}_{m}=2$ then the $e$-polynomials
$e_{m},e_{-m}$ do not exist in $\cP.$

Let us assume that
either 1) $q^{2\la_m}=t$ or 2) $q^{2\la_m}=t^{-1}.$ 
Then dim$\widehat{V}_{-m}=1$ and $e_{-m}\in \widehat{V}_{-m}.$
If here dim$\widehat{V}_{m}=2,$ then respectively 

\centerline{
either\ 1) $(T+t^{-1/2})e_{-m}\ $\hfil or \ 2) $(T-t^{1/2})e_{-m}$} 

\noindent
is nonzero and proportional to
the unique $e$-polynomial in the space 
$\widehat{V}_{m}.$ If dim$\widehat{V}_{m}=1,$
then the vector 1) or 2) is zero. \sq
\end{lemma}

If $V=\cP/J$ is a proper (neither $\{0\}$ nor $\cP$)
irreducible 
$\HH\,$-quotient, then $J$
is an ideal. We may set $J=(e)$ for either
$$(\a):\ 
e=X^l+\ldots +cX^{-l}\hbox{\ \, or\ \,}
(\b):\ e=X^{l+1}+\ldots+ cX^{-l},
$$ 
where $2N>l>0,c\neq 0,$
$$
(\a):\ \hbox{dim}\, V=2l,\ \  
(\b):\ \hbox{dim}\, V=2l+1.
$$
 
In either case, $e$ is an eigenvector of $Y.$ Indeed, 
$Y(e)\in J$ has the same type as $e$ with the same
degrees. A proper
linear combination of $e$ and $Y(e)$ would be lower than
$e,$ unless they are proportional. We set $e=e^0+e^1,$ where
$e^\alpha(-X)=(-1)^\alpha e^\alpha(X).$ Note that this decomposition
is always nontrivial in case $(\b).$ 

Since $Y$ commutes with the $\iota_x$ sending
$X\mapsto -X,$ then $\{e^0,e^1\}$ are $Y$-eigenvectors of
top $X$-degrees 
$$(\a):\ \{l,m\} \hbox{\ as\ } 0<m< l,\hbox{\ for\ odd\ } l-m,\ \ 
(\b):\ \{l+1,l\}.
$$ 
The $Y$-weights of $e^0,e^1$ must coincide mod $ N.$ 
They can be readily calculated:
$$(\a):\ \{(l+k)/2, \pm (m+k)/2\},\ \  (\b):\ 
\{-(l+1+k)/2, (l+k)/2\}.$$
\noindent Here it suffices to know the leading term of $e^\alpha$
with respect to the ordering $\succ:$ $X^{\pm m}\Rightarrow$ 
$\mp (m+k)/2$ for $m>0.$
 
The coincidence mod $ N$ immediately gives that only 
the plus\~sign is possible in $(\a).$ We arrive at the
following relations: 
$$(\a):\ k=-(l+m)/2\mod N,\ \  (\b):\ k=-1/2-l\mod N.
$$
\noindent In either case, $k\in 1/2+\BZ$ if the
decomposition $e=e^0+e^1$ is nontrivial.

We obtain that $(\b)$ results in  claims from (iii) 
from the theorem.
Indeed, $0\le l< 2N,$ and $l$
is nothing   but $n=k+1/2.$ The representations
$V^\pm_{2n+1}$ exist and remain irreducible in this case.
\smallskip

We now need to examine $(\a).$ We will
impose this condition till the end of the proof.

The argument above, which ensured the $Y$\~invariance
of $e,$ can be used for $T$ as well. It preserves the type
of $e$ (note that it does not  hold in case $(\b)$). 
We obtain that 
$T(e)=\pm t^{\pm 1/2} e.$ Therefore $e=e_{-l}.$ 
Indeed, $e-e_{-l}$ must be proportional to $e_{m}$ 
for positive $m\le l.$ However, such small $e_m$ 
are never  eigenvectors of $T$ since
the chain of intertwiners remains invertible in this range.

We have four subcases of $(\a)$:
\begin{align}
&(\alpha)\ T(e)=t^{1/2}e=Y(e),\ \ \ \ \,
(\beta)\ T(e)=t^{1/2}e=-Y(e),\notag
\\
&(\gamma)\ T(e)=-t^{-1/2}e=Y(e),\ 
(\delta)\ T(e)=-t^{-1/2}e=-Y(e).\notag
\end{align}   

Concerning $(\alpha,\beta),$ we obtain that $s(e)=e,$ so $e$ is
even and the decomposition $e=e^0+e^1$ is trivial.
Then $(l+k)/2=k/2\mod N$ for $(\alpha)$ and  
$(l+k)/2=k/2+N/2\mod N$ for $(\beta).$ The former equality
results in $l=0\mod 2N,$ which is impossible. 
As for $(\beta):$
$$
\pi(e)=YT^{-1}(e)=-e\ \Rightarrow\ e(x+1/2)=-e(x)
\ \Rightarrow\ e=X^N+X^{-N}.
$$
Therefore $(\beta)$ implies
that the quotient $\cP/(X^N+X^{-N})$ is an $\HH\,$-module
and is irreducible unless it is covered by (iii) of the theorem,
which was already considered.
\smallskip

$(\gamma)$ We obtain  $(l+k)/2=N/2-k/2\mod N$
and $l=N-2k\mod 2N.$ Recall that $e=e_{-l}.$
If $2k>0,$ we arrive exactly exactly at (ii) of the theorem.
In this case, $l=n.$
All polynomials $e_m$ are well defined and
have nonzero $e_m(-k/2)$ in this case.
Otherwise $2k<0.$  

$(\gamma_1)$ If $2k<0,$ then $k$ cannot be a half-integer.
Indeed, $e_{-l}(-k/2)=0$ in this case (use the
evaluation formula), as well as for the generator of
the ideal for $V^+_{2n+1}$ in case (b). We see that these two 
ideals together do not generate the whole $\cP.$
Since $V$ is irreducible, it has to coincide with $V^+_{2n+1},$
which is impossible because we assumed that $e$ is
in the form (a). Therefore $k$ must be an integer.

$(\gamma_2)$ We consider the integers $-N/2\le k<0,$
examine
the chain 

$$\widehat{V}_0=\C,
\widehat{V}_{1}=\C X,\widehat{V}_{-1},\cdots,
\widehat{V}_m,\cdots $$
from Lemma \ref{epolyn}. See also (\ref{fchainnew}).
The following takes place ($m>0$):
\smallskip

\noindent
0) the $\widehat{V}$-spaces are one-dimensional
till $m=|k|;$

\noindent
1) the intertwiner is trivial ("non-existent")
at $m=|k|,$ and dim$\widehat{V}_m =2$ for
$|k|< m\le 2|k|;$

\noindent
2) the $B_m$ kills $1\in \widehat{V}_m$ at $m=2|k|$, and
dim$\widehat{V}_m =1$ for $2|k|<m\le N;$ 

\noindent
3) $B_m$ becomes proportional to $(T+t^{-1/2})$ at $m=N,$
$e_{-N}=X^N+X^{-N},$ and dim$\widehat{V}_m =1$ for
$N< m\le N+|k|;$

\noindent
4) it is trivial ("non-existent") again at $m=N+|k|,$
and dim$\widehat{V}_m =2$ for
$N+|k|< m\le N+2|k|;$

\noindent
5) the $B$-intertwiner kills $e_{-N}$ at $m=N+2|k|,$
and its image becomes proportional to $e_{-N-2|k|},$
which has the same $Y$-eigenvalue as $e_{N}.$
\medskip

The polynomials $\eps_{-n-2|k|}$ and $\eps_N$ exist
and the difference $e=\eps_N-\eps_{-N-2|k|}$ belongs
to the radical $Rad.$ The lemma gives
that the vector $(T+t^{1/2})e_{-N-2|k|}$
is proportional to $e_{-N},$ as well as 
$(T+t^{1/2})e_{N}.$ Therefore the difference $e$ satisfies
$(T+t^{1/2})e=0.$ 
\smallskip

The lemma gives
that between 2) and 3), the $e_m$ exist and their
images linearly generate the $Y$-semisimple part of $V.$
Otherwise, there will be Jordan $2$-blocks with respect to $Y.$
Indeed, we obtain the  
2-dimensional irreducible representation
of $\cH_Y$ in the space of the corresponding
generalized eigenvectors at step 1). Then we shall apply the
intertwiners to this space
(the weights will go back) and eventually will
obtain the two-dimensional generalized 
eigenspace for the starting
weight $\la=-|k|/2.$ 

The intertwiner at 2) will then make the latter space
one-dimensional (i.e., $Y$-semisimple). 
It will remain one-dimensional
until 3). After 3), we obtain the Jordan blocks similar to  
steps 0)\~ 2). 

It readily results in the irreducibility of $V.$ 
We arrive at (v).
\medskip

Let us consider the last subcase. 

$(\delta)$ Now $l=2N-2k.$ If $2k>0,$ then $e_{-l}$ belongs
to the ideal $(e_{-n}).$ Indeed, its image in $V_{2n}$ is zero 
since there are no eigenvectors of weight $-k/2$ in this module.
That contradicts to the irreduciblity. 
Negative $k$ are impossible because $l<2N$ by assumption. 
So this case is empty.
\smallskip

Concerning the "duality" claims involving $\iota$ 
from (iii),
it holds because $e=e_{-2n-1}$ generating
the kernel $(e)$ of the map $\p\to $
$V_{2n+1}^+\oplus V_{2n-1}^-$  
satisfies the relation $(\delta):$
$T(e)=-t^{-1/2}e=-Y(e).$ Recall that there is
a reservation about the boundary value $2|k|=N.$
 
As for the duality (ii)$\leftrightarrow$(iv), the vector
$e=e_{-n}=e_{-N-2k}$ generating the 
kernel of the map $\p\to V_{2N-4k}$
satisfies the relation $(\gamma):$
$T(e)=-t^{-1/2}e=Y(e).$
The same relation holds for $e=\eps_{N}-
\eps_{-N-2|k|}$ from (iv). Note the reservation about
the boundary value $k=0$ in this case.
\sq 

\smallskip
\subsection{Perfect representations} 
We are going the find out which irreducible modules
are $PGL_2(\BZ)$\~invariant. It is clear that the
exceptional representations (ii) - (iv) from the theorem
are invariant because they can be distinguished by the
dimensions. It is also not difficult to describe 
$PGL_2(\BZ)$\~invariant modules from scratch
without using the classification.
Motivated by the Verlinde algebras
\cite{Ve} (see also \cite{KL2}) we will
require more structures.

\begin{definition}
We say that a finite dimensional irreducible
$\HH$-module $V$ is 
\index{perfect representation ($A_1$)} {\dfont perfect}, 
we also call it 
\index{nonsymmetric Verlinde algebra} 
\index{Verlinde algebra ($A_1$)} 
a {\dfont nonsymmetric Verlinde algebra}
if the following conditions hold:

(a) $V$ is \index{spherical irrep ($A_1$)} {\dfont spherical}, 
i.e., there exists a surjection
$\cP\to V;$

(b) it is $\{ \tau_\pm , \ep \} -$invariant,
which means that there are pullbacks of the $\HH$- automorphisms
$\tau_\pm$ and  $\ep$ to  $V$ satisfying the
relations 
$$\tau_+\tau_-^{-1}\tau_+=
\tau_-^{-1}\tau_+\tau_-^{-1},\  \ep^2=1,\
\ep \tau_+=\tau_-\ep;
$$

(c) and also $X$-\index{pseudo-unitary module} 
{\dfont pseudo\~unitarity}: there exists a 
nondegenerate form
$\lr \cdot , \cdot \rr$ on the space $V$ such that
corresponding anti-involution of $\HH$ is $\star$ and
$\lr e, e\rr \ne 0$ for any $X$-eigenvector $e\in V$.
\end{definition}

Note that condition (b) implies that one can
replace $X$-eigenvectors by $Y$-eigenvectors in (c).
Also $V$ is $X$-semisimple and $Y$-semisimple simultaneously.
By  condition (a), the module $V$ has a
structure of a commutative algebra (since it is a quotient
of the commutative algebra $\cP$). Obviously,
$V^{sym}\equal \{ v\in V| Tv=t^{1/2}v\}$ is a
subalgebra of $V,$ a generalization of the Verlinde
algebra.

\begin{proposition}\label{GREATSI}
Let $k\in \BZ/2, |k|<N/2.$ The notation is from Theorem \ref{tsph}.
The greatest $\si$-invariant quotient of $\cP$
is $V^2$ for half-integral $k,$ and $V^{-2}$
for integral $k.$ Recall that
the module $V^2$ is always an extension of
$V_{2N}$ by its $\varsigma_y$-image $\varsigma_y(V_{2N}).$

There are exact sequences
\begin{align*}
&0\to\iota\varsigma_y(V_{2N+4k})\to V^{-2}\to V_{2N-4k}\to 0
\for k\in \Z_+,\\
&0\to \iota(V^{+}_{2k}\oplus V^{-}_{2k})
\to V_{2N} \to V_{2N-4k}\to 0
\for k\in 1/2+\Z_+.
\end{align*} 
The arrows must be reversed for $k<0:$ 
\begin{align*}
&0\to \iota\varsigma_y(V_{2N-4|k|})\to V^{-2}\to 
V_{2N+4|k|}\to 0 \for k\in -1-\Z_+,\\
&0\to \iota(V_{2N-4|k|})\to V_{2N}\to
V^{+}_{2|k|}\oplus V^{-}_{2|k|}\to 0 
\for k\in -1/2-\Z_+.
\end{align*}

The modules $V^{\pm}_{2|k|},$ $V_{2N-4|k|},$ 
$V_{2N+4|k|}$ are irreducible and $\si$ -invariant. The latter 
module is
$Y$ -non-semisimple. The other three are semisimple.
\end{proposition}
{\em Proof.} We use Theorem \ref{tsph},
the structure of the $Y$-spectrum
in the polynomial representation, and that the dimensions
of the Jordan $Y$-blocks are no greater than $2.$
\sq

We obtain that up to $\iota,\varsigma,$
there are {\em three} 
different series
of $\si$-invariant spherical representations at roots
of unity, namely, $V_{2N-4k}$ (integral $N/2>k>0$),
$V_{2|k|}$ (half-integral $-N/2<k<0$), and 
$V_{2N+4|k|}$ (integral $-N/2<k<0$). 
The latter is non-semisimple.

If $k=1,$ the subalgebra $V_{2N-4}^{sym}$ of 
dimension $N-1$ is
isomorphic to the usual {\dfont Verlinde algebra}
\index{Verlinde algebra ($A_1$)} of
$\widehat{\mathfrak{sl}_2}$ of level $N$ (central charge $=N-2$).
The symmetric polynomials 
$p_m^{(1)}$ are the classical characters of finite
dimensional representations of $\mathfrak{sl}_2.$
In $V_{2N-4}^{sym},$ these characters are
considered as functions at roots of unity.

\rmk
1) Recently a non-semisimple variant of the Verlinde algebra
appeared
in connection with the fusion procedure for
the $(1,p)$ Virasoro algebra, although the connection
with the fusion is still not justified. 
Generally speaking, the fusion procedure for the 
Virasoro-type algebras and the so-called $W$-algebras
can lead to non-semisimple Verlinde algebras. At least,
there are no reasons to expect the existence of a
positive hermitian inner product there like the Verlinde
pairing for the conformal blocks, because the corresponding
physics theories are expected {\em massless}.
Presumably it is connected with $V_{2N+4|k|}$ 
for $k=-1.$

2) The latter module and its multi-dimensional
generalizations are also expected to be connected
with the important problem of describing the
complete tensor category of the representations
of the Lusztig quantum group at roots of unity.
The Verlinde algebra, the symmetric part of
$V_{2N-4|k|}$ for $k=1$, describes the so-called
{\em reduced category}, in a sense, corresponding to the
Weyl chamber. The nonsemisimple modules of type 
$V_{2N+4|k|}$ are supposed to appear in the so-called
{\em case of the parallelogram}.  
\sq 

\begin{theorem}\label{t103}
The perfect represenations $V$ are exactly
the cases 
\begin{align}
&(a):\  2k\in \BZ_+, \, 0< k<N/2,\hbox {\ and\ }\notag\\
&(b):\  k=-1/2-n,\, n\in \BZ,\, 0\le n<N/2\notag
\end{align}
from Theorem \ref{tsph}.
Both $V$ are quotients of $\cP$ by the radical $Rad$
of the form $\{ f,g\}=f(Y^{-1})(g)(-k/2).$
They are 
isomorphic to Funct$(\,\bowtie'\,)$:
$$(a) \ \bowtie'\,=\Bigl\{ \frac{k-N+1}{2},
\frac{k-N+2}{2}, \ldots , -\frac{k}{2}, \frac{k+1}{2}, \ldots
\frac{N-k}{2}\Bigr\},
$$  
$$(b)\ \ \bowtie'\,=\Bigl\{
\frac{1/2-n}{2}, \frac{3/2-n}{2}, \ldots, -\frac{1}{4}, 
\frac{1}{4},\frac{3}{4},
\ldots, \frac{1/2+n}{2}\Bigr\}.
$$
The action of $\HH\,$ is defined here by the formulas from 
Theorem \ref{t81}.
The invariant form is (\ref{fmuprime}):
\begin{align}
&\lr f,g\rr'\equal 
\lr f\bar g\mu^1\rr',\ \bar g(m_\#)
=\overline{g(m_\#)},\ \lr f\rr'=\sum_{m_\#\in \bowtie\,'}
f(m_\#),  \label{fmuprimnew}\\ 
&\mu^1(m_\#)=\mu^1((1-m)_\#)=q^{-k(m-1)}
\prod_{i=1}^{m-1}\frac{1-q^{2k+i}}{1-q^i}\hbox{\ as\ } m>0.
\notag
\end{align}
\end{theorem}
{\em Proof.} We first check that
$V=\cP/\Rad$ if $V$ is perfect. Using $\lr\,,\,\rr$
and $\ep$ on $V,$ we introduce the form
$\{ \cdot , \cdot \}_\ep$ on $\cP$ which is
the pullback of the form $\lr f, \ep g\rr.$
This form is symmetric since $\ep$ commutes
with $\star.$ The corresponding anti-involution
is $\,\ep\, \star=\phi$, i.e.,
the same as for the form $\{ \cdot, \cdot \}$.
A proper linear combination $\{ \cdot, \cdot \}_\dag$ of
$\{ \cdot, \cdot \}$ and $\{ \cdot, \cdot \}_\ep$
satisfies $\{ 1,1\}_\dag=0$. Therefore it is zero identically
and  $\{ \cdot, \cdot \},$  $\{ \cdot, \cdot \}_\ep$
are proportional. This gives the desired. 

Second, using the chain of intertwiners
from (\ref{fchain}) and (\ref{fchainnew}),
we obtain that all $\eps_m$ for $m_\#\in -\bowtie'$
are well defined and their images in $V$
linearly generate an irreducible $\HH\,$-submodule $V'$. 
Indeed, the polynomial $e_{-n}$ exists,
$e_{-n}(-k/2)=0,$ and therefore it
belongs to Rad. Equivalently, one
may check that the formulas from Theorem \ref{t81}
define an $\HH\,$-action on Funct$(\,\bowtie'\,).$

The orthogonal complement of $V'$ in $V$
with respect to $\lr\,,\,\rr$ (or
$\{ \cdot , \cdot \}_\ep$) intersect $V$ by zero.

Third, the operator $X$ has a simple spectrum
on $V$. The existence of the $X$-pseudo-unitary
pairing guarantees that $X$ is semisimple.
Since $V$ is a cyclic module over
$\BC[X^{\pm 1}],$ each eigenvalue of $X$ has a unique
Jordan block. Therefore the spectrum of $X$ in $V$ is simple.
Applying $\ep,$ the same holds for $Y.$ For instance,
the image of $1$ in $V$ can belong either to $V'$ or 
to its orthogonal complement constructed above.
This is impossible because
it is a generator of $V.$

Fourth, we can define $\tau_+$ as the multiplication by
the restricted Gaussian $\gamma(m_\#)=q^{m_\#^2}$ in the
realization $V=$Funct$(\,\bowtie'\,).$ The automorphism
$\ep$ in $V$ is the Fourier transform $\bE:$
$\widehat{f}(m_\#) =\lr \bar{f}\eps_m\mu^1\rr'$
from  (\ref{fbsbe}).\sq

Note that the action of $\tau_+,\sigma$ in
$V_{2n}^{sym}$ for $n=N-2k$ and the integral $k$
was introduced and calculated for the first time in
\cite{Ki} on the basis of the interpretation
of Rogers' polynomials for the integral positive $k$
discovered by Etingof and Kirillov.
Verlinde considered only $k=1.$

We will finish this section with the 
following generalization
of Theorem \ref{t32}. We follow Theorem \ref{masterpr},
where negative half-integres $k$ were considered.
See also \cite{C7}.

\begin{theorem}\label{t103r}
Let $0<2k<N,$ $l_\#,m_\#\in \,\bowtie'.$ Then
\begin{align}
&\lr \eps_l\overline{\eps_m} \gamma\mu^1\rr'
=q^{-\frac{m^2+n^2+2k(|m|+|n|)}{4}}\eps_l(m_\#)
\lr\gamma\mu^1\rr',\label{fmastposit}\\
&\lr\gamma\mu^1\rr'=
\prod_{j=1}^k\frac{1}{1-q^j}\sum_{j=0}^{2N-1}q^{j^2/4}
\hbox{\ \ where\ } k\in \BZ,\label{fmastnat}\\
& \sum_{j=0}^{2N-1}q^{j^2/4}=(1+\imath)\sqrt{N} \hbox{\ as\ } 
q^{1/4}=\exp(\pi\imath/2N),\notag
\end{align}
\begin{align}
& \lr\gamma\mu^1\rr'=
2q^{1/16}\prod_{j=1}^{k-1/2}
\frac{1}{1-q^{1/2-j}}\hbox{ \ for\ } k=1/2+\BZ.\label{fmasthalf}
\end{align}
\end{theorem}
\smallskip

\rmk
To conclude, let us note an interesting connection of
these formulas with Theorem \ref{masterpra} describing
the {\dfont additional series}
\index{additional series ($A_1$)}
for non-cyclotomic $q$
and $t=-q^{-n/2}.$
The formulas (\ref{fmastprodd}) and  (\ref{fmastprev})
actually can be used to deduce 
(\ref{fmastposit}), and the other way
round. Indeed, setting in the last formulas
$k=N/2-n/2$ for $n\in \BN$ we obtain the
relation $t=-q^{-n/2}.$ If $N$ is sufficiently big than
the corresponding $q$ is actually "generic".
The only thing we need to establish the connection 
is rewriting the Gaussian sums using 
Corollary \ref{c31}:
\begin{align*}
&\lr\gamma\mu^1\rr'=
q^{n^2/4}\prod_{j=k+1}^n(1-q^j)&\mbox{\ if\ }&N=2n,\\
&\lr\gamma\mu^1\rr'=
q^{n^2/4}(1+q^{(2n-1)/4})
\prod_{j=k+1}^n(1-q^j)&\mbox{\ if\ }&N=2n+1.
\end{align*}
Then we can chose big $N$ providing the desired parity.
The {\em additional} representation $V_{2n}$
becomes $V_{2N-4k}$ for such $N.$ 
This is essentially what was done in \cite{C4}
(Appendix) in a particular case.
\sq 


\bigskip
\setcounter{equation}{0}
\section{DAHA and the $p$-adic theory}
This section is devoted to double  Hecke
algebras in the general setting and their connection
with the classical $p$-adic spherical transform.
In contrast to  previous sections, we state the
results without proofs. The paper \cite{C10} contains
a complete theory.

\subsection{Affine Weyl group} Let $R\subset \BR^n$ be a
simple reduced root system, $R_+\subset R$ the
set of positive roots, $\{ \alpha_1, \ldots,
\alpha_n\} \subset R_+$ the corresponding set
of simple roots, and $\theta \in R_+$ 
the maximal {\em coroot}. 
We normalize scalar product on $\BR^n$ by
the condition $(\theta,\theta)=2,$
so $\theta$ belongs to $R$ as well.
For any $\alpha \in R,$ its dual is 
$\alpha^\vee $ $ =\frac
{2\alpha}{(\alpha,\alpha)} $ $=\alpha/\nu_{\alpha}.$ 
So $\nu_\alpha=(\alpha,\alpha)/2=1,2,3.$

The affine roots are 
$$
R^a=\{\tilde{\alpha}=[\alpha,\nu_\alpha j]\},\ j\in \BZ.
$$
Note the appearance of $\nu_\alpha$ here. It is because of our
nonstandard choice of $\theta.$
We identify nonaffine roots $\alpha$ 
with $[\alpha,0]$
and set $\alpha_0=[-\theta,1].$
For $\tilde{\alpha} \in R^a,$ 
$s_{\tilde{\alpha}} \in $ $ \Aut(\BR^{n+1})$ is the
reflection 
\begin{align}\label{fwaffine}
&s_{\tilde{\alpha}}([x,\zeta])=[x,\zeta]-
2\frac{(x,\alpha)}{(\alpha,\alpha)}\tilde{\alpha},
\hbox{\ and\ }\notag\\ 
&W=\lr s_\alpha\, \mid \, \alpha \in R\rr,\  
W^a=\lr s_{\tilde{\alpha}}\, \mid \, \tilde{\alpha} \in R^a\rr.
\end{align}

We set $s_i=s_{\alpha_i}.$ It is well-known 
that $W^a$ is a Coxeter group with the  generators
$\{s_i\}$. 

Let $\omega_i \subset \BR^n$ be the
fundamental weights: $(\omega_i ,\alpha^\vee_j)=\delta_{ij}$
for $1\le i,j\le n,$  $P=
\oplus_{i=1}^n\BZ \omega_i^\vee$ the weight
lattice, and  $P_+ \equal \BZ_+\omega_i$ the
cone of dominant weights. We call $\hat W\equal W\lsmash P^\vee$
the {\em extended affine Weyl group}:
$$
wb([x,\zeta])=[w(x),\zeta-(b,x)]\hbox{\ for\ } b\in P,x\in \BR^n.
$$

{\em The length function} 
$l: \hat W\to \BZ_+$ is given by the formula 
\begin{equation}\label{falength}
l(wb)=
\sum_{\begin{array}{c}\alpha \in R_+\\
w(\alpha^\vee)\in R_+\end{array}}|(b,\alpha^\vee)|+
\sum_{\begin{array}{c}\alpha \in R_+\\
w(\alpha)\in -R_+\end{array}}|(b,\alpha^\vee)+1|,
\end{equation}
where $w\in W$, $b\in P$. Let $Q =
\oplus_{i=1}^n\BZ \alpha_i$ be the root
lattice. Then $W^a=W\lsmash Q \subset \hat W$ is a normal
subgroup and $\hat W/W^a=P /Q.$ Moreover,
$\hat W$ is the semidirect product $\Pi \lsmash W^a,$ where
$$
\Pi =\{ \pi\in \hat W\,\mid\, l(\pi)=0\}=
\{ \pi\in \hat W\,\mid\, \pi:\{\alpha_i\}\mapsto
\{\alpha_i\}\, \},\ 0\le i\le n.
$$
It is isomorphic to $P/Q$ and acts naturally on
the affine Dynkin diagram for $R^\vee$ 
with the reversed arrows. It is not
just the standard Dynkin diagram
because of our choice of $\theta.$
We set $\pi(i)=j$ as $\pi(\alpha_i)=\alpha_j.$ 

\subsection{Affine Hecke algebra} 
We denote it by $\cH.$  
Its generators are $T_i$ for  $i=0, \ldots, n$ and
 $\pi \in \Pi.$ The relations are
\begin{align}\label{ftt}
&\underbrace{T_iT_jT_i\ldots}_{\mbox{order of}\; s_is_j}
= \underbrace{T_jT_iT_j\ldots}_{\mbox{order of}\; s_is_j},\ \,
\pi T_i\pi^{-1}=T_{\pi (i)},\\
&(T_i-t^{1/2})(T_i+t^{1/2})=0 \hbox{\ for\ }
0\le i\le n,\ \pi\in \Pi.\notag
\end{align}
If $\hat w=\pi s_{i_l}\ldots s_{i_1}\in \hat W$ is a reduced
expression, i.e., $l(\pi^{-1}w)=l,$ 
we set $T_{\hat w}=\pi T_{i_l}
\ldots T_{i_1}$. The elements
$T_{\hat w}$ are well defined and form a basis of $\cH$.

Let $G$ be the adjoint split $p$-adic 
simple group corresponding to $R$. Then $\cH$ is 
the convolution algebra of compactly supported
functions on $G$ which are
left-right invariant with respect to the Iwahori subgroup $B,$
due to Iwahori and Matsumoto. Namely, $T_i$ are the 
characteristic functions of the double cosets
$Bs_iB,$ where we use a natural embedding
$W\to G.$ Generally, it is not a homomorphism.

To be more exact, the $p$-adic quadratic equations
are in the form $(T_i-1)(T_i+t_i)=0$ for
the standard normalization of the Haar measure. Here
$t_i$ may depend on the length of $\alpha_i$ 
and are given in terms of the cardinality
of the residue field. We  
will stick to our normalization of $T,$ and assume that
the parameters $t$ coincide to simplify formulas.
We will also  use
$$\delta_{\hat w}\equal t^{-l(\hat w)/2}T_{\hat w},
$$ 
which satisfy the quadratic equation with $1.$

Let $\Delta$ be the left regular representation
of $\cH$. In the basis $\{\delta_{\hat w}\},$ the representation
$\Delta$ is given by
\begin{equation}\label{fthecke}
T_i\delta_{\hat w}=\left\{ \begin{array}{cc}t^{1/2}
\delta_{s_i\hat w}& \mbox{if\ } l(s_i\hat w)=l(\hat w)+1,\\ 
t^{-1/2}\delta_{s_i\hat w}+
(t^{1/2}-t^{-1/2})\delta_{\hat w}& 
\mbox{if\ }
l(s_i\hat w)=l(\hat w)-1,\end{array}\right.
\end{equation}
and obvious relations $\pi \delta_{\hat w}=\delta_{\pi\hat w},$
where $\pi\in \Pi,\ \hat w\in\hat W.$

The {\em spherical representation}
appears as follows.
Let 
$$\delta^+\equal (\sum_{w\in W}t^{l(w)})^{-1}
\sum_{w\in W}t^{l(w)}\delta_w\in \cH.
$$ 
One readily checks
that $T_i\delta^+=t^{1/2}\delta^+$ for $i=1, \ldots, n,$
and $(\delta^+)^2=\delta^+$. We call $\delta^+$ 
the $t$-symmetrizer.
Then
$$\Delta^+\equal \Delta \delta^+=\oplus_{b\in P}\BC 
\delta^+_{b},\ \delta^+_{\hat w}=\delta_{\hat w}\delta^+,
$$
is an $\cH$-submodule of $\Delta.$

It is nothing but Ind$_H^\cH
(\BC_{t^{1/2}}),$ where $H\subset \cH$ is the subalgebra
generated by $T_i, i=1, \ldots, n$ and $\BC_{t^{1/2}}$
is the one-dimensional representation of $H$ defined by
$T_i\mapsto t^{1/2}.$

Due to Bernstein, Zelevinsky, and Lusztig
(see e.g., \cite{L}), we
set $Y_a=T_a$ for $a\in P_+$ and extend it to the whole $P$
using $Y_{b-a}=Y_bY_a^{-1}$ for dominant $a,b.$
These elements are well defined and pairwise commutative. 
They form
the subalgebra $\cY\cong \BC[P^\vee]$ inside $\cH.$ 

Using $Y,$ one can omit $T_0.$ Namely, 
the algebra $\cH$ is generated by $\{T_i,i>0,\, Y_b\}$
with the following relations:
\begin {align}
&T^{-1}_iY_b T^{-1}_i=Y_b Y_{\alpha_i}^{-1}\hbox{\ if\ }  
(b,\alpha^\vee_i)=1, \label{ftyt}
\\
&T_iY_b=Y_b T_i \hbox{\ if\ } (b,\alpha^\vee_i)=0,\  1 \le i\le  n.
\notag
\end{align}

The PBW-theorem for $\cH$ gives that 
the spherical representation
can be canonically identified with $\cX,$ as $\delta^+$
goes to $1.$ The problem is to calculate this
isomorphism explicitly. It will be denoted by $\Phi.$ 
We come to the definition of Matsumoto's spherical 
functions:
$$
\phi_{a}(\Lambda)=\Phi(\delta^+_a)(Y\mapsto\Lambda^{-1}) ,
\ a\in P.$$
Here by $Y\mapsto \Lambda^{-1},$ we mean that $\Lambda_b^{-1}$ 
substitutes for $Y_b.$
See \cite{Mat}. 
In the symmetric ($W$-invariant) case, the spherical
functions are due to Macdonald.

By construction,
$$\phi_a=t^{-l(a)}\Lambda_a^{-1}=t^{-(\rho^\vee,a)}
\Lambda_{a}^{-1}\hbox{\ for\ } a\in P_+,\ 
\rho^\vee\equal (1/2)\sum_{\alpha>0}\alpha^\vee.
$$
The formula $l(a)=(\rho^\vee,a)$ readily results
from (\ref{falength}).
So the actual problem
is to calculate $\phi_a$ for non-dominant
$a.$ 
\medskip

{\bf Example.} 
Consider a root system of type $A_1$. In
this case, $P=\BZ\omega,$ where $\omega\in \BZ$ is the fundamental
weight, $Q=2\BZ,$ $\omega=\pi s$ for $s=s_1.$
We identify $\Delta^+$ and $\cY,$ so $\delta^+=1.$ 
Letting $Y=Y_\omega,T=T_1,$ we get $Y_{m}\equal Y_{m\omega}=Y^m,$ and
$$\phi_m\equal \phi_{m\omega}=
t^{-m/2}\Lambda^{-m},\hbox{\ for\ } m\ge 0. 
$$ 
Note that $TY^{-1}T=Y$ and $\pi=Y T^{-1}.$ 

Let us check that
\begin{equation}\label{fxphi}
\Lambda\phi_{-m}=
t^{1/2}\phi_{-m-1}-(t^{1/2}-t^{-1/2})\phi_{m+1},\ m>0.
\end{equation}
Indeed, $\phi_{-m}=t^{-m/2}(T\pi)^m(1)
\mid_{Y\mapsto \Lambda^{-1}}$ and $Y^{-1}\phi_{-m}=$
\begin{align}
&=t^{-m/2}(T^{-1}\pi)(T\pi)^m(1)=
t^{-m/2}(T-(t^{1/2}-t^{-1/2}))\pi(T\pi)^m(1)\notag\\
&=t^{-m/2}(T\pi)^{m+1}(1)-
-t^{-m/2}(t^{1/2}-t^{-1/2})(\pi T)^m\pi(1)\notag\\
&=t^{1/2}\phi_{-m-1}(Y^{-1})
-t^{-m/2}(t^{1/2}-t^{-1/2})(\pi T)^m
(t^{-1/2}\pi T)(1)\notag\\
&=t^{1/2}\phi_{-m-1}(Y^{-1})-
(t^{1/2}-t^{-1/2})\phi_{m+1}(Y^{-1}).\notag
\end{align}

\subsection{Deforming $p$-adic formulas} 
The following chain of theorems represents 
a new vintage of the classical theory. 
We are not going to prove
them here. Actually all claims which are beyond the 
classical theory of affine Hecke algebras can be
checked by direct and not very difficult
calculations, with a reservation about 
Theorems \ref{tinverse} and \ref{tgeps}.

\begin{theorem}\label{t111}
Let $\xi \in \BC^n$ be a
fixed vector and let $q\in \BC^*$ be a fixed scalar. We
represent $\hat w=bw,$ where $w\in W, b\in P$.
In $\Delta_q^\xi\equal \oplus_{\hat w\in \hat W}
\BC\delta_w^\xi,$ 
the formulas
$\pi \delta_{\hat w}^\xi =\delta_{\pi \hat w}^\xi$ and
$$T_i\delta_{\hat w}^\xi=\left\{ \begin{array}{cl}
\frac{t^{1/2}q^{(\alpha_i,w(\xi)+b)}-t^{-1/2}}
{q^{(\alpha_i,w(\xi)+b)}-1}\delta_{s_i\hat w}^\xi -
\frac{t^{1/2}-t^{-1/2}}{q^{(\alpha_i,w(\xi)+b)}-1}
\delta_{\hat w}^\xi&\mbox{if}\ i> 0, \\
\frac{t^{1/2}q^{1-(\theta,w(\xi)+b)}-t^{-1/2}}
{q^{1-(\theta,w(\xi)+b)}-1}\delta_{s_0\hat w}^\xi -
\frac{t^{1/2}-t^{-1/2}}{q^{1-(\theta,w(\xi)+b)}-1}
\delta_{\hat w}^\xi&\mbox{if}\ i=0\end{array}\right.$$
define a representation of the
algebra $\cH,$ provided that all denominators are nonzero,
i.e., $q^{(\alpha,b+\xi)}\neq 1$ for 
all $\alpha\in R,\, b\in P.$ 
\end{theorem}

The regular representation $\Delta$ with the
basis $\delta_{\hat w}$ is the limit of representation
$\Delta_{q}^\xi$ as $q\to \infty,$ provided that
$\xi$ lies in the fundamental
alcove:
\begin{equation}\label{falcove}
(\xi,\alpha_i)>0 \hbox{\ for\ } i=1,
\ldots, n,\ (\xi,\theta)<1.
\end{equation}

We see that the representation $\Delta_q^\xi$
is a flat deformation of $\Delta$ for such $\xi.$
Moreover, $\Delta_q^\xi\cong \Delta.$
This will readily follow from the next theorem.

Note that taking $\xi$ in other alcoves, we
get other 
limits of $\Delta_{q}^\xi$ as $q\to \infty.$
They are isomorphic to  
the same regular representation, however the formulas
do depend on the particular alcove.
We see that the regular representation has rather many
remarkable systems of basic vectors. They are not 
quite new in the theory of affine Hecke algebras, but 
were not studied systematically.

\begin{theorem}\label{tdaha}
(i) We set  $X_b(\delta^\xi_{\hat w})=
q^{(\alpha_i,w(\xi)+b)}\delta^\xi_{\hat w}\,$ for $\hat w=bw,$
where we use the notation
$X_{[b,j]}=q^j X_b.$ 
These operators have a simple spectrum 
in $\Delta^\xi_q$ under the conditions of the
theorem and
satisfy the relations dual to (\ref{ftyt}):
\begin {align}
&T^{-1}_iX_b T^{-1}_i=X^{-1}_b X_{\alpha_i}\hbox{\ if\ }  
(b,\alpha^\vee_i)=1, \label{ftxt}
\\
&T_iX_b=X_b T_i \hbox{\ if\ } (b,\alpha^\vee_i)=0,\  1 \le i\le  n,
\notag\\
&\hbox{and\ moreover,\ }\pi X_b \pi^{-1}=
X_{\pi(b)} \hbox{ \ as\ } \pi\in \Pi.\notag
\end{align}

(ii) Introducing the double affine Hecke algebra
$\HH\, $ by (\ref{ftt}) and (\ref{ftxt}), 
the representation $\Delta^\xi_q$ is nothing 
but the induced representation
$$\Ind_{\cX}^{\HH\,}(\BC\delta^\xi_{id}),\ \cX\equal \BC[X_b],
$$
which is isomorphic to $\Delta$ as an $\cH$-module.
\end{theorem}

{\bf Comment.}
The operators $X_b$ are in a way the coordinates
of the Bruhat-Tits buildings corresponding to the
$p$-adic group $G.$ In the classical theory, we use
only their combinatorial variants, namely, the distances
between vertices, which are integers.
The $X$-operators clarify dramatically
the theory of the $p$-adic spherical Fourier transform,
because, as we will see, they are the "missing" 
Fourier-images of the $Y$-operators. Obviously, the
$X_b$ do not survive in the limit $q\to \infty,$
however they do not collapse completely. Unfortunately
the Gaussian, which is $q^{x^2/2}$ as $X_b=q^{(x,b)},$
does.
\sq\medskip

The theorem is expected to be directly connected
with \cite{HO2} and via this paper with \cite{KL1}.
Here we will stick to the spherical representations. 
The following theorem establishes a connection
with $\Delta^+.$ 

\begin{theorem}\label{t113}
We set $q^\xi=t^{-\rho},$ i.e.,
$q^{(\xi,b)}\mapsto t^{-(\rho,b)}$ for all
$b\in P.$ The corresponding representation 
will be denoted by $\widetilde{\Delta}.$ 
It is well defined for generic $q,t.$
For any $b\in P,$ let $\pi_b$
be the minimal length representative in the set
$\{bW\}.$ It equals $\,bu_b^{-1}$ for the length-minimum
element $u_b\in W$ such that $u_b(b)\in -P_+.$
Setting $\delta_b^\#\equal 
\delta_{\pi_b}$ in $\widetilde{\Delta},$
the space
$\Delta^\#\equal 
\oplus_{b\in P}\BC \delta_b^\#$ is an
$\HH\,$-submodule of $\widetilde{\Delta}.$ 
It is isomorphic to $\Delta^+$ as an $\cH$-module.
\end{theorem}

The representation $\Delta^\#$ is described by the
same formulas from Theorem \ref{t111}, which 
vanish {\em automatically} on $s_i\pi_b$ not in the
form $\pi_{c}$ thanks to the special choice of $q^\xi.$
It results directly from the following:
$$s_i\pi_b=\pi_c\Leftrightarrow 
(\alpha_i,b+d)\neq 0,\hbox{\ where\ } (\alpha_i,d)\equal \delta_{i0}.
$$
Here $c=b-((\alpha_i,b+d)\alpha_i^\vee$
for $\alpha_0^\vee\equal -\theta.$

We define the action $(\!(\ )\!)$ of $\hat W$ on $\BR^n$
by the formulas $wa(\!(x)\!)=w(a+x).$
The above $c$ is $s_i(\!(b)\!).$
This action is constantly used in the theory 
of Kac-Moody algebras. It is very convenient 
when dealing with $\Delta^\#.$ Note that
$$
\pi_{b}(\!(c)\!) =bu_b^{-1}(\!(c)\!) =u_b^{-1}(c)+b
\hbox{\ for\ } b,c\in P.
$$

Let us calculate the formulas from Theorem \ref{t111}
upon $q^\xi\mapsto t^{-\rho}$ as
$q\to \infty.$ This substitution changes the consideration
but not too much:
\begin{equation}\label{fdeltas}
T_i\delta_b^\#=\left\{ \begin{array}{cl}
 t^{1/2}\delta_{s_i(\!(b)\!)}^\#
&\mbox{\ if\ } (\alpha_i,b+d)>0, \\
t^{-1/2}\delta_{s_i(\!(b)\!)}^\# + 
(t^{1/2}-t^{-1/2})\delta_{b}^\#\ 
&\mbox{\ as\ } (\alpha_i,b+d)<0.
\end{array}\right.
\end{equation}
Otherwise it is zero.
The formulas $\pi \delta_{b}^\# =$ 
$\delta_{\pi(\!(b)\!)}^\#$ hold for arbitrary
$\pi\in \Pi, b\in P.$

Since the calculation is different from that for generic 
$\xi,$ it is not surprising that (\ref{fdeltas}) 
do not coincide
with (\ref{fthecke}) restricted to $\hat w=b$ 
and multiplied on the right by the $t$-symmetrizer $\delta^+.$
The representations $\lim_{q\to\infty} \Delta^\#$ and
$\Delta^+$ are equivalent, but the $T$--formulas with
respect to the limit of the basis 
$\{\delta^\# \}$ are 
different from those in terms of the classical basis
$\{\delta_b^+=\delta_b\delta^+\}.$ 

\subsection{Fourier transform}
In the first place,
Macdonald's nonsymmetric polynomials generalize
the Matsumoto spherical functions. 
We use $\pi_b=bu_b^{-1}=$ $\hbox{Min-length}\, \{bw,\,w\in W\}.$

\begin{theorem}\label{tpolyn}
(i) Let $\cP$ be the representation of the double
affine Hecke algebra $\HH\,$ in the space of
Laurent polynomials $\cP=\BC[X_b]:$ 
\begin{align}\label{ftonx}
&T_i\  = \  t_i ^{1/2} s_i\ +\ 
(t_i^{1/2}-t_i^{-1/2})(X_{\alpha_i}-1)^{-1}(s_i-1),
\ 0\le i\le n,\\
&X_b(X_c)=X_{b+c},\ \pi(X_b)=X_{\pi(b)},\, \pi\in \Pi,
\hbox{\ where\ }
X_{[b,j]}=q^jX_b.\notag
\end{align}

(ii) For generic $q,t,$
the polynomials $\eps_b$ are uniquely defined from
the relations:
\begin{align} \label{fepsb}
&Y_a(\eps_b)=t^{(u(\rho),a)}q^{-(b,a)}\eps_b,\
\hbox{\ where\ } \pi_b=bu \hbox{\ for\ } u\in W \\
&\eps_b(t^{-\rho})=1 \hbox{\ where\ } 
X_b(t^{-\rho})=t^{-(b,\rho)}.\notag
\end{align} 

(iii) Setting
$X_b^*=X_{-b},\ q^*=q^{-1},$  $ t^*=t^{-1},$
the limit of $\eps^*_b(X\mapsto \Lambda)$ as $q\to \infty$
coincides with $\phi_b(\Lambda)$ for $b\in P.$
\end{theorem}

{\bf Comment.}
Note that the $\Delta^\xi$--formulas from
Theorem \ref{t111} are actually the evaluations
of (\ref{ftonx}) at $q^\xi.$ To be more exact, there
is an $\HH\,$-homomorphism from $\cP$ to the $\HH'$-module
of functions on $\Delta^\xi.$  
For instance,  Theorem \ref{t111} can be
deduced from Theorem \ref{tpolyn}.
The formulas for the polynomial representation of the
double affine Hecke algebra $\HH\,$ are nothing but
the Demazure-Lusztig operations in the affine setting.
\sq\medskip

We are going to establish a Fourier-isomorphism
$\Delta^\#\to \cP$ and 
generalize the Macdonald--Matsumoto
inversion formula. We use the constant term functional on Laurent
series and polynomials  denoted by $\lr\,\rr.$  

The first step is to make both representations unitary using
\begin{align}
&\mu\ =\ \prod_{\tilde{\alpha} \in R^a}
\frac{1-X_{\tilde{\alpha}}}{1-tX_{\tilde{\alpha}}},\ 
\mu^0=\mu/\lr\mu\rr,
\label{fmugen}\\
&\mu^1(\pi_b)=\mu(t^{-\pi_b(\!(\rho)\!)})/\mu(t^{-\rho}),
\ \pi_b=bu_b^{-1}.\notag
\end{align}
Here we treat $\mu$ as a Laurent series to define $\mu^0.$
The coefficients of $\mu^0$ are rational functions in terms 
of $q,t.$
The values $\mu^1(\pi_b)$  are rational function in terms
of $q,t^{1/2}.$ 

The corresponding pairings are 
\begin{align}
&\lr f\, ,\, g\rr_{\,pol}=\lr f\, T_{w_0}w_0(g(X^{-1})\, \mu^0\rr,
\ f,g\in \cP,\notag\\
&\lr \sum f_b\delta^\#_b ,\sum g_b\delta^\#_b\rr_{\,del}=
\sum (\mu^1(\pi_b))^{-1}\,f_b g_b.\notag
\end{align} 
Here $w_0$ is the longest element in $W.$ Note that
the element $T_{w_0}^2$
is central in the nonaffine Hecke algebra $H.$ Both pairings
are well defined and symmetric. Let us give the formulas for the 
corresponding anti-involutions:
\begin{align}
T_i\mapsto T_i,\, X_b\mapsto X_b,\, T_0\mapsto T_0,\,
Y_b\mapsto T_{w_0} Y^{-1}_{w_0(b)}T_{w_0}^{-1} \ \ \ 
&\hbox{\ in\ }\Delta^+,\notag\\
T_i\mapsto T_i, \, Y_b\mapsto Y_b,\, 
T_0\mapsto T_{s_\theta}^{-1}Y_{\theta},\,
X_b\mapsto T_{w_0}^{-1}X^{-1}_{w_0(b)}T_{w_0} 
&\hbox{\ in\ }\cP,\notag
\end{align}
where $1\le i\le n,\, b\in P.$

\begin{theorem}\label{tinverse}
(i) Given $f=\sum_b f_b\delta^\#_b \in \Delta^\#,$ we
set $\widehat{f}=\sum_b\, f_b\,\eps_c^*\in \cP,$
where $X_b^*=X_{b}^{-1}, \ q^*=q^{-1},t^*=t^{-1}.$
The inversion of this transform is as follows:
\begin{equation}\label{fhinv}
f_b\ =\ t^{l(w_0)/2}(\mu^1(\pi_b))^{-1}
\lr\, \widehat{f},\, \eps^*\rr_{\,pol}.
\end{equation}

(ii) The Plancherel formula reads as 
\begin{equation}\label{fhplan}
\lr f\,,\, g\rr_{\,del}\ =\ 
t^{l(w_0)/2}\lr \widehat{f}\,,\, \widehat{g}\rr_{\,pol}.
\end{equation}
Both pairings are positive definite over $\BR$ if
$\,t=q^k,\,q>0$ and $k>-1/h$ for the Coxeter number 
$h=(\rho,\theta)+1.$

(iii) The transform $f=\sum_b f_b\delta^\#_b\mapsto
\widetilde{f}= \sum_b f_b^*\delta^\#_b$ is an involution:
$\widetilde{(\widetilde{f})}=f.$ To apply it for the second
time, we need to replace
$\eps_b^*$ by the corresponding $\delta$-function
$\sum_c \eps_b^*(\pi_c)\mu^1(\pi_c)\delta^\#_c.$
\end{theorem}

Recall that $\eps^*_b$ becomes the Matsumoto spherical
function $\phi_b$ in the limit $q\to \infty$
upon the substitution $X\mapsto \Lambda.$ 
It is easy
to calculate the limits of $\mu^0$ and $\mu^1(\pi_b).$
We come to a variant of the Macdonald--Matsumoto formula.
Claim (iii) has no counterpart in the $p$-adic theory.
Technically, it is because the conjugation $\ast$ sends 
$q\mapsto q^{-1}$ and is not compatible with the limit
$q\to \infty.$ It is equivalent to the non$-p$-adic
self-duality $\eps_b(\pi_c)=\eps_c(\pi_b)$ of the
nonsymmetric Macdonald polynomials.
The following theorem also has no 
$p$-adic counterpart because the Gaussian is missing.

\begin{theorem} \label{tgeps}
We set $\gamma(\pi_b)=
q^{(\,\pi_b(\!(k\rho)\!)\,,\,\pi_b(\!(k\rho)\!)\,)/2},$
where  $t=q^k,$ use $\nu_\alpha=(\alpha,\alpha)/2,$
and $\rho=(1/2)\sum_{\alpha>0}\alpha.$
For arbitrary $b,c\in P,$ 
\begin{align}\label{fhatmu} 
&\lr \eps_b^*\, , \eps_c^*\,\gamma\rr_{del}\ =\
\gamma(\pi_0)^2\gamma(\pi_b)^{-1}\gamma(\pi_c)^{-1}  
\eps_c^*(\pi_b)\lr 1\, ,\,\gamma\rr_{del},\\
&\lr 1\, ,\, \gamma\rr_{del} = 
(\sum_{a\in P} \gamma(\pi_a))\, 
\prod_{\alpha\in R_+}\prod_{ j=1}^{\infty}\Bigl(
\frac{1- t^{(\rho,\alpha)}q^{j\nu_\alpha}}
{1-t^{(\rho,\alpha)-1}q^{j\nu_\alpha}}\Bigr). 
\end{align}
\end{theorem}

\subsection{The $A_1$-case} 
Let us use the Pieri formula (\ref{fpierin})
from Corollary \ref{c73} for $X^{-1}$ and $m\ge 0$
upon the conjugation $\ast$ :
\begin{align}
&X\eps^*_{-m}=\frac{t^{-1/2-1}q^{-m-1}-t^{1/2}}
{t^{-1} q^{-m-1}-1}\eps^*_{-m-1}-
\frac{t^{-1/2}-t^{1/2}}{t^{-1} q^{-m-1}-1}\eps^*_{1+m}.
\label{fpierinc}
\end{align}
Under the limit $q\to\infty,$ we get exactly (\ref{fxphi})
for $\phi_{-m}(\Lambda)\mapsto \eps^*_{-m}(X):$ 
\begin{equation}\label{fxphix}
X\eps^*_{-m}=
t^{1/2}\eps^*_{-m-1}-(t^{1/2}-t^{-1/2})\eps^*_{m+1}.
\end{equation}
The generalization is Theorem \ref{tpolyn}, (iii).

Also $\lim_{q\to \infty}\mu_k^0(x)=
\frac{1-X}{1-tX}$. So the limit of the pairing $\lr \, ,\, \rr_{pol}$
reads as  
$$
\lr f,g\rr^\infty_{pol}=
\lr f\, T(g)\frac{1-X}{1-tX}\rr,\ \, f,g\in \cP.
$$

Substituting $m$ for $m\omega$ 
in the indices: $\Delta^\#=$ $\oplus_{m}\BC \delta_m^\#.$
In the limit $q\to \infty,$ the
operators $T$ and $\pi$ act here as follows: 
$$T\delta_m^\#=t^{1/2}\delta_{-m}^\#,\;
T\delta_{-m}^\#=t^{-1/2}\delta_m^\# +
(t^{1/2}-t^{-1/2})\delta_{-m}^\#,$$
and
$\pi \delta_m^\#=\delta_{1-m}^\#$ for $m\ge 0.$

Concerning $\lr \cdot, \cdot \rr_{del},$ we use
(\ref{fmuone}) for $m>0:$ 
\begin{align}\label{fmuonelim}
&\mu^1(m_\#)=\mu^1((1-m)_\#)=t^{-(m-1)}
\prod_{j=1}^{m-1}\frac{1-t^2q^{j}}{1-q^j}\\
=&\prod_{j=1}^{m-1} \frac{tq^{j/2}-t^{-1}q^{-j/2}}
{q^{j/2}-q^{-j/2}}\ \to \ t^{m-1} \hbox{\ as\ }
q\to\infty.  \notag
\end{align}
Here $\mu^1(m_\#)=\mu^1(\pi_{m\omega}).$
Therefore the limit of this scalar product is simply
$$\lr \delta_m^\#, \delta_n^\#\rr_{del}^\infty=
\delta_{mn}
\left\{ \begin{array}{cl}t^{1-m}&\mbox{if}\; m>0,\\
t^m&\mbox{if}\; m\le 0.\end{array}\right.
$$

{\bf Comment.}
The latter inner product 
is different from the classical one,
which is calculated as follows. We define the pairing 
$$(T_{\hat w}, T_{\hat u})=\hbox {\ Constant\ Term\ }
( T_{\hat w}T_{\hat{u}^{-1}} )$$
on the affine Hecke algebra $\lr Y^{\pm 1},T\rr,$
where the constant term is with respect to the
decomposition via $T_{\hat w}.$
It is simply $\delta_{\hat{u},\hat{w}}.$
Then we switch to
$\delta_{\hat w}=t^{-l(\hat w)/2}T_{\hat w}$ and finally
calculate $(\delta_m^+,\delta_n^+),$
which is $\delta_{mn}t^{-|m|}/(1+t).$
\sq\medskip

The Fourier transform $f\mapsto \widehat{f}$
from Theorem \ref{tinverse} is compatible with the
limit $q\to \infty.$ The inversion
formula (\ref{fhinv}) and the Plancherel formula
(\ref{fhplan}) survive as well. We get a minor reformulation
of the Matsumoto formulas. Upon symmetrization, 
$T$ disappears from the inversion formula and we
come to the Macdonald inversion.

What is completely missing in the limit is (\ref{fhatmu}). 
One of the main applications of the double Hecke algebra
is adding the Gaussian to the classical $p$-adic theory.
Technically, one does not need $\HH\, $ to do this. 
The $\xi$-deformation of the Iwahori-Matsumoto formulas 
(Theorem \ref{t111}) is the main tool. Its justification
is elementary. It is surprising
that it had not been discovered well before the
double Hecke algebras were introduced.


\end{document}